\documentclass{amsart}

\usepackage{amssymb}

\numberwithin{equation}{section}

\newtheorem{lemma}{Lemma}[section]
\newtheorem{corollary}{Corollary}[section]
\newtheorem{conjecture}{Conjecture}[section]
\newtheorem{theorem}{Theorem}[section]
\newtheorem{proposition}{Proposition}[section]

\theoremstyle{definition}
\newtheorem{remark}{Remark}[section]

\newcommand{\abs}[1]{\lvert#1\rvert}
\newcommand{\qbin}[2]{\genfrac{[}{]}{0pt}{}{#1}{#2}}
\newcommand{\Z}{\mathbb{Z}}
\newcommand{\A}{\mathcal{A}}
\newcommand{\B}{\mathcal{B}}
\newcommand{\C}{\mathcal{C}}
\newcommand{\f}{\tilde{f}}

\begin{document}

\title[Positivity preserving transformations]
{Positivity preserving transformations for $q$-binomial coefficients}

\author{Alexander Berkovich}\thanks{First author supported in part by
NSF grant DMS-0088975}
\address{Department of Mathematics, University of Florida,
Gainesville, FL 32611, USA}
\email{alexb@math.ufl.edu}

\author{S.~Ole Warnaar}\thanks{Second author supported by the Australian 
Research Council}
\address{Department of Mathematics and Statistics,
The University of Melbourne, Vic 3010, Australia}
\email{warnaar@ms.unimelb.edu.au}



\begin{abstract}
Several new transformations for $q$-binomial coefficients are found,
which have the special feature that the kernel is a polynomial with 
nonnegative coefficients. 
By studying the group-like properties of these positivity preserving
transformations, as well as their connection with the Bailey lemma, 
many new summation and transformation formulas for basic hypergeometric 
series are found.
The new $q$-binomial transformations are also applied to obtain
multisum Rogers--Ramanujan identities, 
to find new representations for the Rogers--Szeg\"o polynomials, and
to make some progress on Bressoud's generalized Borwein conjecture.
For the original Borwein conjecture we formulate a refinement based 
on a new triple sum representations of the Borwein polynomials.
\end{abstract}

\maketitle

\section{Introduction}\label{sec1}

\subsection{$q$-Binomial transformations}
In the literature on $q$-series one finds numerous transformations
of the type
\begin{equation}\label{qtrafo}
\sum_{\substack{r=0\\r\equiv j\;(2)}}^L
\frac{q^{\frac{1}{4}r^2}(q;q)_L}{(q;q)_{\frac{1}{2}(L-r)}(q;q)_r}
\qbin{r}{\frac{1}{2}(r-j)}=q^{\frac{1}{4}j^2}\qbin{L}{\frac{1}{2}(L-j)}
\end{equation}
and
\begin{equation}\label{qtrafo2}
\sum_{\substack{r=0\\r\equiv j\;(2)}}^L
\frac{q^{\frac{1}{8}r^2}(q;q)_L}
{(q,-q^{\frac{1}{2}(r+1)};q)_{\frac{1}{2}(L-r)}
(q;q)_r}\qbin{r}{\frac{1}{2}(r-j)}=
q^{\frac{1}{8}j^2}\qbin{L}{\frac{1}{2}(L-j)},
\end{equation}
where $j$ and $L$ are integers such that $j\equiv L\pmod{2}$.
(Throughout this paper the notation $a\equiv b\;(c)$ instead of 
$a\equiv b\pmod{c}$ will be used in equations for brevity).
Here
\begin{equation*}
\qbin{L}{a}_q=\qbin{L}{a}=
\begin{cases}\displaystyle
\frac{(q;q)_L}{(q;q)_a(q;q)_{L-a}}& \text{for $a\in\{0,\dots,L\}$} \\[3mm]
0 &\text{otherwise}
\end{cases}
\end{equation*}
is a $q$-binomial coefficient,
\begin{equation*}
(a;q)_n=\prod_{j=0}^{n-1}(1-aq^j)
\end{equation*}
is a $q$-shifted factorial, and
\begin{equation*}
(a_1,\dots,a_k;q)_n=\prod_{j=1}^k (a_j;q)_n.
\end{equation*}

Important features of \eqref{qtrafo} and \eqref{qtrafo2} are 
(i) the sum over a $q$-binomial coefficient multipied by a simple factor 
again yields a $q$-binomial coefficient, (ii) only the lower entries
of the $q$-binomial coefficients and a simple exponential factor on the right
depend on $j$, (iii) they can readily be iterated.

As an example of this last point let us consider the simple $q$-binomial 
identity
\begin{equation}\label{seed}
\sum_{j=-L}^L (-1)^j q^{\binom{j}{2}}\qbin{2L}{L-j}=\delta_{L,0},
\end{equation}
which is a special case of the finite form of Jacobi's triple product 
identity~\cite[p. 49]{Andrews76} (see also \eqref{FJTP}).
Replacing $L$ by $r$, multiplying both sides by 
$$\frac{q^{r^2}(q;q)_{2L}}{(q;q)_{L-r}(q;q)_{2r}}$$ 
and summing over $r$ using
\eqref{qtrafo} with $L\to 2L$, $j\to 2j$ and $r\to 2r$, yields
\begin{equation}\label{Euler}
\sum_{j=-L}^L (-1)^j q^{j^2+\binom{j}{2}}\qbin{2L}{L-j}=(q^{L+1};q)_L.
\end{equation}
This bounded version of Euler's pentagonal number theorem
\cite{Rogers17} is of the same form as \eqref{seed} and we may
repeat the above procedure to find the well-known bounded analogue of the 
first Rogers--Ramanujan identity \cite{Andrews74,Watson29}
\begin{equation}\label{RR}
\sum_{j=-L}^L (-1)^j q^{2j^2+\binom{j}{2}}\qbin{2L}{L-j}=
(q^{L+1};q)_L \sum_{r=0}^L q^{r^2}\qbin{L}{r}.
\end{equation}

Some known $q$-binomial transformations similar to \eqref{qtrafo} and
\eqref{qtrafo2}, but in which the base of the
$q$-binomial coefficient is changed from $q$ to $q^2$ or $q^3$, are
given by 
\begin{equation}\label{qtrafo3}
\sum_{\substack{r=0\\r\equiv j\;(2)}}^L
\frac{q^{\frac{1}{4}r^2}(q;q)_L}{(q^2;q^2)_{\frac{1}{2}(L-r)}(q;q)_r}
\qbin{r}{\frac{1}{2}(r-j)}=
q^{\frac{1}{4}j^2}\qbin{L}{\frac{1}{2}(L-j)}_{q^2},
\end{equation}
\begin{equation}\label{qtrafo4}
\sum_{\substack{r=0\\r\equiv j\;(2)}}^L
\frac{q^{\frac{1}{8}r^2}(q^{\frac{1}{2}(r+1)};q)_{L-r}(q;q)_L}
{(q^2,q^{r+1};q^2)_{\frac{1}{2}(L-r)}(q;q)_r}
\qbin{r}{\frac{1}{2}(r-j)}=
q^{\frac{1}{8}j^2}\qbin{L}{\frac{1}{2}(L-j)}_{q^2},
\end{equation}
and
\begin{equation}\label{qtrafo5}
\sum_{\substack{r=0\\r\equiv j\;(2)}}^L
\frac{q^{\frac{1}{4}r^2}(q;q)_{\frac{1}{2}(3L-r)}}
{(q^3;q^3)_{\frac{1}{2}(L-r)}(q;q)_r}\qbin{r}{\frac{1}{2}(r-j)}=
q^{\frac{1}{4}j^2}\qbin{L}{\frac{1}{2}(L-j)}_{q^3},
\end{equation}
assuming once again that $j\equiv L\pmod{2}$.

All of the above transformations are of the form
\begin{equation}\label{sym}
\sum_{\substack{r=0\\r\equiv j\;(2)}}^L q^{\frac{1}{8}\gamma r^2} f_{L,r}(q)
\qbin{r}{\frac{1}{2}(r-j)}=
q^{\frac{1}{8}\gamma j^2}\qbin{L}{\frac{1}{2}(L-j)}_{q^k}
\end{equation}
with $f_{L,r}(q)$ a polynomial in $q$ or $q^{1/2}$, which for $0\leq r<L$ has 
both positive and negative coefficients.

The issue of positivity of coefficients in polynomial expressions of the 
type given by the left-hand sides of \eqref{seed}, \eqref{Euler} 
and \eqref{RR} has recently received considerable attention in relation to 
conjectures of Borwein~\cite{Andrews95} and Bressoud~\cite{Bressoud96}.
For this reason it is important to find $q$-binomial transformations
a l\`{a} \eqref{sym} with $f_{L,r}(q)$ a polynomial with nonnegative
coefficients. We will refer to such transformations as positivity preserving.
Indeed, applying a positivity preserving transformation to an
identity like \eqref{seed} --- with on the right a polynomial
with nonnegative coefficients --- results in a new identity which again
has a polynomial with nonnegative coefficients on the right.

\subsection{Outline}
In the next section five new, positivity preserving $q$-binomial 
transformations plus two related, rational transformations are proved. 
In order to establish the positivity of one of our results
we generalize nonnegativity theorems of Andrews for
$q$-binomial coefficients and of Haiman for
principally specialized Schur functions.

Group-like relations among our $q$-binomial transformations
and those listed in the introduction are investigated 
in Section~\ref{sec3}. This will give rise to numerous new transformation
formulas for balanced and `almost' balanced basic hypergeometric series.

The inverses of the transformations for $q$-binomial coefficients
are established in Section~\ref{sec4}. Again this will lead to several 
elegant new summation formulas.

The relation between our $q$-binomial transformations and
the Bailey lemma is the subject of Sections~\ref{sec5} and \ref{sec6}. 
The reader may indeed have recognized \eqref{qtrafo} and \eqref{qtrafo2}
as special cases of the ordinary Bailey lemma in its version due to Andrews 
\cite{Andrews84} and Paule \cite{Paule85}, and 
\eqref{qtrafo3}--\eqref{qtrafo5} as special cases of base-changing 
extensions of the Bailey lemma discovered by Bressoud, Ismail 
and Stanton \cite{BIS00}.
In Section~\ref{sec5} we show that our new transformations
correspond to new types of base-changing Bailey lemmas.
In Section~\ref{sec6} this is exploited to yield some new (and old)
transformations for basic hypergeometric series.

The Sections~\ref{sec7} and \ref{sec8} deal with simple applications of the
$q$-binomial transformations of section~\ref{sec2}. In Section~\ref{sec7}
new single and multisum identities of the Rogers--Ramanujan identities
are proved and in Section~\ref{sec8} we obtain a remarkable new representation
of the Rogers--Szeg\"o polynomials.

Finally, in Section~\ref{sec9}, we use the positivity preserving nature of
our results to make some progress on Bressoud's generalized Borwein 
conjecture. In the last section we also prove new triple-sum representation
for the Borwein polynomials and use this to formulate a new conjecture that
implies the original Borwein conjecture.

\section{Positivity preserving $q$-binomial transformations}\label{sec2}

The reason that none of the transformations of the previous
section preserves positivity is not a very deep one.
Setting $q=1$ in \eqref{sym} yields
\begin{equation*}
\sum_{\substack{r=0\\r\equiv j\;(2)}}^L f_{L,r}(1)
\binom{r}{\frac{1}{2}(r-j)}=\binom{L}{\frac{1}{2}(L-j)},
\end{equation*}
which has the unique solution $f_{L,r}(1)=\delta_{L,r}$.
Hence the only polynomial solution to \eqref{sym} that 
preserves positivity is the less-than-exciting
$f_{L,r}(q)=\delta_{L,r}$ for $k=1$ and $\gamma=0$.
To get around this problem we need to modify \eqref{sym}, and in the
following we look for polynomials $f_{L,r}(q)$ with nonnegative
coefficients that satisfy
\begin{equation}\label{sym2}
\sum_{\substack{r=0\\r\equiv j\;(2)}}^L q^{\frac{1}{4}\gamma r^2} f_{L,r}(q)
\qbin{r}{\frac{1}{2}(r-j)}_{q^k}=
q^{\frac{1}{4}\gamma j^2}\qbin{2L}{L-j}, \qquad k\geq 1
\end{equation}
or small variations hereof (see \eqref{t3} below).

To see that from a positivity point of view \eqref{sym2} is indeed more 
promising than \eqref{sym}, let us again set $q=1$. Multiplying both
sides by $x^{2j}$ and summing over $j$ using the binomial theorem gives
$$\sum_{r=0}^L f_{L,r}(1)(x^2+x^{-2})^r=(x+x^{-1})^{2L}.$$
This is readily solved to yield 
\begin{equation}\label{qis1}
f_{L,r}(1)=2^{L-r}\binom{L}{r},
\end{equation}
a solution that may well have $q$-analogues free of minus signs.

In the remainder of the paper we will make extensive use of 
basic hypergeometric series, and before presenting our solutions to
\eqref{sym2} we need to introduce some further notation \cite{GR90}.
First,
\begin{align*}
{_r\phi_s}\biggl[\genfrac{}{}{0pt}{}
{a_1,\dots,a_r}{b_1,\dots,b_s};q,z\biggr]&=
{_r\phi_s}(a_1,\dots,a_r;b_1,\dots,b_r;q,z) \\
&=\sum_{k=0}^{\infty} \frac{(a_1,\dots,a_r;q)_k}{(q,b_1,\dots,b_s;q)_k}
\Bigl[(-1)^k q^{\binom{k}{2}}\Bigr]^{s-r+1}z^k. \notag
\end{align*}
Here it is assumed that the $b_i$ are such that none of the factors in the
denominator is zero, $q\neq 0$ if $r>s+1$ and $\abs{q}<1$ whenever the
${_r\phi_s}$ is nonterminating. Moreover, if the series does not 
terminate then $r\leq s+1$ with $\abs{z}<1$ if $r=s+1$. If it does
however terminate one can reverse the order of summation
as discussed in \cite[Exercise 1.4]{GR90}.
An ${_{r+1}\phi_r}$ series is called balanced if $z=q$ and 
$a_1\cdots a_{r+1}q=b_1\dots b_r$, well-poised if
$qa_1=a_2b_1=\cdots=a_{r+1}b_r$ and very-well-poised if
it is well-poised and $a_2=-a_3=a_1^{1/2}q$. We will always
abbreviate such very-well-poised series by 
$_{r+1}W_r(a_1;a_4,\dots,a_{r+1};q,z)$.
Whenever one of the numerator parameters in a $q$-hypergeometric series 
is $q^{-n}$ we assume $n$ to be a nonnegative integer.
(Hence, provided the base of the series is $q$ (or $q^{1/2}$, $q^{1/3}$ etc.), 
the series will terminate).
After these definitions we return to \eqref{sym2}.
For $k=1$ it is not hard to see that there are no factorizable
solutions (two non-factorizable or non-$q$-hypergeometric solutions are 
given in Section~\ref{sec9}), and all our results will involve a change of
base. There is of course ample precedent for base-changing transformations
see, e.g., \cite{AV84,AB02,BIS00,GR86,GR90b,GR90,VJ80,VJ82,W02}. 
 
Our first result is of a quadratic nature assuming $k=\gamma=2$.
\begin{lemma}\label{lem1}
For $L$ and $j$ integers there holds
\begin{equation}\label{t2}
\sum_{\substack{r=0\\r\equiv j\;(2)}}^L
q^{\frac{1}{2}r^2}(-q;q)_{L-r}\qbin{L}{r}_{q^2}
\qbin{r}{\frac{1}{2}(r-j)}_{q^2}=q^{\frac{1}{2}j^2}\qbin{2L}{L-j}.
\end{equation}
\end{lemma}
This corresponds to 
\begin{equation}\label{fLr}
f_{L,r}(q)=(-q;q)_{L-r}\qbin{L}{r}_{q^2},
\end{equation}
which is about the simplest imaginable $q$-analogue of \eqref{qis1}. 
Since the $q$-binomial coefficient on the right 
is a polynomial with nonnegative 
coefficients \cite{Andrews76} so is $f_{L,r}(q)$.

By the substitution $q\to 1/q$ and the simple identities
\begin{equation*}
\qbin{L}{a}_{q^{-1}}=q^{-a(L-a)}\qbin{L}{a}_q \quad \text{and}
\quad (a;q^{-1})_n=(-1)^n a^n q^{-\binom{n}{2}}(a^{-1};q)_n
\end{equation*}
we obtain the following corollary of Lemma~\ref{lem1}.
\begin{corollary}
For $L$ and $j$ integers there holds
\begin{equation}\label{t2dual}
\sum_{\substack{r=0\\r\equiv j\;(2)}}^L
q^{\binom{L-r}{2}}(-q;q)_{L-r}\qbin{L}{r}_{q^2}
\qbin{r}{\frac{1}{2}(r-j)}_{q^2}=\qbin{2L}{L-j}.
\end{equation}
\end{corollary}
This corresponds to \eqref{sym2} with $k=2$ and $\gamma=0$.

\begin{proof}[Proof of Lemma \ref{lem1}]
Without loss of generality we may assume that $0\leq j\leq L$.
After shifting $r\to 2r+j$ the identity \eqref{t2} correspond to 
\begin{equation}\label{phi21sum}
{_2\phi_1}(q^{-n},q^{1-n};aq;q^2,aq^{2n})
=\frac{(a;q^2)_n}{(a;q)_n}
\end{equation}
with $(a,n)\to(q^{2j+1},L-j)$.
(Throughout this paper we denote the simultaneous variable changes 
$a_1\to b_1,\dots,a_k\to b_k$ by $(a_1,\dots,a_k)\to (b_1,\dots,b_k)$.)
Equation \eqref{phi21sum} readily follows from the $q$-Gauss sum
\cite[Eq. (II.8)]{GR90}
\begin{equation}\label{qGauss}
{_2\phi_1}(a,b;c;q,c/ab)=\frac{(c/a,c/b;q)_{\infty}}{(c,c/ab;q)_{\infty}}.
\qedhere
\end{equation}
\end{proof}

Our next result is a somewhat more complicated quadratic transformation,
in accordance with \eqref{sym2} for $k=2$ and $\gamma=1$.
\begin{lemma}\label{lem2}
For $L$ and $j$ integers there holds
\begin{equation}\label{t2b}
(1+q^L)\sum_{\substack{r=0\\r\equiv j\;(2)}}^L
q^{\frac{1}{4}r^2}(-q^{r+2};q^2)_{L-r-1}
\qbin{L}{r}\qbin{r}{\frac{1}{2}(r-j)}_{q^2}=
q^{\frac{1}{4}j^2}\qbin{2L}{L-j}.
\end{equation}
\end{lemma}
To make sense of the above lemma we need to extend our earlier definition
of the $q$-shifted factorial, and for nonnegative $n$ we set
$(a;q)_{-n}=1/(aq^{-n};q)_n$. Note that this implies that $1/(q)_{-n}=0$.
With this definition it is once again clear that the corresponding
polynomial $f_{L,r}(q)$ has nonnegative coefficients. 

Before proving \eqref{t2b} we state a variation that is not of the 
form \eqref{sym2}.
\begin{lemma}\label{lem3}
For $L$ and $j$ integers there holds
\begin{equation}\label{t2c}
\sum_{\substack{r=0\\r\equiv j\;(2)}}^L
\frac{q^{\frac{1}{4}r(r+2)}}{1+q^r}(-q^{r+1};q^2)_{L-r}
\qbin{L}{r}\qbin{r}{\frac{1}{2}(r-j)}_{q^2}
=\frac{q^{\frac{1}{4}j(j+2)}}{1+q^j}\qbin{2L}{L-j}.
\end{equation}
\end{lemma}

\begin{proof}[Proof of Lemmas \eqref{lem2} and \eqref{lem3}]
Without loss of generality we may assume that $0\leq j\leq L$.
Shifting $r\to 2r+j$ the summations \eqref{t2b} and \eqref{t2c}
correspond to
\begin{equation}\label{phi32sum}
{_3\phi_2}\biggl[\genfrac{}{}{0pt}{}{a/b,q^{-n},q^{1-n}}
{aq,q^{2-2n}/b};q^2,q^2\biggr]
=\frac{(b;q)_n(a;q^2)_n}{(a;q)_n(b;q^2)_n}
\end{equation}
with $(a,b,n)\to (q^{2j+1},-q^j,L-j)$ and
$(a,b,n)\to (q^{2j+1},-q^{j+1},L-j)$, respectively.
Equation \eqref{phi32sum} follows from the $q\to q^2$ case of the
$q$-Pfaff--Saalsch\"utz sum \cite[Eq. (II.24)]{GR90}
written in the form
\begin{equation}\label{qPS}
{_3\phi_2}\biggl[\genfrac{}{}{0pt}{}{a,b,c}{d,abcq/d};q,q\biggr]
=\frac{(q/d,abq/d,acq/d,bcq/d;q)_{\infty}}
{(aq/d,bq/d,cq/d,abcq/d;q)_{\infty}},
\end{equation}
provided the ${_3\phi_2}$ terminates.
\end{proof}

Our final solution to \eqref{sym2} provides a positivity preserving
transformation of a quartic nature.
\begin{lemma}\label{lem4}
For $L$ and $j$ integers there holds
\begin{equation}\label{t4}
\sum_{\substack{r=0\\r\equiv j\;(2)}}^L
q^{L-r}(-q^{-1};q^2)_{L-r}
\qbin{L}{r}_{q^2}\qbin{r}{\frac{1}{2}(r-j)}_{q^4}=\qbin{2L}{L-j}.
\end{equation}
\end{lemma}
Once again we state a variation that is not of the form \eqref{sym2}.
\begin{lemma}\label{lem5}
For $L$ and $j$ integers there holds
\begin{equation}\label{t4b}
\sum_{\substack{r=0\\r\equiv j\;(2)}}^L
\frac{q^r}{1+q^{2r}}(-q;q^2)_{L-r}
\qbin{L}{r}_{q^2}\qbin{r}{\frac{1}{2}(r-j)}_{q^4}
=\frac{q^j}{1+q^{2j}}\qbin{2L}{L-j}.
\end{equation}
\end{lemma}

\begin{proof}[Proof of Lemma \ref{lem4}]
Without loss of generality we may assume that $0\leq j\leq L$. 
After shifting $r\to 2r+j$ the identity \eqref{t4} corresponds to
\begin{equation}\label{phi43sum}
{_4\phi_3}\biggl[\genfrac{}{}{0pt}{}{aq,aq^3,q^{-2n},q^{2-2n}}
{a^2q^2,-q^{3-2n},-q^{5-2n}};q^4,q^4\biggr]
=q^{-n}\frac{(-q;q)_n(-a;q^2)_n}{(-q^{-1};q^2)_n(-a;q)_n}
\end{equation}
with $(a,n)\to(-q^{2j+1},L-j)$. The above equation follows 
from \cite[Eq. (2.1)]{BIS00} by the substitution
$(C,D,m,q)\to (-q^{1-2n},aq,\lfloor n/2\rfloor,q^2)$.
Unfortunately, the proof of \cite[Eq. (2.1)]{BIS00} as stated in 
\cite{BIS00} appears to be incomplete and below we provide the full 
details of the derivation of \eqref{phi43sum}.

First recall Sears' ${_4\phi_3}$ transformation \cite[Eq. (III.15)]{GR90},
which we write in the form 
\begin{multline}\label{Sears}
{_4\phi_3}\biggl[\genfrac{}{}{0pt}{}{a,b,c,d}
{e,f,abcdq/ef};q,q\biggr] \\
=\frac{(q/f,abq/f,acdq/ef,bcdq/ef;q)_{\infty}}
{(aq/f,bq/f,cdq/ef,abcdq/ef;q)_{\infty}}
{_4\phi_3}\biggl[\genfrac{}{}{0pt}{}{a,b,e/c,e/d}
{e,abq/f,ef/cd};q,q\biggr],
\end{multline}
provided both series terminate.
Letting 
$$(a,b,c,d,e,f,q) \to (q^{-2n},q^{2-2n},aq,aq^3,a^2q^2,
-q^{3-2n},q^4)$$
\eqref{phi43sum} can be written as
\begin{equation}\label{phi43sumb}
{_4\phi_3}\biggl[\genfrac{}{}{0pt}{}{aq^{-1},aq,q^{-2n},q^{2-2n}}
{a^2q^2,-q^{1-2n},-q^{3-2n}};q^4,q^4\biggr]
=\frac{(-q;q)_n(-a;q^2)_n}{(-q;q^2)_n(-a;q)_n}.
\end{equation}
At first sight it may appear that little progress has been made,
but upon closer inspection one may note that the parameters in this new 
${_4\phi_3}$ series are tuned to allow the application 
of Singh's quadratic transformation \cite[Eq. (III.21)]{GR90}
\begin{equation}\label{Singh}
{_4\phi_3}\biggl[\genfrac{}{}{0pt}{}{a^2,b^2,c,d}
{abq^{1/2},-abq^{1/2},-cd};q,q\biggr]=
{_4\phi_3}\biggl[\genfrac{}{}{0pt}{}{a^2,b^2,c^2,d^2}
{a^2b^2q,-cd,-cdq};q^2,q^2\biggr],
\end{equation}
true provided both series terminate. Indeed, utilizing this transform with 
$$(a,b,c,d,q)\to ((a/q)^{1/2},(aq)^{1/2},q^{-n},q^{1-n},q^2)$$
we arrive at \eqref{phi32sum} with $(a,b)\to (-a,-q)$.

Equation \eqref{phi43sumb} may also be derived from the summation
\cite[Eq. (4.3) with $b=1$]{AV84} (rediscovered in \cite[Eq. (2.2)]{BIS00})
by making the substitutions
$(a,b,w,m,q)\to (aq^{-1},1,-aq^{2\lfloor (n+1)/2\rfloor},
\lfloor n/2\rfloor,q^2)$.
\end{proof}

\begin{proof}[Proof of Lemma \ref{lem5}]
Without loss of generality we may assume that $0\leq j\leq L$.
Afer shifting $r\to 2r+j$ the sum \eqref{t4b} correspond to
\eqref{phi43sumb} with $(a,n)\to(-q^{2j+1},L-j)$.
\end{proof}

Our final transformation for $q$-binomial coefficients takes a
form that is slightly different from \eqref{sym2}.
\begin{lemma}\label{lem6}
For $L$ and $j$ integers such that $j\equiv L\pmod{2}$ there holds
\begin{equation}\label{t3}
\sum_{\substack{r=0\\r\equiv j\;(2)}}^{\lfloor L/3\rfloor}
\frac{q^{\frac{3}{4}r^2}(q^3;q^3)_{\frac{1}{2}(L-r-2)}(1-q^L)}
{(q^3;q^3)_r(q;q)_{\frac{1}{2}(L-3r)}}\qbin{r}{\frac{1}{2}(r-j)}_{q^3}
=q^{\frac{3}{4}j^2}\qbin{L}{\frac{1}{2}(L-3j)}.
\end{equation}
\end{lemma}
When $r=L=0$ the factor multiplying the $q$-binomial coefficient
in the summand on the left should be taken to be $1$.

\begin{proof}[Proof of Lemma~\ref{lem6}]
Shifting $r\to 2r+j$ and defining $n=(L-3j)/2$ we arrive at the 
$(a,b,c,d,q)\to (q^{-n},q^{1-n},q^{2-n},q^{3j+3},q^3)$ instance of
\eqref{qPS}.
\end{proof}

Again an important question is whether the polynomial
\begin{equation}\label{fcubic}
f_{L,r}(q)=\frac{(q^3;q^3)_{\frac{1}{2}(L-r-2)}(1-q^L)}
{(q^3;q^3)_r(q;q)_{\frac{1}{2}(L-3r)}}
\end{equation}
for $0\leq 3r\leq L$ and $r\equiv L\pmod{2}$ has nonnegative coefficients.
To answer this is not entirely trivial and we need a generalization of a 
result of Andrews \cite{Andrews97} that arose in connection with 
a monotonicity conjecture of Friedman, Joichi and Stanton \cite{FJS94}.
\begin{theorem}\label{thmA}
Let $k$ and $n$ be positive integers, $j\in\{0,\dots,n\}$ and
$g=\gcd(n,j)$. Then
\begin{equation*}
A_{n,j,k}(q)=\frac{1-q^k}{1-q^n}\qbin{n}{j}
\end{equation*}
is a reciprocal polynomial of degree $j(n-j)+k-n$ with nonnegative 
coefficients if $k\equiv 0\pmod{g}$.
\end{theorem}
For $k=1$ this is Andrews' result \cite[Thm. 2]{Andrews97}.

Assuming the theorem it is not difficult to show that $f_{L,r}(q)$
given by \eqref{fcubic} is a polynomial with nonnegative coefficients.
First we note that for $r=0$ or $3r=L$ this is obvious; $f_{3r,r}(q)=1$ and
$f_{2L,0}(q)=(1+q^L)(q^3;q^3)_{L-1}/(q;q)_{L-1}$ ($L>0$), where
the positivity of the second polynomial follows from $(1-q^3)/(1-q)=1+q+q^2$.
In the following we may therefore assume $0<3r<L$, which implies that
$k:=\gcd((L-r)/2,r)\leq (L-3r)/2$ as follows. There holds $r=uk$ and 
$(L-r)/2=vk$ with $v>u$ and $\gcd(u,v)=1$. Hence
$(L-3r)/2=(v-u)k$ so that $k\leq (L-3r)/2$.
Next we observe the decomposition
\begin{equation*}
f_{L,r}(q)= \frac{(1-q^k)(q^3;q^3)_{(L-3r)/2}}{(1-q^{3k})(q;q)_{(L-3r)/2}}
\times \frac{1-q^L}{1-q^k} \times A_{(L-r)/2,r,k}(q^3),
\end{equation*}
where all three factors on the right are polynomials with nonnegative 
coefficients. The first term because $k\leq (L-3r)/2$ so that
\begin{equation*}
\frac{(1-q^k)(q^3;q^3)_{(L-3r)/2}}{(1-q^{3k})(q;q)_{(L-3r)/2}}
=\sum_{\substack{j=1 \\ j\neq k}}^{(L-3r)/2}(1+q^j+q^{2j}),
\end{equation*}
the second term because $k\mid L$, and the last term because of 
Theorem~\ref{thmA} with $k=g$.

It is possible to arrive at Theorem~\ref{thmA} by modifying Andrews'
proof for $k=1$. Instead, however, we will establish a
more general theorem generalizing results of Haiman \cite[\S 2.5]{Haiman94}
that he used to show polynomiality and nonnegativity of a conjectured 
expression for a specialization of the Frobenius series 
$\mathcal{F}(q,t)$ of `diagonal harmonics'. For most of the terminology
and notation used below we refer to \cite{Macdonald95,Stanley99}.

Let $s_{\lambda}$ be the Schur function labelled by the partition $\lambda$
and define
\begin{equation*}
B_{\lambda,d,k}(q)=\frac{1-q^k}{1-q^d}\: s_{\lambda}(1,q,\dots,q^{d-1}).
\end{equation*}

\begin{theorem}\label{thmB}
Let $d$ and $k$ be positive integers and $\lambda$ a partition such that
$l(\lambda)\leq d$. Set $g=\gcd(d,\abs{\lambda})$.
Then $B_{\lambda,d,k}(q)$ is a reciprocal polynomial of
degree $k-d+\sum_{i=1}^{l(\lambda)}(d-i)\lambda_i$ with nonnegative 
coefficients for every $\lambda$ if $k\equiv 0\pmod{g}$. 
\end{theorem}
For $k=1$ this is due to Haiman.

Before proving the theorem let us show that it includes the previous
theorem as special case. For notational convenience we set 
$q^{\delta}=(q^{d-1},\dots,q,1)$ ($\delta=(d-1,\dots,1,0)$) so that
for $f$ a symmetric function $f(1,q,\dots,q^{d-1})$ may be written as
$f(q^{\delta})$.
Now we choose $\lambda=(j)$ and use that 
\cite[Ch. 1.3, Example 1]{Macdonald95}, \cite[Prop. 2.19.12]{Stanley99}
$$s_{(j)}(q^{\delta})=\qbin{j+d-1}{j}.$$ Therefore
\begin{equation*}
B_{(j),n-j,k}(q)=\frac{1-q^k}{1-q^{n-j}}\qbin{n-1}{j}
=\frac{1-q^k}{1-q^n}\qbin{n}{j}=A_{n,j,k}(q).
\end{equation*}
By Theorem~\ref{thmB} the statement of Theorem~\ref{thmA}
now follows, be it that $j\in\{0,\dots,n-1\}$ and $g=\gcd(j,n-j)$.
Since $\gcd(n-j,j)=\gcd(n,j)$ and since Theorem~\ref{thmA}
is trivially true for $j=n$ this completes our derivation.

\begin{proof}[Proof of Theorem~\ref{thmB}]
Let $\lambda'$ be the conjugate of the partition $\lambda=
(\lambda_1,\dots,\lambda_d)$. Then we have
\cite[Ch. 1.3, Example 1]{Macdonald95}, \cite[Thm. 7.21.2]{Stanley99}
\begin{equation}\label{Schurqd}
s_{\lambda}(q^{\delta})=q^{n(\lambda)}\prod_{x\in \lambda}
\frac{1-q^{d+c(x)}}{1-q^{h(x)}}.
\end{equation}
Here for each $x=(i,j)\in\lambda$ (a partition and its diagram are
identified) the hook-length and content of $x$ are
given by $h(x)=\lambda_i+\lambda_j'-i-j+1$ and $c(x)=j-i$, respectively,
and $n(\lambda)=\sum_{i=1}^d(i-1)\lambda_i$.
To proceed further we need the following lemma, communicated
to us by Richard Stanley.
\begin{lemma}\label{lemcore}
Let $i\mid d$, and let $\omega_i$ be an $i$th primitive
root of unity. Then for $l(\lambda)\leq d$, 
$s_{\lambda}(\omega_i^{\delta})=0$ iff
$\lambda$ has a non-empty $i$-core.
\end{lemma}
To prove this we note that $i\mid d$ and \eqref{Schurqd} imply that
$s_{\lambda}(\omega_i^{\delta})=0$ iff the number of hook-lengths $h(x)$
divisible by $i$ is strictly less than the number of contents $c(x)$
divisible by $i$. Next we recall that the $i$-core of $\lambda$
is obtained from $\lambda$ by repeated removal of border strips of 
length $i$ from the diagram of $\lambda$ until no further strips of 
length $i$ can be removed \cite[Ch. 1.1, Example 8(c)]{Macdonald95},
\cite[Exercise 7.59.d]{Stanley99}.
It is straightforward to verify that each time a border strip is removed,
the number of hook-lengths and the number of contents divisible by $i$
is decreased by one. When we finally reach the $i$-core of $\lambda$ the
number of hook-lengths divisible by $i$ becomes zero. On the other hand,
unless the $i$-core is empty, there will still be a content divisible by
$i$, for example, $c(1,1)=0$. This completes the proof of the lemma.

\begin{remark}
If $i\nmid \abs{\lambda}$ then $\lambda$ has a non-empty $i$-core.
If $i\mid\abs{\lambda}$ and either $\lambda$ or $\lambda'$ consists of a
single row, then $\lambda$ has an empty $i$-core.
However, in general, the $i$-core of $\lambda$ is not necessarily
empty when $i\mid\abs{\lambda}$.
\end{remark}

Next, since $s_{\lambda}(q^{\delta})$ is a polynomial, the
only potential poles of $B_{\lambda,d,k}(q)$ are poles of
$R_{k,d}(q):=(1-q^k)/(1-q^d)$.
Clearly, $R_{k,d}(q)$ has first order poles at each $i$th primitive
root of unity $\omega_i$, provided $i>1$, $i\mid d$, but $i\nmid k$.
Now, if $k\equiv 0\pmod{g}$, then $i\nmid\abs{\lambda}$ and, as a result,
the $i$-core of $\lambda$ is not empty. Hence,
by Lemma~\ref{lemcore}, $s_{\lambda}(\omega_i^{\delta})=0$.
Thus, if $k\equiv 0\pmod{g}$, every pole of $R_{k,d}(q)$ is cancelled by a
zero of $s_{\lambda}(q^{\delta})$, and consequently $B_{\lambda,d,k}(q)$
is polynomial if $k\equiv 0\pmod{g}$.

In the remainder we assume that $k\equiv 0\pmod{g}$.

The degree of $B_{\lambda,d,k}(q)$ immediately follows from the degree
of $s_{\lambda}(q^{\delta})$ given in \cite[Ch. 1.3, Example 1]{Macdonald95}.
To show that the polynomial $B_{\lambda,d,k}(q)$ has nonnegative
coefficients and is reciprocal we use that $s_{\lambda}(q^{\delta})$ is 
a reciprocal, unimodal polynomial with nonnegative coefficients
\cite[Ch. 1.3, Example 1, Ch. 1.8, Example 4]{Macdonald95},
\cite[Exercise 7.75.c]{Stanley99}. This immediately implies
the reciprocality of $B_{\lambda,d,k}(q)$. To see that it also
implies nonnegativity we denote the degree of 
$s_{\lambda}(q^{\delta})$  by $D$ and note that it suffices to show
positivity for $k=g$ thanks to 
$1-q^k=(1-q^g)(1+q^g+\cdots+q^{k-g})$ for $k=mg$.
Now, by the unimodality and nonnegativity of
$s_{\lambda}(q^{\delta})$, it follows that 
$(1-q)s_{\lambda}(q^{\delta})$ is a polynomial of degree $D+1$ with
nonnegative coefficients up to the coefficient of 
$q^{\lfloor (D+1)/2\rfloor}$.
Hence 
\begin{equation*}
B_{\lambda,d,g}(q)=
\frac{1+q+\cdots+q^{g-1}}{1-q^d}(1-q)s_{\lambda}(q^{\delta})
\end{equation*}
is a polynomial of degree $D+g-d$ with nonnegative coefficients up 
to the coefficient of $q^{\lfloor (D+1)/2\rfloor}$. But by its 
reciprocality and by the fact that 
$\lfloor (D+1)/2\rfloor\geq \lfloor (D+g-d)/2\rfloor$ it follows that all
its coefficients must be nonnegative.
\end{proof}
We conclude this section with the following remarks.

\begin{remark}
Lemma~\ref{lemcore} is closely related to 
\cite[Ch. 1.3, Example 17(a)]{Macdonald95}.
It is also a straightforward corollary of \cite[Lem. 2]{SW03}.
\end{remark}
\begin{remark}
It is important to realize that $B_{\lambda,d,k}(q)$
can be a polynomial in $q$ for $k\not\equiv 0\pmod{g}$.
Indeed, the argument given above suggests that
$B_{\lambda,d,k}(q)$ is a polynomial as long as the $i$-core of
$\lambda$ is not empty for any $i$ that divides $d$ but not $k$.
For example, consider the $5$-core partition $\mu=(5,2,2,1)$.
Then
\begin{equation*}
B_{\mu,d,k}(q)=q^9\frac{1-q^5}{1-q}\frac{1-q^7}{1-q}\frac{1-q^k}{1-q}
\frac{1-q^4}{1-q^2}\frac{1-q^9}{1-q^3}
\end{equation*}
is a polynomial for any positive $k$. Note, however,
that when $\lambda=(j)$, $\lambda$ cannot have a non-empty $i$-core if
$i\mid j$, $i>1$. Hence,
$B_{(j),d,k}(q)=A_{d+j,j,k}(q)$ is a polynomial in $q$ iff
$k\equiv 0\pmod{g}$. For $k=1$ this is due to Andrews 
\cite[Thm. 2]{Andrews97}.
\end{remark}

\section{Group-like relations}\label{sec3}
\subsection{Preliminaries}\label{sec31}
Not all of the $q$-binomial transformations of the previous two sections
are independent, and many relations of various degree of
complexity can be found. Such relations are important 
because they often imply new summation or transformation formulas.
For the results of Section~\ref{sec1} the ocurrence of
relations was first investigated by Bressoud \textit{et al.} \cite{BIS00}
and later studied in more detail by Stanton \cite{Stanton01} who introduced 
the notion of the Bailey--Rogers--Ramanujan group. 

For notational reasons we write $q^{\frac{1}{8}\gamma r^2} f_{L,r}(q)$
in \eqref{sym} as $F_{L,r}(q)$ and add as a superscript the relevant
equation number. For example, 
\begin{equation*}
F^{\eqref{qtrafo}}_{L,r}(q)=
\frac{q^{\frac{1}{4}r^2}(q;q)_L}{(q;q)_{(L-r)/2}(q;q)_r}.
\end{equation*}
Likewise we write $q^{\frac{1}{4}\gamma r^2} f_{L,r}(q)$ in \eqref{sym2}
as $F_{L,r}(q)$ and again add equation numbers, and we write
$F^{\eqref{t3}}_{L,r}(q)$ for the kernel of \eqref{t3}. For instance,
\begin{equation*}
F^{\eqref{t2}}_{L,r}(q)=
q^{\frac{1}{2}r^2}(-q;q)_{L-r}\qbin{L}{r}_{q^2}.
\end{equation*}
With this notation we quote from \cite{BIS00,Stanton01}:
\begin{subequations}\label{FCBIS}
\begin{align}
\sum_{\substack{s=r\\s\equiv r\;(2)}}^L
F^{\eqref{qtrafo2}}_{L,s}(q)F^{\eqref{qtrafo2}}_{s,r}(q)&=
F^{\eqref{qtrafo}}_{L,r}(q), \\
\sum_{\substack{s=r\\s\equiv r\;(2)}}^L
F^{\eqref{qtrafo4}}_{L,s}(q)F^{\eqref{qtrafo2}}_{s,r}(q)&=
F^{\eqref{qtrafo3}}_{L,r}(q),
\end{align}
\end{subequations}
and the more complicated
\begin{subequations}\label{FF=FF}
\begin{equation}\label{FF=FF1}
\sum_{\substack{s=r\\s\equiv r\;(2)}}^L
F^{\eqref{qtrafo}}_{L,s}(q)F^{\eqref{qtrafo2}}_{s,r}(q)
=\sum_{\substack{s=r\\s\equiv r\;(2)}}^L
F^{\eqref{qtrafo2}}_{L,s}(q)F^{\eqref{qtrafo}}_{s,r}(q),
\end{equation}
\begin{align}\label{FF=FF2}
\sum_{\substack{s=r\\s\equiv r\;(2)}}^L
F^{\eqref{qtrafo4}}_{L,s}(q)F^{\eqref{qtrafo}}_{s,r}(q)
&=\sum_{\substack{s=r\\s\equiv r\;(2)}}^L
F^{\eqref{qtrafo3}}_{L,s}(q)F^{\eqref{qtrafo2}}_{s,r}(q) \\
&=\sum_{\substack{s=r\\s\equiv r\;(2)}}^L
F^{\eqref{qtrafo2}}_{L,s}(q^2)F^{\eqref{qtrafo4}}_{s,r}(q), \notag
\end{align}
\begin{equation}\label{FF=FF3}
\sum_{\substack{s=r\\s\equiv r\;(2)}}^L
F^{\eqref{qtrafo3}}_{L,s}(q)F^{\eqref{qtrafo}}_{s,r}(q)
=\sum_{\substack{s=r\\s\equiv r\;(2)}}^L
F^{\eqref{qtrafo2}}_{L,s}(q^2)F^{\eqref{qtrafo3}}_{s,r}(q),
\end{equation}
\begin{equation}\label{FF=FF4}
\sum_{\substack{s=r\\s\equiv r\;(2)}}^L
F^{\eqref{qtrafo4}}_{L,s}(q^3)F^{\eqref{qtrafo5}}_{s,r}(q)
=\sum_{\substack{s=r\\s\equiv r\;(2)}}^L
F^{\eqref{qtrafo5}}_{L,s}(q^2)F^{\eqref{qtrafo4}}_{s,r}(q).
\end{equation}
\end{subequations}
(It seems that \eqref{FF=FF1}, \eqref{FF=FF2} and \eqref{FF=FF4} 
are actually missing in \cite{BIS00,Stanton01}.)
The relations in equation \eqref{FCBIS} correspond to summations 
and the relations in \eqref{FF=FF} to transformations for 
basic hypergeometric series.
For example, after shifting $s\to 2s+r$ and replacing 
$(L-r)/2$ by $n$, then using a polynomial argument to
replace $q^{(r+1)/2}$ by the indeterminate $a$, 
and finally using a polynomial argument to replace $q^{-n}$ by $b$, 
\eqref{FF=FF4} becomes the balanced transformation 
\begin{multline}\label{trafo23}
{_5\phi_4}\biggl[\genfrac{}{}{0pt}{}
{a,aq,b^2,b^2\omega,b^2\omega^2}
{a^2,-a^2,-a^2 q,b^6q^2/a^4};q^2,q^2\biggr] \\
=\frac{(a^4/b^6;q^2)_{\infty}(a^3;q^3)_{\infty}(a^6q^3;q^6)_{\infty}}
{(a^4;q^2)_{\infty}(a^3/b^6;q^3)_{\infty}(a^6q^3/b^6;q^6)_{\infty}} \\
\times {_5\phi_4}\biggl[\genfrac{}{}{0pt}{}
{a^2,a^2q,a^2q^2,b^3,-b^3}
{a^3,a^3q^{3/2},-a^3q^{3/2},b^6q^3/a^3};q^3,q^3\biggr],
\end{multline}
provided both series terminate, i.e., provided $a$ or $b$ is of the form
$q^{-n}$. Here $\omega=\exp(2\pi \text{i}/3)$.
For a simple proof of \eqref{FCBIS} and \eqref{FF=FF},
and hence for a proof of the above new transformation we refer
to the next (sub)section.

In the following we extend the analysis of \cite{BIS00,Stanton01} and present
two sets of relations, one of the type $\sum FF=F$ as in
\eqref{FCBIS} and one of the type $\sum FF=\sum FF$ as in \eqref{FF=FF}.
Especially the transformations implied by the second set
are interesting as many appear to be new.

\subsection{Relations of the type $\sum FF=F$}\label{sec32}
Our first set of results, which should be read as five different
ways to decompose $F^{\eqref{t2}}_{L,r}(q)$, is given by
\begin{align}\label{F25}
F^{\eqref{t2}}_{L,r}(q^2)
&=\sum_{\substack{s=r\\s\equiv r\;(2)}}^L
F^{\eqref{t4}}_{L,s}(q^2)F^{\eqref{qtrafo3}}_{s,r}(q^4)
=\sum_{s=r}^L F^{\eqref{qtrafo3}}_{2L,2s}(q)F^{\eqref{t4}}_{s,r}(q) \\
&=\sum_{\substack{s=r\\s\equiv r\;(2)}}^L
F^{\eqref{t2dual}}_{L,s}(q^2)F^{\eqref{qtrafo}}_{s,r}(q^4)
=\sum_{\substack{s=r\\s\equiv r\;(2)}}^L
F^{\eqref{t2b}}_{L,s}(q^2)F^{\eqref{qtrafo2}}_{s,r}(q^4) \notag \\
&=\sum_{s=r}^L F^{\eqref{qtrafo2}}_{2L,2s}(q^2)F^{\eqref{t2dual}}_{s,r}(q^2).
\notag
\end{align}
Similarly, there are three different decompositions of
$F^{\eqref{t2b}}_{L,r}(q)$,
\begin{align}\label{F29}
F^{\eqref{t2b}}_{L,r}(q^2)&=
\sum_{\substack{s=r\\s\equiv r\;(2)}}^L
F^{\eqref{t4}}_{L,s}(q^2)F^{\eqref{qtrafo4}}_{s,r}(q^4)=
\sum_{s=r}^L F^{\eqref{qtrafo4}}_{2L,2s}(q)F^{\eqref{t4}}_{s,r}(q) \\
&=\sum_{\substack{s=r\\s\equiv r\;(2)}}^L
F^{\eqref{t2dual}}_{L,s}(q^2)F^{\eqref{qtrafo2}}_{s,r}(q^4). \notag
\end{align}

\begin{proof}
Since all of the above eight results (and those of Sections~\ref{sec31}
and \ref{sec33}) arise in similar fashion we only show how to prove the 
very first relation.

First take \eqref{qtrafo3} and make the substitution
$(L,q)\to (s,q^2)$. Then multiply this by 
$F^{\eqref{t4}}_{L,s}(q)$ and sum over $s$ to arrive at
\begin{equation*}
\sum_{\substack{s=0\\s\equiv j\;(2)}}^L
\sum_{\substack{r=0\\r\equiv j\;(2)}}^s
F^{\eqref{t4}}_{L,s}(q)F^{\eqref{qtrafo3}}_{s,r}(q^2)
\qbin{r}{\frac{1}{2}(r-j)}_{q^2}
=q^{\frac{1}{2}j^2} \sum_{\substack{s=0\\s\equiv j\;(2)}}^L
F^{\eqref{t4}}_{L,s}(q) \qbin{s}{\frac{1}{2}(s-j)}_{q^4}.
\end{equation*}
Now change the order of summation on the left and
apply \eqref{t4} on the right to get
\begin{equation*}
\sum_{\substack{r=0\\r\equiv j\;(2)}}^L 
\qbin{r}{\frac{1}{2}(r-j)}_{q^2}
\sum_{\substack{s=r\\s\equiv r\;(2)}}^L 
F^{\eqref{t4}}_{L,s}(q) F^{\eqref{qtrafo3}}_{s,r}(q^2)
=q^{\frac{1}{2}j^2}\qbin{2L}{L-j}.
\end{equation*}
Comparing this with \eqref{t2} yields
\begin{equation}\label{FFminF}
\sum_{\substack{r=0\\r\equiv j\;(2)}}^L 
\qbin{r}{\frac{1}{2}(r-j)}_{q^2}
\sum_{\substack{s=r\\s\equiv r\;(2)}}^L 
\Bigl[F^{\eqref{t4}}_{L,s}(q) F^{\eqref{qtrafo3}}_{s,r}(q^2)
-F^{\eqref{t2}}_{L,r}(q)\Bigr]=0,
\end{equation}
which should hold for all integers $L$ and $j$ such that 
$0\leq \abs{j}\leq L$.

The above equation is of the form 
\begin{equation*}
\sum_{\substack{r=0\\r\equiv j\;(2)}}^L
\qbin{r}{\frac{1}{2}(r-j)}h_{L,r}(q)=0,
\end{equation*}
where, without loss of generality, it may be assumed that $0\leq j\leq L$.
Hence the lower bound on the sum may be replaced by $j$.
Recursively it can be seen that $h_{L,r}(q)=0$ is the unique solution.
Indeed, by taking $j=L$ and $j=L-1$ it follows that 
$h_{L,L}(q)=h_{L,L-1}(q)=0$. Next taking $j=L-2$ and $j=L-3$ it in turn 
follows that $h_{L,L-2}(q)=h_{L,L-3}(q)=0$. Repeatedly decreasing $j$
by $2$ it thus follows after $\lfloor L/2\rfloor+1$
steps that all $h_{L,r}(q)$ for $0\leq r\leq L$ must vanish.
Applying this reasoning to \eqref{FFminF} yields the desired
$F^{\eqref{t4}}_{L,s}(q) F^{\eqref{qtrafo3}}_{s,r}(q^2)
-F^{\eqref{t2}}_{L,r}(q)=0$.
\end{proof}

Like \eqref{FCBIS}, the relations of \eqref{F25} and \eqref{F29} (which
should all be read as the left-hand side being equal to one of the
right-hand side expressions) imply summation formulas. 
The only one of these that is possibly new corresponds to
the second relation of \eqref{F29}. After the replacement  
$(q^{r+1/2},L-r)\to (a,n)$ this sum can be stated as
\begin{equation}\label{phi43new}
{_4\phi_3}\biggl[\genfrac{}{}{0pt}{}
{\text{i}q^{-1/2},-\text{i}q^{-1/2},q^{-n},-q^{-n}}
{-q,a,q^{1-2n}/a};q,q^2\biggr]
=\frac{1+a^2q^{2n-1}}{1+a^2q^{-1}}
\frac{(-a^2q^{-1};q^4)_n}{(a;q)_{2n}}.
\end{equation}
Presumably this follows by contiguity from the $b=\text{i}q^{1/2}$ case
of the easily established
\begin{equation*}
{_4\phi_3}\biggl[\genfrac{}{}{0pt}{}
{b,q/b,q^{-n},-q^{-n}}{-q,a,q^{1-2n}/a};q,q\biggr]
=\frac{(ab,aq/b;q^2)_n}{(a;q)_{2n}}
\end{equation*}
or from \eqref{quartic} on page \pageref{pref}
with $bq=-1$. In Section~\ref{sec61} we rederive 
\eqref{phi43new} from the more general transformation 
formula \eqref{BL4BL0}.

\subsection{Relations of the type $\sum FF=\sum FF$}\label{sec33}
This time there are rather a large number of results, all of which
can be proved using the method detailed in Section~\ref{sec32}.
Those relations that imply base-changing transformations from $q$ to $q^k$
for fixed $k$ have been grouped together. 
Here $k$ will be an element of the set $\{1,4/3,3/2,2,4,6,9,12\}$.

\subsubsection{Linear transformations}
There are just two linear relations
\begin{subequations}\label{rel11}
\begin{align}\label{rel11a}
\sum_{s=r}^{\lfloor L/2\rfloor}
F^{\eqref{t2dual}}_{L,2s}(q)F^{\eqref{t2}}_{s,r}(q^2)&=
\sum_{s=r}^{\lfloor L/2\rfloor}
F^{\eqref{t2b}}_{L,2s}(q)F^{\eqref{t2dual}}_{s,r}(q^2), \\
\sum_{s=r}^{\lfloor L/2\rfloor}
F^{\eqref{t2}}_{L,2s}(q)F^{\eqref{t2dual}}_{s,r}(q^2)&=
\sum_{s=r}^{\lfloor L/2\rfloor}
F^{\eqref{t2b}}_{L,2s}(q)F^{\eqref{t2}}_{s,r}(q^2),
\end{align}
\end{subequations}
which are dual in the sense of $q\leftrightarrow 1/q$.
The corresponding $q$-hypergeometric transformations are nothing
but specializations of the identity obtained by equating the
right-hand sides of the Jackson transformations \cite[Eq. (III.4)]{GR90}
and \cite[Eq. (III.5)]{GR90}.

\subsubsection{Transformations from $q$ to $q^{4/3}$}
Much more interesting than \eqref{rel11} are the generalized 
commutation relations
\begin{align*}
\sum_{s=r}^{\lfloor L/3\rfloor}
F^{\eqref{t3}}_{2L,2s}(q)F^{\eqref{t4b}}_{s,r}(q^3)
=\sum_{\substack{s=3r\\s\equiv r\;(2)}}^L
F^{\eqref{t4b}}_{L,s}(q)F^{\eqref{t3}}_{s,r}(q^4),  \\
\sum_{s=r}^{\lfloor L/3\rfloor}
F^{\eqref{t3}}_{2L,2s}(q)F^{\eqref{t4}}_{s,r}(q^3)
=\sum_{\substack{s=3r\\s\equiv r\;(2)}}^L
F^{\eqref{t4}}_{L,s}(q)F^{\eqref{t3}}_{s,r}(q^4).
\end{align*}
Making the variable change $s\to s+r$ on the left and $s\to 2s+3r$ on the
right and then substituting $(q^{6r},L-3r)\to (a,n)$,
the above relations imply the balanced and `almost' balanced formulas
\begin{multline*}
{_5\phi_4}\biggl[\genfrac{}{}{0pt}{}{
\text{i}q^{3/2},-\text{i}q^{3/2},q^{-n},q^{1-n},q^{2-n}}
{-q^3,a^{1/2}q^{3/2},-a^{1/2}q^{3/2},q^{3-3n}/a};q^3,q^3\biggr] \\
=\frac{(-q,a;q^2)_n}{(-q;q)_n(a;q^3)_n}
{_5\phi_4}\biggl[\genfrac{}{}{0pt}{}{
a^{2/3},a^{2/3}\omega,a^{2/3}\omega^2,q^{-2n},q^{2-2n}}
{a,aq^2,-q^{1-2n},-q^{3-2n}};q^4,q^4\biggr]
\end{multline*}
and
\begin{multline*}
{_5\phi_4}\biggl[\genfrac{}{}{0pt}{}{
\text{i}q^{-3/2},-\text{i}q^{-3/2},q^{-n},q^{1-n},q^{2-n}}
{-q^3,a^{1/2}q^{3/2},-a^{1/2}q^{3/2},q^{3-3n}/a};q^3,q^6\biggr] \\
=q^n\frac{(-q^{-1},a;q^2)_n}{(-q;q)_n(a;q^3)_n}
{_6\phi_5}\biggl[\genfrac{}{}{0pt}{}{
a^{2/3},a^{2/3}\omega,a^{2/3}\omega^2,-a q^4,q^{-2n},q^{2-2n}}
{a,-a,aq^2,-q^{3-2n},-q^{5-2n}};q^4,q^4\biggr],
\end{multline*}
respectively. To the best of our knowledge these are the
first examples of a transformations relating base $q^3$ and $q^4$.

\subsubsection{Transformations from $q$ to $q^{3/2}$}
Again there are two results, not dissimilar to the previous pair;
\begin{align*}
\sum_{s=r}^{\lfloor L/3\rfloor}
F^{\eqref{t3}}_{2L,2s}(q)F^{\eqref{t2b}}_{s,r}(q^3)&=
\sum_{\substack{s=3r\\s\equiv r\;(2)}}^L
F^{\eqref{t2b}}_{L,s}(q)F^{\eqref{t3}}_{s,r}(q^2), \\
\sum_{s=r}^{\lfloor L/3\rfloor}
F^{\eqref{t3}}_{2L,2s}(q)F^{\eqref{t2c}}_{s,r}(q^3)&=
\sum_{\substack{s=3r\\s\equiv r\;(2)}}^L
F^{\eqref{t2c}}_{L,s}(q)F^{\eqref{t3}}_{s,r}(q^2).
\end{align*}
Making the same variable change as above and then substituting
$(-q^{3r},L-3r)\to (a,n)$, this yields
\begin{multline*}
{_5\phi_4}\biggl[\genfrac{}{}{0pt}{}
{a^{2/3},a^{2/3}\omega,a^{2/3}\omega^2,q^{-n},q^{1-n}}
{a,-a,-aq,q^{2-2n}/a};q^2,q^2\biggr] \\
=\frac{(a^2;q^3)_n}{(a;q^2)_n(-a;q)_n}
{_5\phi_4}\biggl[\genfrac{}{}{0pt}{}
{a^{1/2},-a^{1/2},q^{-n},q^{1-n},q^{2-n}}
{a,aq^{3/2},-aq^{3/2},q^{3-3n}/a^2};q^3,q^3\biggr]
\end{multline*}
and 
\begin{multline*}
{_5\phi_4}\biggl[\genfrac{}{}{0pt}{}
{a^{2/3},a^{2/3}\omega,a^{2/3}\omega^2,q^{-n},q^{1-n}}
{aq,-aq,-a,q^{1-2n}/a};q^2,q^2\biggr] \\
=\frac{1-a^2q^{2n}}{1-a^2}
\frac{(a^2;q^3)_n}{(aq;q^2)_n(-aq;q)_n} 
\\ \times
{_5\phi_4}\biggl[\genfrac{}{}{0pt}{}
{a^{1/2}q^{3/2},-a^{1/2}q^{3/2},q^{-n},q^{1-n},q^{2-n}}
{aq^3,aq^{3/2},-aq^{3/2},q^{3-3n}/a^2};q^3,q^3\biggr].
\end{multline*}
Both these results should be compared with \eqref{trafo23}.

\subsubsection{Quadratic transformations}
There are quite a number of different relations of a quadratic nature.
First, 
\begin{subequations}
\begin{align}\label{rel12a}
\sum_{s=r}^L F^{\eqref{qtrafo}}_{2L,2s}(q)F^{\eqref{t2dual}}_{s,r}(q)&=
\sum_{s=r}^L F^{\eqref{qtrafo2}}_{2L,2s}(q)F^{\eqref{t2}}_{s,r}(q) \\
&=\sum_{\substack{s=r\\s\equiv r\;(2)}}^L
F^{\eqref{t2}}_{L,s}(q)F^{\eqref{qtrafo}}_{s,r}(q^2), \notag
\end{align}
\begin{align}\label{rel12b}
\sum_{s=r}^L F^{\eqref{qtrafo2}}_{2L,2s}(q)F^{\eqref{t2b}}_{s,r}(q)&=
\sum_{\substack{s=r\\s\equiv r\;(2)}}^L
F^{\eqref{t2b}}_{L,s}(q)F^{\eqref{qtrafo}}_{s,r}(q^2) \\
&=\sum_{\substack{s=r\\s\equiv r\;(2)}}^L
F^{\eqref{t2}}_{L,s}(q)F^{\eqref{qtrafo2}}_{s,r}(q^2). \notag
\end{align}
\end{subequations}
The first equality in \eqref{rel12a} corresponds to a specialization 
of the transformation \cite[Eq. (III.4)]{GR90}, and the second equality 
implies the $(a,b,c,n)\to (q^{r+1/2},\infty,0,L-r)$ specialization of
\begin{multline}\label{trafo1}
{_4\phi_3}\biggl[\genfrac{}{}{0pt}{}{b,c,-c,q^{-n}}
{a,c^2,-bq^{1-n}/a};q,q\biggr] \\
=\frac{(a^2/b;q)_n(c^2;q^2)_n}{(a,c^2,-a/b;q)_n}
{_4\phi_3}\biggl[\genfrac{}{}{0pt}{}{a^2/b^2,a^2/c^2,q^{-n},q^{1-n}}
{a^2/b,a^2 q/b,q^{2-2n}/c^2};q^2,q^2\biggr],
\end{multline}
proved in Section~\ref{sec62}.
Similarly, the second equality in \eqref{rel12b} corresponds to a 
specialization of the transformation \cite[Eq. (III.12)]{GR90}, 
and the first equality implies the 
$(a,b,c,n)\to (q^{r+1/2},\infty,\text{i}q^{r/2},L-r)$ 
specialization of \eqref{trafo1}.

Next are the four closely related results
\begin{subequations}
\begin{equation}\label{rel12c}
\sum_{s=r}^{\lfloor L/2\rfloor}
F^{\eqref{t2}}_{L,2s}(q)F^{\eqref{t4b}}_{s,r}(q^2)=
\sum_{s=r}^{\lfloor L/2\rfloor}
F^{\eqref{t4b}}_{L,2s}(q)F^{\eqref{t2}}_{s,r}(q^4),
\end{equation}
\begin{equation}
\sum_{s=r}^{\lfloor L/2\rfloor}
F^{\eqref{t2dual}}_{L,2s}(q)F^{\eqref{t4b}}_{s,r}(q^2)=
\sum_{s=r}^{\lfloor L/2\rfloor}
F^{\eqref{t4b}}_{L,2s}(q)F^{\eqref{t2dual}}_{s,r}(q^4),
\end{equation}
\begin{equation}\label{rel12e}
\sum_{s=r}^{\lfloor L/2\rfloor}
F^{\eqref{t2}}_{L,2s}(q)F^{\eqref{t4}}_{s,r}(q^2)=
\sum_{s=r}^{\lfloor L/2\rfloor}
F^{\eqref{t4}}_{L,2s}(q)F^{\eqref{t2}}_{s,r}(q^4),
\end{equation}
\begin{equation}
\sum_{s=r}^{\lfloor L/2\rfloor}
F^{\eqref{t2dual}}_{L,2s}(q)F^{\eqref{t4}}_{s,r}(q^2)=
\sum_{s=r}^{\lfloor L/2\rfloor}
F^{\eqref{t4}}_{L,2s}(q)F^{\eqref{t2dual}}_{s,r}(q^4),
\end{equation}
\end{subequations}
The first as well as the last two relations
are dual in the sense of $q\leftrightarrow 1/q$.
After the substitution $(q^{4r+2},L-2r)\to (a,n)$
equation \eqref{rel12c} implies 
\begin{equation*}
{_4\phi_3}\biggl[\genfrac{}{}{0pt}{}{\text{i}q,-\text{i}q,
q^{-n},q^{1-n}}{-q^2,a^{1/2},-a^{1/2}};q^2,aq^{2n-1}\biggr]
=\frac{(-q;q^2)_n}{(-q;q)_n}
{_4\phi_3}\biggl[\genfrac{}{}{0pt}{}{0,-aq^{-2},q^{-2n},q^{2-2n}}
{a,-q^{1-2n},-q^{3-2n}};q^4,q^4\biggr]
\end{equation*}
and equation \eqref{rel12e} implies
\begin{multline*}
{_4\phi_3}\biggl[\genfrac{}{}{0pt}{}{\text{i}q^{-1},-\text{i}q^{-1},
q^{-n},q^{1-n}}{-q^2,a^{1/2},-a^{1/2}};q^2,aq^{2n+1}\biggr]\\
=q^n\frac{(-q^{-1};q^2)_n}{(-q;q)_n}
{_4\phi_3}\biggl[\genfrac{}{}{0pt}{}{0,-aq^2,q^{-2n},q^{2-2n}}
{a,-q^{3-2n},-q^{5-2n}};q^4,q^4\biggr].
\end{multline*}

Finally there holds
\begin{subequations}
\begin{align}
\sum_{s=r}^{\lfloor L/2\rfloor} F^{\eqref{t2b}}_{L,2s}(q)
F^{\eqref{t4b}}_{s,r}(q^2) &=
\sum_{s=r}^{\lfloor L/2\rfloor} F^{\eqref{t4b}}_{L,2s}(q)
F^{\eqref{t2b}}_{s,r}(q^4), \\
\sum_{s=r}^{\lfloor L/2\rfloor} F^{\eqref{t2b}}_{L,2s}(q)
F^{\eqref{t4}}_{s,r}(q^2) &=
\sum_{s=r}^{\lfloor L/2\rfloor} F^{\eqref{t4}}_{L,2s}(q)
F^{\eqref{t2b}}_{s,r}(q^4).
\end{align}
\end{subequations}
After the replacement $(q^{2r},L-2r)\to (a,n)$ these yield
\begin{multline}\label{tr4p}
{_4\phi_3}\biggl[\genfrac{}{}{0pt}{}{\text{i}q,-\text{i}q,
q^{-n},q^{1-n}}{-q^2,aq,-q^{2-2n}/a};q^2,q^2\biggr] \\
=\frac{(-a;q)_n(-q;q^2)_n}{(-q;q)_n(-a;q^2)_n}
{_4\phi_3}\biggl[\genfrac{}{}{0pt}{}{
\text{i}a,-\text{i}a,q^{-2n},q^{2-2n}}
{a^2q^2,-q^{1-2n},-q^{3-2n}};q^4,q^4\biggr]
\end{multline}
and 
\begin{multline}\label{tr4}
{_4\phi_3}\biggl[\genfrac{}{}{0pt}{}{\text{i}q^{-1},-\text{i}q^{-1},
q^{-n},q^{1-n}}{-q^2,aq,-q^{2-2n}/a};q^2,q^4\biggr] \\
=q^n\frac{(-q^{-1};q^2)_n(-a;q)_n}{(-q;q)_n(-a;q^2)_n}
{_5\phi_4}\biggl[\genfrac{}{}{0pt}{}{
\text{i}a,-\text{i}a,-a^2 q^4,q^{-2n},q^{2-2n}}
{-a^2,a^2q^2,-q^{3-2n},-q^{5-2n}};q^4,q^4\biggr],
\end{multline}
respectively.
It is not hard to see that \eqref{tr4p} is a special case of
\begin{multline}\label{comp}
{_4\phi_3}\biggl[\genfrac{}{}{0pt}{}{b,q^2/b,
q^{-n},q^{1-n}}{-q^2,aq,-q^{2-2n}/a};q^2,q^2\biggr] \\
=\frac{(-a;q)_n(-q;q^2)_n}{(-q;q)_n(-a;q^2)_n}
{_4\phi_3}\biggl[\genfrac{}{}{0pt}{}{
aq/b,ab/q,q^{-2n},q^{2-2n}}
{a^2q^2,-q^{1-2n},-q^{3-2n}};q^4,q^4\biggr],
\end{multline}
which generalizes \eqref{phi43sumb} and follows by first applying
Singh's quadratic transformation \eqref{Singh}
to the right-side and then using Sears' ${_4\phi_3}$
transformation \eqref{Sears}. (Equation \eqref{comp} also follows
from \cite[Eq. (4.3)]{AV84} by a single use of Sears' transform.)
Because of the ${_5\phi_4}$ series on the right, it is unclear whether 
\eqref{tr4} admits a similar kind of generalization.

\subsubsection{Quartic transformations}
Our list of quartic relations begins with
\begin{subequations}\label{qu1}
\begin{align}\label{rel3c}
\sum_{s=r}^L F^{\eqref{qtrafo3}}_{2L,2s}(q)F^{\eqref{t2dual}}_{s,r}(q)&=
\sum_{s=r}^L F^{\eqref{qtrafo4}}_{2L,2s}(q)F^{\eqref{t2}}_{s,r}(q) \\
&=\sum_{\substack{s=r\\s\equiv r\;(2)}}^L
F^{\eqref{t2b}}_{L,s}(q^2)F^{\eqref{qtrafo3}}_{s,r}(q^2), \notag 
\end{align}
\begin{equation}\label{com2717}
\sum_{s=r}^L F^{\eqref{qtrafo4}}_{2L,2s}(q)F^{\eqref{t2b}}_{s,r}(q)=
\sum_{\substack{s=r\\s\equiv r\;(2)}}^L
F^{\eqref{t2b}}_{L,s}(q^2)F^{\eqref{qtrafo4}}_{s,r}(q^2).
\end{equation}
\end{subequations}
The first equality in \eqref{rel3c} once again corresponds to a 
specialization of \cite[Eq. (III.12)]{GR90}. More interesting
are the second equality in \eqref{rel3c} and the generalized
commutation relation \eqref{com2717}.
These prove the $(a,b,n)\to(q^{r+1/2},0,L-r)$, respectively,
$(a,b,n)\to(q^{r+1/2},-q^r,L-r)$ case of \label{pref}
\begin{multline}\label{quartic}
{_4\phi_3}\biggl[\genfrac{}{}{0pt}{}
{b^{1/2},-b^{1/2},q^{-n},-q^{-n}}
{a,b,q^{1-2n}/a};q,q\biggr] \\
=\frac{1+a^2q^{2n-1}}{1+a^2q^{-1}}
\frac{(-a^2q^{-1};q^4)_n}{(a;q)_{2n}}
{_4\phi_3}\biggl[\genfrac{}{}{0pt}{}{-bq,-bq^3,q^{-2n},q^{2-2n}}
{-a^2q^3,b^2q^2,-q^{5-4n}/a^2};q^4,q^4\biggr],
\end{multline}
established in Section~\ref{sec62}.
As a variation on the above there also holds
\begin{subequations}\label{qu2}
\begin{align}\label{com1623}
\sum_{s=r}^L F^{\eqref{qtrafo3}}_{2L,2s}(q)F^{\eqref{t2}}_{s,r}(q)&=
\sum_{\substack{s=r\\s\equiv r\;(2)}}^L
F^{\eqref{t2}}_{L,s}(q^2)F^{\eqref{qtrafo3}}_{s,r}(q^2), \\
\label{rel6737}
\sum_{s=r}^L F^{\eqref{qtrafo3}}_{2L,2s}(q)F^{\eqref{t2b}}_{s,r}(q)&=
\sum_{\substack{s=r\\s\equiv r\;(2)}}^L
F^{\eqref{t2}}_{L,s}(q^2)F^{\eqref{qtrafo4}}_{s,r}(q^2).
\end{align}
\end{subequations}
Here \eqref{com1623}, respectively, \eqref{rel6737}
imply the $(a,b,n)\to(q^{r+1/2},0,L-r)$ and
$(a,b,n)\to(q^{r+1/2},-q^r,L-r)$ instances of
\begin{multline}\label{quartic2}
{_4\phi_3}\biggl[\genfrac{}{}{0pt}{}
{b^{1/2},-b^{1/2},q^{-n},-q^{-n}}
{a,-a,b};q,-a^2q^{2n}\biggr] \\
=\frac{(-a^2q;q^2)_n}{(a^2;q^2)_n}
{_4\phi_3}\biggl[\genfrac{}{}{0pt}{}{-bq,-bq^3,q^{-2n},q^{2-2n}}
{-a^2q,-a^2q^3,b^2q^2};q^4,a^4q^{4n}\biggr].
\end{multline}
Once again this is proved in Section~\ref{sec62}.
By making the substitution $(q^{r+1/2},L-r)\to (a,n)$ the relation
\begin{equation}\label{rel7556}
\sum_{s=r}^L F^{\eqref{qtrafo4}}_{2L,2s}(q)F^{\eqref{t2dual}}_{s,r}(q)=
\sum_{\substack{s=r\\s\equiv r\;(2)}}^L
F^{\eqref{t2dual}}_{L,s}(q^2)F^{\eqref{qtrafo3}}_{s,r}(q^2),
\end{equation}
yields the $b\to\infty$ limit of the quartic transformation
\begin{multline}\label{eqfail}
{_4\phi_3}\biggl[\genfrac{}{}{0pt}{}{b^{1/2},-b^{1/2},q^{-n},-q^{-n}}
{a,b,q^{1-2n}/a};q,q\biggr] \\
=\frac{(-a^2q;q^2)_n(b^2;q^4)_n}{(a;q)_{2n}(b^2;q^2)_n}
{_4\phi_3}\biggl[\genfrac{}{}{0pt}{}{a^2/b,a^2q^2/b,q^{-2n},q^{2-2n}}
{-a^2q,-a^2q^3,q^{4-4n}/b^2};q^4,q^4\biggr].
\end{multline}
It is not hard to prove this identity by applying Sears' ${_4\phi_3}$
transformation \eqref{Sears} with
$$(a,b,c,d,e,f,q)\to (q^{-2n},q^{2-2n},-bq,-bq^3,
-a^2q^3,b^2q^2,q^4)$$
to the right-hand side of \eqref{quartic}.

Next is the pair
\begin{align*}
\sum_{s=r}^L F^{\eqref{qtrafo2}}_{2L,2s}(q)F^{\eqref{t4b}}_{s,r}(q) 
&=\sum_{\substack{s=r\\s\equiv r\;(2)}}^L
F^{\eqref{t4b}}_{L,s}(q)F^{\eqref{qtrafo2}}_{s,r}(q^4), \\
\sum_{s=r}^L F^{\eqref{qtrafo2}}_{2L,2s}(q)F^{\eqref{t4}}_{s,r}(q)
&=\sum_{\substack{s=r\\s\equiv r\;(2)}}^L
F^{\eqref{t4}}_{L,s}(q)F^{\eqref{qtrafo2}}_{s,r}(q^4).
\end{align*}
After the substitutions $(q^{r+1/2},L-r)\to (a,n)$ these lead to
\begin{equation}\label{fail1}
{_3\phi_2}\biggl[\genfrac{}{}{0pt}{}
{\text{i}q^{1/2},-\text{i}q^{1/2},q^{-n}}
{-q,a};q,-aq^n\biggr]
=\frac{(-q;q^2)_n}{(-q,a;q)_n}
{_3\phi_2}\biggl[\genfrac{}{}{0pt}{}{-a^2q^{-1},q^{-2n},q^{2-2n}}
{-q^{1-2n},-q^{3-2n}};q^4,q^4\biggr]
\end{equation}
and
\begin{multline}\label{fail2}
{_3\phi_2}\biggl[\genfrac{}{}{0pt}{}
{\text{i}q^{-1/2},-\text{i}q^{-1/2},q^{-n}}
{-q,a};q,-aq^{n+1}\biggr] \\
=q^n\frac{(-q^{-1};q^2)_n}{(-q,a;q)_n}
{_3\phi_2}\biggl[\genfrac{}{}{0pt}{}{-a^2q^3,q^{-2n},q^{2-2n}}
{-q^{3-2n},-q^{5-2n}};q^4,q^4\biggr],
\end{multline}
which we failed to generalize to the level of ${_4\phi_3}$ (or
${_5\phi_4}$) series. It is however not hard to see that
by applying Singh's transformation
\eqref{Singh} to the right-hand side, \eqref{fail1}
becomes the $(b,c)\to (\infty,\text{i}q^{1/2})$ limit of
\eqref{trafo1}. It is also possible to arrive at
\eqref{fail1} and \eqref{fail2} (with $A$ replaced by $a$)
by taking the $a,b\to\infty$ limit
in \eqref{tr1p} and \eqref{tr1} such that $A=-bq^{1-n}/a$ is fixed,
and by then transforming the resulting ${_3\phi_2}$ series on the right
using \cite[Eq. (III.13)]{GR90}.

Our last two quartic commutation relations are rather interesting,
\begin{subequations}
\begin{align}\label{rel14a}
\sum_{s=r}^L F^{\eqref{qtrafo}}_{2L,2s}(q)F^{\eqref{t4b}}_{s,r}(q) 
&=\sum_{\substack{s=r\\s\equiv r\;(2)}}^L
F^{\eqref{t4b}}_{L,s}(q)F^{\eqref{qtrafo}}_{s,r}(q^4), \\
\label{rel14b}
\sum_{s=r}^L F^{\eqref{qtrafo}}_{2L,2s}(q)F^{\eqref{t4}}_{s,r}(q)
&=\sum_{\substack{s=r\\s\equiv r\;(2)}}^L
F^{\eqref{t4}}_{L,s}(q)F^{\eqref{qtrafo}}_{s,r}(q^4).
\end{align}
\end{subequations}
Equation \eqref{rel14a} implies the $(a,b,n)\to (q^{2r+1},\infty,L-r)$ 
instance of
\begin{multline}\label{tr1p}
{_5\phi_4}\biggl[\genfrac{}{}{0pt}{}{\text{i}q^{1/2},-\text{i}q^{1/2},
b^{1/2},-b^{1/2},q^{-n}}{-q,a^{1/2},-a^{1/2},-bq^{1-n}/a};q,q\biggr] \\
=\frac{(-q,a^2/b;q^2)_n}{(-q,-a/b;q)_n(a;q^2)_n}
{_5\phi_4}\biggl[\genfrac{}{}{0pt}{}{a^2/b^2,-aq^{-1},-aq,q^{-2n},q^{2-2n}}
{a^2/b,a^2q^2/b,-q^{1-2n},-q^{3-2n}};q^4,q^4\biggr].
\end{multline}
This result, which will be proved in Section~\ref{sec62},
simplifies to \eqref{phi43sumb} for $b=1$.
Similarly, \eqref{rel14b} corresponds to the
$(a,b,n)\to (q^{2r+1},\infty,L-r)$ case of
\begin{multline}\label{tr1}
{_5\phi_4}\biggl[\genfrac{}{}{0pt}{}{\text{i}q^{-1/2},-\text{i}q^{-1/2},
b^{1/2},-b^{1/2},q^{-n}}{-q,a^{1/2},-a^{1/2},-bq^{1-n}/a};q,q^2\biggr] \\
=q^n\frac{(-q^{-1},a^2/b;q^2)_n}{(-q,-a/b;q)_n(a;q^2)_n}
{_5\phi_4}\biggl[\genfrac{}{}{0pt}{}{a^2/b^2,-aq,-aq^3,q^{-2n},q^{2-2n}}
{a^2/b,a^2q^2/b,-q^{3-2n},-q^{5-2n}};q^4,q^4\biggr].
\end{multline}
When $b=1$ this simplifies to \eqref{phi43sum} and when $aq=-1$
(and $b\to b^2)$ to
\begin{equation}\label{p32}
{_3\phi_2}(b,-b,q^{-n};-q,b^2q^{2-n};q,q^2)=
q^n\frac{(q^{-2}/b^2;q^2)_n}{(-q,q^{-1}/b^2;q)_n}
\end{equation}
needed shortly. The proof of \eqref{tr1} can again be found
in Section~\ref{sec62}.

Both \eqref{tr1p} and \eqref{tr1} may be further manipulated into new
quadratic transformations as follows.
The left-hand side of \eqref{tr1p} simplifies to a ${_3\phi_2}$ series
by the $(b,x,y)\to (\text{i}(bq/a)^{1/2},(a/b)^{1/2},
\text{i}(a/q)^{1/2})$ case of \cite[Eq (3.5.2); $a\to q^{-n}$]{GR90}
\begin{equation}\label{eq352}
{_5\phi_4}\biggl[\genfrac{}{}{0pt}{}{bx,-bx,by,-by,q^{-n}}
{-q,bxy,-bxy,b^2q^{-n}};q,q\biggr]
=\frac{(q^2/b^2;q^2)_n}{(-q,q/b^2;q)_n}
{_3\phi_2}\biggl[\genfrac{}{}{0pt}{}{x^2,y^2,q^{-2n}}
{b^2x^2y^2,b^2q^{-2n}};q^2,b^2q^3\biggr].
\end{equation}
After the further substitution $(a,b,q)\to (aq,aq/b,-q)$ this leads to
\begin{equation}\label{fin54}
{_5\phi_4}\biggl[\genfrac{}{}{0pt}{}{a,aq^2,b^2,q^{-2n},q^{2-2n}}
{abq,abq^3,q^{1-2n},q^{3-2n}};q^4,q^4\biggr]
=\frac{(aq,bq;q^2)_n}{(q,abq;q^2)_n}
{_3\phi_2}\biggl[\genfrac{}{}{0pt}{}{a,b,q^{-2n}}
{aq,q^{1-2n}/b};q^2,-\frac{q^2}{b}\biggr].
\end{equation}
When $(b,q)\to (-aq,-q)$ the left is summable by 
\eqref{phi43sumb} and we infer the further identity
\begin{equation}\label{phi18}
{_2\phi_1}(a,q^{-2n};q^{-2n}/a;q^2,q/a)=
\frac{(-q,aq;q)_n}{(aq^2;q^2)_n},
\end{equation}
which also follows from \cite[Exercise 1.8]{GR90}.
By the usual polynomial argument equation \eqref{fin54} may
also be stated as
\begin{equation*}
{_5\phi_4}\biggl[\genfrac{}{}{0pt}{}{a,aq^2,b^2,c,cq^2}
{abq,abq^3,cq,cq^3};q^4,q^4\biggr]
=\frac{(q,cq/a,cq/b,q/ab;q^2)_{\infty}}{(cq,q/a,q/b,cq/ab;q^2)_{\infty}}
{_3\phi_2}\biggl[\genfrac{}{}{0pt}{}{a,b,c}
{aq,cq/b};q^2,-\frac{q^2}{b}\biggr],
\end{equation*}
provided both series terminate. For $c=aq$ the ${_3\phi_2}$
series on the right becomes a ${_2\phi_1}$ which
precisely takes the form of the sum side of the
Bailey-Daum summation \cite[Eq. (II.9)]{GR90}.
\begin{remark}
When $a=q^{2j}$ with $j\geq 1$, 
equation \eqref{phi18} may be put in the form
\begin{equation}\label{qpartid}
\sum_{k=0}^n q^k\qbin{k+j-1}{k}_{q^2}\qbin{n-k+j}{j}_{q^2}=
\qbin{n+2j}{n}.
\end{equation}
This has the following elegant partition theoretic interpretation.
The expression
\begin{equation*}
q^k\qbin{k+j-1}{k}_{q^2}
\end{equation*}
is the generating function of partitions of exactly $k$ parts, with all parts
being odd and no parts exceeding $2j-1$.
The expression 
\begin{equation*}
\qbin{n-k+j}{j}_{q^2}
\end{equation*}
is the generating function of partitions of at most $n-k$ parts, with all
parts being even and no parts exceeding $2j$.
Hence the summand on the left of \eqref{qpartid} is the generating function
of partitions of at most $n$ parts, with no parts exceeding $2j$ and exactly
$k$ odd parts. When summed over the number of odd parts this gives
the generating function of partitions of at most $n$ parts with
no parts exceeding $2j$, in accordance with the right-hand side 
of \eqref{qpartid}.
\end{remark}

To also rewrite \eqref{tr1} as a quadratic transformation
requires a bit more work.
Indeed, in order to trade the ${_5\phi_4}$ on the left for a ${_3\phi_2}$ 
we need to prove the following companion to \eqref{eq352}:
\begin{multline}\label{var352}
{_5\phi_4}\biggl[\genfrac{}{}{0pt}{}{bx,-bx,by,-by,q^{-n}}
{-q,bxy,-bxy,b^2q^{2-n}};q,q^2\biggr] \\
=q^n \frac{(q^{-2}/b^2;q^2)_n}{(-q,q^{-1}/b^2;q)_n}
{_3\phi_2}\biggl[\genfrac{}{}{0pt}{}{x^2,y^2,q^{-2n}}
{b^2x^2y^2,b^2q^{4-2n}};q^2,b^2q^3\biggr].
\end{multline}
Using this with $(b,x,y)\to (\text{i}(b/aq)^{1/2},(a/b)^{1/2},
\text{i}(aq)^{1/2})$ and making the further substitution 
$(a,b,q)\to (aq^{-1},aq^{-1}/b,-q)$ yields 
\begin{multline}\label{phi54phi32}
{_5\phi_4}\biggl[\genfrac{}{}{0pt}{}{a,aq^2,b^2,q^{-2n},q^{2-2n}}
{abq^{-1},abq,q^{3-2n},q^{5-2n}};q^4,q^4\biggr] \\
=\frac{(aq^{-1},bq^{-1};q^2)_n}{(q^{-1},abq^{-1};q^2)_n}
{_3\phi_2}\biggl[\genfrac{}{}{0pt}{}{a,b,q^{-2n}}
{aq^{-1},q^{3-2n}/b};q^2,-\frac{q^2}{b}\biggr].
\end{multline}
When $(a,b,q)\to (aq,-a,-q)$ the sum on the left can be carried out by 
\eqref{phi43sum} leading to
\begin{equation}\label{phi21sumb}
{_2\phi_1}(aq,q^{-2n};q^{3-2n}/a;q^2,q^2/a)=
q^{-n}\frac{(-q,a;q)_n}{(a/q;q^2)_n}.
\end{equation}
This sum, which is in fact \eqref{phi18} with order of summation reversed,
will be needed in Section~\ref{sec8}.
Again we may replace $q^{-2n}$ in \eqref{phi54phi32} by $c$ to find
\begin{multline*}
{_5\phi_4}\biggl[\genfrac{}{}{0pt}{}{a,aq^2,b^2,c,cq^2}
{abq^{-1},abq,cq^3,cq^5};q^4,q^4\biggr] \\
=\frac{(q^3,cq^3/a,cq^3/b,q^3/ab;q^2)_{\infty}}
{(cq^3,q^3/a,q^3/b,cq^3/ab;q^2)_{\infty}}
{_3\phi_2}\biggl[\genfrac{}{}{0pt}{}{a,b,c}
{aq^{-1},cq^3/b};q^2,-\frac{q^2}{b}\biggr],
\end{multline*}
provided both series terminate. 

To the best of our knowledge \eqref{fin54} and
\eqref{phi54phi32} are new, and
the result closest to these transformations that we were able to obtain using
just elementary results from \cite{GR90} is
\begin{equation*}
{_4\phi_3}\biggl[\genfrac{}{}{0pt}{}{a,aq,q^{-n},q^{1-n}}
{b^2q,-aq^{1-n}/b,-aq^{2-n}/b};q^2,q^2\biggr]
=\frac{(b^2;q^2)_n}{(b^2,-b/a;q)_n}
{_2\phi_1}\biggl[\genfrac{}{}{0pt}{}{b,q^{-n}}{q^{1-n}/b};q,\frac{q}{a}\biggr].
\end{equation*}
This generalizes \cite[Exercise 1.6 (i)]{GR90} obtained when $a$ tends to
$0$, and follows from \cite[Exercise 3.4]{GR90} and \cite[Eq, (III.8)]{GR90}.
\begin{proof}[Proof of \eqref{var352}]
Take \eqref{p32} and let $j$ be the summation variable in the ${_3\phi_2}$
series. Replace $n$ by $n-k$, shift $j\to j-k$ and multiply both sides by 
\begin{equation*}
\frac{(x^2,y^2;q^2)_k}{(q^2,b^2x^2y^2;q^2)_k} 
\frac{(q^{-n};q)_k}{(b^2 q^{2-n};q)_k} (bq)^{2k}.
\end{equation*}
Next sum $k$ from $0$ to $n$ and interchange the order of the sums over
$j$ and $k$ on the left. This gives, after some tedious but elementary
manipulations involving $q$-shifted factorials,
\begin{multline*}
\sum_{j=0}^n 
{_3\phi_2}\biggl[\genfrac{}{}{0pt}{}{x^2,y^2,q^{-2j}}
{b^2x^2y^2,q^{2-2j}/b^2};q^2,q^2\biggr]
\frac{(b^2;q^2)_j(q^{-n};q)_j}
{(q^2;q^2)_j(b^2 q^{2-n};q)_j}q^{2j}\\
=q^n \frac{(q^{-2}/b^2;q^2)_n}{(-q,q^{-1}/b^2;q)_n}
{_3\phi_2}\biggl[\genfrac{}{}{0pt}{}{x^2,y^2,q^{-2n}}
{b^2x^2y^2,b^2q^{4-2n}};q^2,b^2q^3\biggr].
\end{multline*}
The ${_3\phi_2}$ can be summed by \eqref{qPS} resulting in \eqref{var352}.
\end{proof}

\subsubsection{Sextic transformations}

Both our results take the form of generalized commutation relations.
First,
\begin{equation*}
\sum_{s=r}^L F^{\eqref{qtrafo5}}_{2L,2s}(q)F^{\eqref{t2b}}_{s,r}(q)=
\sum_{\substack{s=r\\s\equiv r\;(2)}}^L
F^{\eqref{t2b}}_{L,s}(q^3)F^{\eqref{qtrafo5}}_{s,r}(q^2),
\end{equation*}
which, by the substitution $(-q^r,L-r)\to (a,n)$, yields
\begin{multline*}
{_5\phi_4}\biggl[\genfrac{}{}{0pt}{}
{a^{1/2},-a^{1/2},q^{-n},\omega q^{-n},\omega^2 q^{-n}}
{a,aq^{1/2},-aq^{1/2},q^{-3n}/a^2};q,q\biggr] \\
=\frac{1-a^3q^{3n}}{1-a^3}
\frac{(a^3;q^6)_n(-a^3q^3;q^3)_n}{(a^2q;q)_{3n}} \\ \times
{_5\phi_4}\biggl[\genfrac{}{}{0pt}{}
{a^2q^2,a^2q^4,a^2q^6,q^{-3n},q^{3-3n}}
{-a^3q^3,a^3q^6,-a^3q^6,q^{6-6n}/a^3};q^6,q^6\biggr].
\end{multline*}
Second,
\begin{equation*}
\sum_{\substack{s=r\\s\equiv r\;(2)}}^{\lfloor L/3\rfloor}
F^{\eqref{t3}}_{L,s}(q^2)F^{\eqref{qtrafo4}}_{s,r}(q^3)=
\sum_{\substack{s=3r\\s\equiv r\;(2)}}^L
F^{\eqref{qtrafo4}}_{L,s}(q)F^{\eqref{t3}}_{s,r}(q),
\end{equation*}
which, by the substitution $(q^{3r},(L-3r)/2)\to (a^2,n)$, yields
\begin{multline*}
{_5\phi_4}\biggl[\genfrac{}{}{0pt}{}
{a^{2/3},a^{2/3}\omega,a^{2/3}\omega^2,q^{-n},-q^{-n}}
{a,-a,aq^{1/2},q^{1/2-2n}/a};q,q\biggr] \\
=\frac{1-a^4q^{6n}}{1-a^4}
\frac{(aq^{1/2};q)_{2n}(a^2q^2;q^2)_n}{(a^4;q^6)_n} \\
\times
{_5\phi_4}\biggl[\genfrac{}{}{0pt}{}
{aq^{3/2},aq^{9/2},q^{-2n},q^{2-2n},q^{4-2n}}
{a^2q^3,-a^2q^3,-a^2q^6,q^{6-6n}/a^4};q^6,q^6\biggr].
\end{multline*}

\subsubsection{Transformation from $q$ to $q^9$}
As our second-last last relation there holds
\begin{equation*}
\sum_{\substack{s=3r\\s\equiv r\;(2)}}^L
F^{\eqref{qtrafo5}}_{L,s}(q)F^{\eqref{t3}}_{s,r}(q)=
\sum_{\substack{s=r\\s\equiv r\;(2)}}^{\lfloor L/3\rfloor}
F^{\eqref{t3}}_{L,s}(q^3)F^{\eqref{qtrafo5}}_{s,r}(q^3).
\end{equation*}
After replacing $(q^{3r},(L-3r)/2)\to (a^2,n)$ this becomes
\begin{multline*}
{_6\phi_5}\biggl[\genfrac{}{}{0pt}{}
{a^{2/3},a^{2/3}\omega,a^{2/3}\omega^2,q^{-n},\omega q^{-n},\omega^2 q^{-n}}
{a,-a,aq^{1/2},-aq^{1/2},q^{-3n}/a^2};q,q\biggr] \\
=\frac{1-a^6q^{6n}}{1-a^6}\frac{(a^6;q^6)_n}{(a^2q;q)_{3n}} 
{_6\phi_5}\biggl[\genfrac{}{}{0pt}{}
{a^2q^3,a^2q^6,a^2q^9,q^{-3n},q^{3-3n},q^{6-3n}}
{a^3q^{9/2},-a^3q^{9/2},a^3q^9,-a^3q^9,q^{9-9n}/a^6};q^9,q^9\biggr].
\end{multline*}
To the best of our knowledge this is the first transformation
between the bases $q$ and $q^9$.

\subsubsection{Transformation from $q$ to $q^{12}$}
Also our very last relation is an isolated result because
$F^{\eqref{t4b}}$ commutes with all but $F^{\eqref{qtrafo5}}$;
\begin{equation*}
\sum_{s=r}^L F^{\eqref{qtrafo5}}_{2L,2s}(q)F^{\eqref{t4}}_{s,r}(q)=
\sum_{\substack{s=r\\s\equiv r\;(2)}}^L 
F^{\eqref{t4}}_{L,s}(q^3)F^{\eqref{qtrafo5}}_{s,r}(q^4).
\end{equation*}
Making the replacement $(q^{2r+1},L-r)\to (a,n)$ this corresponds to
\begin{multline*}
{_5\phi_4}\biggl[\genfrac{}{}{0pt}{}{\text{i}q^{-1/2},-\text{i}q^{-1/2},
q^{-n},\omega q^{-n},\omega^2 q^{-n}}
{-q,a^{1/2},-a^{1/2},q^{1-3n}/a};q,q^2\biggr] \\
=q^{3n}\frac{(-q^{-3},a^3q^3;q^6)_n}{(-q^3;q^3)_n(a;q)_{3n}}
{_5\phi_4}\biggl[\genfrac{}{}{0pt}{}
{a^2q^2,a^2q^6,a^2q^{10},q^{-6n},q^{6-6n}}
{a^3q^3,a^3q^9,-q^{9-6n},-q^{15-6n}};q^{12},q^{12}\biggr].
\end{multline*}
We believe this to be the first example of a 
transformation relating base $q$ to base $q^{12}$.

\section{Inverse transformations}\label{sec4}

\subsection{Main results}
When iterating any of the transformations of Section~\ref{sec2}
it is often important to start with an as simple as possible 
$q$-binomial identity as seed. One possible way to determine whether
a potential seed can actually be reduced is by applying the inverses 
of the transformations of Lemmas~\ref{lem1}--\ref{lem6}.

For a transformation of the type \eqref{sym2} we consider a formula
of the form
\begin{equation*}
q^{-\frac{1}{4}\gamma L^2}
\sum_{r=0}^L \f_{L,r}(q)\qbin{2r}{r-j}=
q^{-\frac{1}{4}\gamma j^2}\qbin{L}{\frac{1}{2}(L-j)}_{q^k}
\chi(L\equiv j\;(2))
\end{equation*}
as its inverse. Here $\chi$ is the truth function; $\chi(\text{true})=1$ and 
$\chi(\text{false})=0$.
Indeed, replacing $(L,r)\to (r,s)$ in \eqref{sym2} and then
using this to eliminate the $q$-binomial coefficient in the above summand
yields
\begin{equation*}
\sum_{\substack{s=0\\s\equiv j\;(2)}}^L 
q^{\frac{1}{4}\gamma(s^2-L^2)}
\qbin{s}{\frac{1}{2}(s-j)}_{q^k}
\sum_{r=s}^L \f_{L,r}(q)f_{r,s}(q)
=\qbin{L}{\frac{1}{2}(L-j)}_{q^k} \chi(L\equiv j\;(2)).
\end{equation*}
This is obviously satisfied if the inverse relations
\begin{subequations}\label{invrel}
\begin{align}\label{invrela}
\sum_{r=s}^L \f_{L,r}(q)f_{r,s}(q)=\delta_{L,s}, \\
\sum_{r=s}^L f_{L,r}(q)\f_{r,s}(q)=\delta_{L,s}
\end{align}
\end{subequations}
hold. Here the second equation follows from the first and
the fact that $f_{L,r}(q)$ is nonzero if and only if $0\leq r\leq L$.
Inverse relations like \eqref{invrel}
have been much studied in the theory of basic hypergeometric
series. Most importantly, they are related to the Bailey 
transform \cite{Andrews79,Andrews01,AB02,Bressoud83,Bressoud88,W01b}, 
the problem of $q$-Lagrange inversion \cite{GS83,GS86} and summations
and transformations of $q$-hypergeometric series 
\cite{AV84,Carlitz73,Chu94,Chu95,Gasper89,Krattenthaler96,Schlosser02}.

The first inverse is that of Lemma~\ref{lem1}.
\begin{lemma}\label{lem7}
For $L$ and $j$ integers there holds
\begin{multline*}
q^{-\frac{1}{2}L^2}\sum_{r=0}^L(-1)^{r+L}q^{\binom{L-r}{2}}
(-q;q)_{L-r}\qbin{L}{r}_{q^2}\qbin{2r}{r-j} \\
=q^{-\frac{1}{2}j^2}\qbin{L}{\frac{1}{2}(L-j)}_{q^2}\chi(L\equiv j\;(2)).
\end{multline*}
\end{lemma}

\begin{proof}
All we need to do is show that
\begin{equation*}
\f_{L,r}(q)=(-1)^{r+L}q^{\binom{L-r}{2}}
(-q;q)_{L-r}\qbin{L}{r}_{q^2}
\end{equation*}
and $f_{L,r}(q)$ as given by \eqref{fLr} satisfy \eqref{invrela}. 
Shifting $r\to r+s$ this becomes the $n\to L-s$ case of
${_1\phi_0}(q^{-n};\text{---};q,q)=\delta_{n,0}$, which
follows from the $q$-binomial theorem \cite[Eq. (II.4)]{GR90}
\begin{equation}\label{qbthm}
{_1\phi_0}(q^{-n};\text{---};q,z)=(zq^{-n};q)_n.
\end{equation}

Alternatively we can prove Lemma~\ref{lem7} without resorting to 
inverse relations. Assuming $0\leq j\leq L$ and shifting $r\to r+j$ 
the identity of the lemma becomes the $(a,c,n)\to(0,q^{j+1/2},L-j)$
instance of \cite[Eq. (II.17)]{GR90}
\begin{equation}\label{A17}
{_4\phi_3}\biggl[\genfrac{}{}{0pt}{}
{a^2q,c,-c,q^{-n}}{c^2,aq^{1-n/2},-aq^{1-n/2}};q,q\biggr]=
\frac{(q,c^2/a^2;q^2)_{n/2}}{(c^2q,1/a^2;q^2)_{n/2}}\chi(n\equiv 0\;(2))
\end{equation}
due to Andrews \cite{Andrews76b}.
\end{proof}

Next is the inverse of Lemma~\ref{lem2}.
\begin{lemma}\label{lem8}
For $L$ and $j$ integers there holds
\begin{multline*}
q^{-\frac{1}{4}L^2}\sum_{r=0}^L(-1)^{r+L} q^{\binom{L-r}{2}}
(-q^{2r-L+2};q^2)_{L-r}\qbin{L}{r}\qbin{2r}{r-j} \\
=q^{-\frac{1}{4}j^2}\qbin{L}{\frac{1}{2}(L-j)}_{q^2}\chi(L\equiv j\;(2)).
\end{multline*}
\end{lemma}

\begin{proof}
Using that \eqref{invrel} remains unchanged if we multiply $f_{L,r}(q)$ by
$x_r(q)y_L(q)$ and divide $\f_{L,r}(q)$ by $x_L(q)y_r(q)$ 
($x_r(q)\neq 0$, $y_L(q)\neq 0$) we this time need to show that
\begin{align*}
f_{L,r}(q)&=(1+q^L)\frac{(-q^{r+2};q^2)_{L-r-1}}{(q;q)_{L-r}}, \\[2mm]
\f_{L,r}(q)&=(-1)^{r+L}q^{\binom{L-r}{2}}
\frac{(-q^{2r-L+2};q^2)_{L-r}}{(q;q)_{L-r}}
\end{align*}
satisfies \eqref{invrela}. Shifting $r\to r+s$ this is \eqref{A17} with
$c=a$ and $(a^2,n)\to(-q^s,L-s)$.

Alternatively, we may assume $0\leq j\leq L$ and shift $r\to r+j$ 
to find that Lemma~\ref{lem8} is \eqref{A17} with 
$(a^2,c,n)\to (-q^j,q^{j+1/2},L-j)$.
\end{proof}

The following lemma, corresponding to the inverse of \eqref{t2c} is 
(literally) the odd one out as the sum on the left does not vanish 
when $L-j$ is odd.
\begin{lemma}\label{lem9}
For $L$ and $j$ integers such that $j\equiv L\pmod{2}$ there holds
\begin{multline*}
q^{-\frac{1}{4}L(L+2)}(1+q^L)\sum_{r=0}^L(-1)^{r+L} q^{\binom{L-r}{2}}
(-q^{2r-L+3};q^2)_{L-r}\qbin{L}{r}\qbin{2r}{r-j} \\
=q^{-\frac{1}{4}j(j+2)}(1+q^j)\qbin{L}{\frac{1}{2}(L-j)}_{q^2}.
\end{multline*}
\end{lemma}

\begin{proof}
The difference with the previous two cases is that the pair
\begin{align*}
f_{L,r}(q)&=\frac{(-q^{r+1};q^2)_{L-r}}{(q;q)_{L-r}},\\[2mm]
\f_{L,r}(q)&=(-1)^{r+L} q^{\binom{L-r}{2}}
\frac{(-q^{2r-L+3};q^2)_{L-r}}{(q;q)_{L-r}}
\end{align*}
only satisfies \eqref{invrel} for $s\equiv L\pmod{2}$.
Indeed, shifting $r\to r+s$ and substituting the above,
\eqref{invrela} becomes
\begin{multline}\label{phi54}
{_5\phi_4}\biggl[\genfrac{}{}{0pt}{}
{a^2,bq,c,-c,q^{-n}}{c^2,b,aq^{1-n/2},-aq^{1-n/2}};q,q\biggr] \\[1mm]
=\begin{cases} \displaystyle
\frac{1-a^2}{1-a^2 q^n}\frac{1-bq^n}{1-b}\;
\frac{(q,c^2/a^2;q^2)_{n/2}}{(c^2q,1/a^2;q^2)_{n/2}}  &
\text{if $n$ is even,} \\[3mm]
\displaystyle
\frac{1-q^n}{1-a^2 q^n}
\frac{a^2-b}{1-b}\;
\frac{(q,c^2q/a^2;q^2)_{(n-1)/2}}{(c^2q,q/a^2;q^2)_{(n-1)/2}}  &
\text{if $n$ is odd,}
\end{cases}
\end{multline}
with $c=a$, $b=0$ and $(a^2,n)\to (-q^{s+1},L-s(\equiv 0\;(2)))$.
Note in particular that for this choice of $a$, $b$ and $c$ the right
side of \eqref{phi54} only trivializes to $\delta_{n,0}$ for even values
of $n$, explaining why $L-s$ must be even. The proof of \eqref{phi54}
is given in the next subsection.

Also the direct proof of the lemma relies on a special case of 
\eqref{phi54}. Assuming $0\leq j\leq L$ and shifting $r\to r+j$
Lemma~\ref{lem9} is \eqref{phi54} with 
$(a^2,b,c,n)\to(-q^{j+1},0,q^{j+1/2},L-j\:(\in 2\Z))$. 
\end{proof}

The inverses of the two quartics transforms
\eqref{t4} and \eqref{t4b} are as follows.
\begin{lemma}\label{lem10}
For $L$ and $j$ integers there holds
\begin{equation*}
\sum_{r=0}^L(-1)^{r+L} (-q;q^2)_{L-r}\qbin{L}{r}_{q^2}\qbin{2r}{r-j}
=\qbin{L}{\frac{1}{2}(L-j)}_{q^4}\chi(L\equiv j\;(2)).
\end{equation*}
\end{lemma}

\begin{lemma}\label{lem11}
For $L$ and $j$ integers there holds
\begin{multline*}
(1+q^{2L})\sum_{r=0}^L(-1)^{r+L}q^{-r}(-q^{-1};q^2)_{L-r}
\qbin{L}{r}_{q^2}\qbin{2r}{r-j}\\
=q^{-j}(1+q^{2j})\qbin{L}{\frac{1}{2}(L-j)}_{q^4}\chi(L\equiv j\;(2)).
\end{multline*}
\end{lemma}

\begin{proof}
The Lemmas \ref{lem10} and \ref{lem11} follow 
from \eqref{t4} and \eqref{t4b} and the
$a=-q^{-1}$ and $a=-q$ instances of the inverse pair
\begin{subequations}\label{ffquartic}
\begin{align}
f_{L,r}(q)&=a^r\frac{(a;q^2)_{L-r}}{(q^2;q^2)_{L-r}}, \\[2mm]
\f_{L,r}(q)&=(1/a)^r\frac{(1/a;q^2)_{L-r}}{(q^2;q^2)_{L-r}}.
\end{align}
\end{subequations}
Shifting $r\to r+s$ in equation \eqref{invrel} this follows from the
$n\to L-s$ case of
\begin{equation}\label{phi21del}
{_2\phi_1}(a,q^{-2n};aq^{2-2n};q^2,q^2)=\delta_{n,0},
\end{equation}
which is a specialization of \eqref{qGauss}.

The direct proof of Lemmas~\ref{lem10} and \ref{lem11} is
only interesting for the latter. Namely, if we assume that
$0\leq j\leq L$ and shift $r\to r+j$ then Lemma~\ref{lem10} is
equation \eqref{A17} with $(a^2,c,n)\to (-q^{j-L-1},q^{j+1/2},L-j)$,
but  Lemma~\ref{lem11} is 
\begin{multline}\label{phi54p}
{_5\phi_4}\biggl[\genfrac{}{}{0pt}{}{c,-c,bq,q^{-n},-q^{-n}}
{c^2q,b,\textup{i}q^{3/2-n},-\textup{i}q^{3/2-n}};q,q^2\biggr] \\[1mm]
= \begin{cases} \displaystyle
\frac{(q^2;q^4)_{n/2}(-c^2q^{-1};q^2)_n} 
{(c^4q^2;q^4)_{n/2}(-q^{-1};q^2)_n} &
\text{if $n$ is even,} \\[3mm]
\displaystyle
\frac{1-q}{1-c^2q} \frac{c^2-b}{1-b}\;
\frac{(q^6;q^4)_{(n-1)/2}(-c^2q;q^2)_{n-1}} 
{(c^4q^6;q^4)_{(n-1)/2}(-q;q^2)_{n-1}} &
\text{if $n$ is odd,}
\end{cases}
\end{multline}
with $b=c^2$ and $(c,n)\to (q^{j+1/2},L-j)$. The identity \eqref{phi54p}
will be proven in Section~\ref{sec4b}.
\end{proof}

\begin{remark}
By \eqref{qGauss} it can also be shown that \eqref{invrel} with 
\eqref{ffquartic} (normalized) is the $b=1/a$ case of $M(a)M(b)=M(ab)$,
with $M(a)$ the infinite-dimensional, lower-triangular matrix 
$M(a)=(M_{i,j}(a))_{i,j\geq 0}$ whose entries are given by 
\begin{equation*}
M_{i,j}(a)=a^j(a;q^2)_{i-j}\qbin{i}{j}_{q^2}.
\end{equation*}
\end{remark}

Finally we state the `inverse' of the cubic transformation of
Lemma~\ref{lem6}.
\begin{lemma}\label{lem12}
For $L$ and $j$ integers such that $j\equiv L\pmod{2}$ there holds
\begin{multline*}\label{inv4}
q^{-\frac{3}{4}L^2}\sum_{\substack{r=L \\r\equiv j\;(2)}}^{3L}
(-1)^{\frac{1}{2}(r+L)} q^{\binom{\frac{1}{2}(3L-r)}{2}}
\frac{(q^{\frac{3}{2}(r-L+2)};q^3)_{\frac{1}{2}(3L-r)}}
{(q;q)_{\frac{1}{2}(3L-r)}}\qbin{r}{\frac{1}{2}(r-3j)} \\
=q^{-\frac{3}{4}j^2}\qbin{L}{\frac{1}{2}(L-j)}_{q^3}.
\end{multline*}
\end{lemma}

\begin{proof}
This case is quite different from the previous ones in that
$f_{L,r}(q)$ corresponding to \eqref{t3} is nonzero if and only if
$0\leq 3r\leq L$. As a consequence only a left-inverse exists,
and we claim that
\begin{align*}
f_{L,r}(q)&=\frac{(aq^3;q^3)_{\frac{1}{2}(L-r-2)}(1-aq^L)}
{(aq^3;q^3)_r(q;q)_{\frac{1}{2}(L-3r)}} \\[2mm]
\f_{L,r}(q)&=(-1)^{\frac{1}{2}(r+L)} 
q^{\binom{\frac{1}{2}(3L-r)}{2}}
\frac{(aq^{\frac{3}{2}(r-L+2)};q^3)_{\frac{1}{2}(3L-r)}}
{(q;q)_{\frac{1}{2}(3L-r)}}
\end{align*}
with $r\equiv L\pmod{2}$ satisfies 
\begin{equation}\label{ffdel}
\sum_{\substack{r=3s \\r\equiv 0\;(2)}}^{3L}
\f_{L,r}(q)f_{r,s}(q)=\delta_{L,s}
\end{equation}
for $s\equiv L\pmod{2}$. Note that this suffices to conclude 
Lemma~\ref{lem12} from \eqref{t3} by taking $a=1$. To prove that 
\eqref{ffdel} indeed holds we repace $r\to 2r+3s$ to arrive at the 
$(b,n)\to (aq^{3s},3(L-s)/2(\equiv 0\;(3)))$ case of 
\begin{equation}\label{del}
\sum_{r=0}^n\frac{1-bq^{2r}}{1-b}\frac{(b;q^3)_r(q^{-n};q)_r q^r}
{(q;q)_r(bq^{3-n};q^3)_r}=\delta_{n,0},
\end{equation}
which is \cite[Eq. (3.6.17); $p\to q^3$, $a\to 0$]{GR90}
due to Bressoud \cite{Bressoud88}, Gasper \cite{Gasper89} 
and Krattenthaler \cite{Krattenthaler96}. 

For a direct proof of the lemma we shift $r\to 2r+j$ to obtain the 
`singular case' $(c,n)\to (q^{3j+1},3(L-j)/2(\equiv 0\;(3)))$ of
\begin{equation}\label{GS432}
\sum_{r=0}^n\frac{(c;q)_{2r}(q^{-n};q)_r q^r}
{(q,c;q)_r(cq^{2-n};q^3)_r}=
\frac{(q,q^2;q^3)_{\frac{1}{3}n}}{(q/c,cq^2;q^3)_{\frac{1}{3}n}}
\chi(n\equiv 0\;(3)),
\end{equation}
which is \cite[Eq. (4.32); $k\to n-k$, $A\to q^{1-2n}/c$]{GS83}
of Gessel and Stanton.
\end{proof}

\subsection{Proofs of \eqref{phi54} and \eqref{phi54p}}\label{sec4b}

Before proving the what-we-believe-to-be new balanced ${_5\phi_4}$ sum 
\eqref{phi54} we note that Andrews' identity \eqref{A17} arises as the 
case $b=a^2$ (or $b=q^{-n}$). Since \eqref{A17} provides a $q$-analogue
of Watson's ${_3F_2}$ summation, \eqref{phi54} also provides a 
generalization of Watson's sum. Specifically, replacing 
$(a,b,c)\to(q^{a/2},q^b,q^c)$ in \eqref{phi54} and then letting $q$
tend to one we find
\begin{equation}\label{F43}
{_4F_3}\biggl[\genfrac{}{}{0pt}{}
{a,b+1,c,-n}{2c,b,\frac{1}{2}(a-n+2)};1\biggr]
=\begin{cases} \displaystyle
\frac{a(b+n)}{b(a+n)}\frac{(\frac{1}{2},c-\frac{1}{2}a)_{n/2}}
{(c+\frac{1}{2},-\frac{1}{2}a)_{n/2}}  &
\text{if $n$ is even,} \\[3.5mm]
\displaystyle
\frac{n(b-a)}{b(a+n)}
\frac{(\frac{1}{2},c-\frac{1}{2}a+\frac{1}{2})_{(n-1)/2}}
{(c+\frac{1}{2},\frac{1}{2}-\frac{1}{2}a)_{(n-1)/2}}  &
\text{if $n$ is odd,}
\end{cases}
\end{equation}
where we employ standard notation for hypergeometric series
\cite{AAR99,GR90,Slater66}.
For $b=a$ this yields Watson's (terminating) ${_3F_2}$ sum.
(Wipple extended Watson's result to nonterminating series, but at the
${_4F_3}$ level this no longer appears to be possible.)
At the end of this section another extension of Watson's sum is be given.

\begin{proof}[Proof of \eqref{phi54}]
It is not hard to establish \eqref{phi54} by application of
the contiguous relation \cite[Eq. (3.8)]{Krattenthaler93}
\begin{multline}\label{cont}
{_r\phi_s}\biggl[\genfrac{}{}{0pt}{}{aq,b,c,(A)}{(B)};q,z\biggr]
=\frac{(1-b)(a-c)}{(1-a)(b-c)}
{_r\phi_s}\biggl[\genfrac{}{}{0pt}{}{a,bq,c,(A)}{(B)};q,z\biggr] \\
-\frac{(1-c)(a-b)}{(1-a)(b-c)}
{_r\phi_s}\biggl[\genfrac{}{}{0pt}{}{a,b,cq,(Aq)}{(B)};q,z\biggr].
\end{multline}
Here $(A)$, $(B)$ and $(Aq)$ are shorthand notations for 
$a_1,\dots,a_{r-3}$ and $b_1,\dots,b_s$ and $a_1 q,\dots,a_{r-3}q$,
respectively.
Utilizing \eqref{cont} with $(a,b,c)\to (b,a^2,q^{-n})$, 
the left-hand side of \eqref{phi54} transforms into the sum of two 
$b$-independent ${_4\phi_3}$ series.
Both are summable by \eqref{A17} to yield the desired right-hand side.
\end{proof}

\begin{proof}[Proof of \eqref{phi54p}]
To show \eqref{phi54p} we split its left-hand side by \eqref{cont} with 
$(a,b,c)\to (b,-q^{-n},q^{-n})$ so that
\begin{multline*}
\text{LHS}\eqref{phi54p}=\frac{(1+q^n)(1-bq^n)}{2q^n(1-b)}
{_4\phi_3}\biggl[\genfrac{}{}{0pt}{}{c,-c,q^{-n},-q^{1-n}}
{c^2q,\textup{i}q^{3/2-n},-\textup{i}q^{3/2-n}};q,q^2\biggr] \\
-\frac{(1-q^n)(1+bq^n)}{2q^n(1-b)}
{_4\phi_3}\biggl[\genfrac{}{}{0pt}{}{c,-c,q^{1-n},-q^{-n}}
{c^2q,\textup{i}q^{3/2-n},-\textup{i}q^{3/2-n}};q,q^2\biggr].
\end{multline*}
Both the ${_4\phi_3}$ series on the right are summable by
\begin{multline}\label{phi43p}
{_4\phi_3}\biggl[\genfrac{}{}{0pt}{}{a^2,c,-c,q^{-n}}
{c^2q,aq^{1-n/2},-aq^{1-n/2}};q,q^2\biggr] \\[1mm]
=\begin{cases} \displaystyle
\frac{c^2-a^2q^n}{1-a^2q^n}
\frac{(q;q^2)_{n/2}(c^2q^2/a^2;q^2)_{n/2-1}}
{(c^2q;q^2)_{n/2}(q^2/a^2;q^2)_{n/2-1}}
& \text{if $n$ is even,} \\[3mm]
\displaystyle
\frac{1-a^2c^2q^n}{1-a^2q^n}
\frac{(q;q^2)_{(n+1)/2}(c^2q/a^2;q^2)_{(n-1)/2}}
{(c^2q;q^2)_{(n+1)/2}(q/a^2;q^2)_{(n-1)/2}}
& \text{if $n$ is odd,}
\end{cases}
\end{multline}
leading to the right side of \eqref{phi54p}.
To complete the proof we need to deal with \eqref{phi43p}.
By \cite[Eq. (2.3)]{Krattenthaler93}
\begin{multline*}
{_{r+1}\phi_r}\biggl[\genfrac{}{}{0pt}{}{(A)}{aq,(B)};q,z\biggr]
={_{r+1}\phi_r}\biggl[\genfrac{}{}{0pt}{}{(A)}{a,(B)};q,z\biggr] \\
-\frac{az}{(1-a)(1-aq)}\frac{\prod_{i=1}^{r+1}(1-A_i)}
{\prod_{i=1}^{r-1}(1-B_i)}
{_{r+1}\phi_r}\biggl[\genfrac{}{}{0pt}{}{(Aq)}{aq^2,(B)};q,z\biggr]
\end{multline*}
with $a \to c^2q$ the left side of \eqref{phi43p} can be written as
the sum of two ${_4\phi_3}$ series, both of which can be summed
by the $b\to\infty$ limit of \eqref{phi54}. This results in the right side 
of \eqref{phi43p}.
\end{proof}

To conclude this section we wish to point out that 
\eqref{phi54} is certainly not the only generalization
of \eqref{A17} that may be obtained using contiguous relations.
For example, by \eqref{A17} and \cite[Eq. (3.3)]{Krattenthaler93}
\begin{multline*}
{_{r+1}\phi_r}\biggl[\genfrac{}{}{0pt}{}{a,(A)}{b,(B)};q,z\biggr]
={_{r+1}\phi_r}\biggl[\genfrac{}{}{0pt}{}{a/q,(A)}{b/q,(B)};q,z\biggr] \\
+\frac{z(a-b)}{(q-b)(1-b)}
\frac{\prod_{i=1}^r(1-A_i)}{\prod_{i=1}^{r-1}(1-B_i)}
{_{r+1}\phi_r}\biggl[\genfrac{}{}{0pt}{}{a,(Aq)}{bq,(Bq)};q,z\biggr]
\end{multline*}
with $(a,b,(A),(B))\to 
(bq,c^2q,(a^2q,c,-c,q^{-n}),(b,aq^{1-n/2},-aq^{1-n/2}))$ 
it follows that
\begin{multline}\label{phi54b}
{_5\phi_4}\biggl[\genfrac{}{}{0pt}{}
{a^2q,bq,c,-c,q^{-n}}{c^2q,b,aq^{1-n/2},-aq^{1-n/2}};q,q\biggr] \\[1mm]
=\begin{cases} \displaystyle
\frac{(q,c^2/a^2;q^2)_{n/2}}{(c^2q,1/a^2;q^2)_{n/2}} &
\text{if $n$ is even,} \\[3mm]\displaystyle
\frac{c^2-b}{1-b}\frac{1-a^2q}{c^2-a^2q}
\frac{(q,c^2q^{-1}/a^2;q^2)_{(n+1)/2}}{(c^2q,q^{-1}/a^2;q^2)_{(n+1)/2}}
& \text{if $n$ is odd.}
\end{cases}
\end{multline}
For $b=c^2$ this simplifies to \eqref{A17} and for
$(a,b,c)\to(q^{a/2-1/2},q^b,q^c)$ together with $q\to 1$ it yields
\begin{equation*}
{_4F_3}\biggl[\genfrac{}{}{0pt}{}
{a,b+1,c,-n}{2c+1,b,\frac{1}{2}(a-n+1)};1\biggr]
=\begin{cases} \displaystyle
\frac{(\frac{1}{2},c-\frac{1}{2}a+\frac{1}{2})_{n/2}}
{(c+\frac{1}{2},\frac{1}{2}-\frac{1}{2}a)_{n/2}}  &
\text{if $n$ is even,} \\[3.5mm]
\displaystyle
\frac{a(b-2c)}{b(a-2c)}
\frac{(\frac{1}{2},c-\frac{1}{2}a)_{(n+1)/2}}
{(c+\frac{1}{2},-\frac{1}{2}a)_{(n+1)/2}}  &
\text{if $n$ is odd.}
\end{cases}
\end{equation*}
This is to be compared with \eqref{F43}. For $b=2c$ this is again Watson's
${_3F_2}$ sum.

Finally we remark that other balanced ${_4\phi_3}$ summations
than \eqref{A17} follow from \eqref{phi54} and \eqref{phi54b}.
Taking $b=c^2/q$ in \eqref{phi54} and $b=a^2q$ in \eqref{phi54b} 
leads to two more such results. Especially the latter is appealing 
as some factors on the left of \eqref{phi54b} nicely cancel leading to
\begin{equation*}
{_4\phi_3}\biggl[\genfrac{}{}{0pt}{}
{a^2,c,-c,q^{-n}}{c^2q,aq^{-n/2},-aq^{-n/2}};q,q\biggr]
=\begin{cases} \displaystyle
\frac{(q,c^2q^2/a^2;q^2)_{n/2}}{(c^2q,q^2/a^2;q^2)_{n/2}} &
\text{if $n$ is even,} \\[3mm]\displaystyle
\frac{(q,c^2q/a^2;q^2)_{(n+1)/2}}{(c^2q,q/a^2;q^2)_{(n+1)/2}}
& \text{if $n$ is odd,}
\end{cases}
\end{equation*}
where we have also replaced $a$ by $a/q$.

\section{The Bailey lemma}\label{sec5}

As alluded to in the introduction, the $q$-binomial transformations of 
the first two sections are closely related to Bailey's lemma. 
Presently we will make this more precise and restate our results 
in terms of transformations on Bailey pairs. 

First we recall the definition of a Bailey pair \cite{Bailey49}.
If $\alpha(a;q)=\{\alpha_L(a;q)\}_{L\geq 0}$ and
$\beta(a;q)=\{\beta_L(a;q)\}_{L\geq 0}$ are sequences such that
\begin{equation*}
\beta_L(a;q)=\sum_{r=0}^L \frac{\alpha_r(a;q)}{(q;q)_{L-r}(aq;q)_{L+r}},
\end{equation*}
then $(\alpha(a;q),\beta(a;q))$
is called a Bailey pair relative to $a$ and $q$. 
The Bailey lemma is the following powerful mechanism for generating new 
Bailey pairs \cite{Andrews84,Andrews85,Andrews01,AB02,Paule85,W01b}.
\begin{lemma}\label{lemBL0}
If $(\alpha(a;q),\beta(a;q))$ is a Bailey pair relative to $a$ and $q$,
then so is $(\alpha'(a;q),\beta'(a;q))$ given by
\begin{subequations}\label{BL0}
\begin{align}
\alpha'_L(a;q)&=\frac{(b,c;q)_L (aq/bc)^L}
{(aq/b,aq/c;q)_L}\alpha_L(a;q), \\
\beta'_L(a;q)&=\frac{(aq/bc;q)_L}{(q,aq/b,aq/c;q)_L} \sum_{r=0}^L 
\frac{(b,c,q^{-L};q)_r q^r}{(bcq^{-L}/a;q)_r}\beta_r(a;q).
\end{align}
\end{subequations}
\end{lemma}
For $(b,c)\to (\infty,\infty)$ and $(b,c)\to(\infty,-(aq)^{1/2})$ 
this is equivalent to \eqref{qtrafo} and \eqref{qtrafo2}, respectively.

Before we state similar such results arising from the transformations
of section~\ref{sec2} we recall the base-changing 
Bailey-pair transformations of Bressoud \textit{et al.} \cite{BIS00}.
The first result is (equivalent to) \cite[Thm. 2.2]{BIS00}.
\begin{lemma}
If $(\alpha(a;q),\beta(a;q))$ is a Bailey pair relative to $a$ and $q$,
then the pair $(\alpha'(a^2;q^2),\beta'(a^2;q^2))$ given by
\begin{subequations}\label{BIS1}
\begin{align}
\alpha'_L(a^2;q^2)&=\frac{(b;q)_L}{(aq/b;q)_L}
\Bigl(\frac{aq}{b}\Bigr)^L (-1)^L q^{\binom{L}{2}}\alpha_L(a;q), \\
\beta'_L(a^2;q^2)&=\frac{(-aq/b;q)_{2L}}
{(-aq;q)_{2L}(q^2,a^2q^2/b^2;q^2)_L}\sum_{r=0}^L 
\frac{(b;q)_r(q^{-2L};q^2)_r q^r}{(-bq^{-2L}/a;q)_r}\beta_r(a;q)
\end{align}
\end{subequations}
forms a Bailey pair relative to $a^2$ and $q^2$. 
\end{lemma}
For $b\to 0$, $b\to\infty$ and $b\to -(aq)^{1/2}$ this yields
the equations (E1), (E2) and Eq. (E3) of \cite{Stanton01}.
By some simple variable changes, (E2) and (E3) can be seen to be
equivalent to \eqref{qtrafo3} and \eqref{qtrafo4}.

The next result is (equivalent to) \cite[Thm. 2.3]{BIS00},
\cite[Eq. (T1)]{Stanton01} and \eqref{qtrafo5}.
\begin{lemma}
If $(\alpha(a;q),\beta(a;q))$ is a Bailey pair relative to $a$ and $q$,
then the pair $(\alpha'(a^3;q^3),\beta'(a^3;q^3))$ given by
\begin{subequations}\label{BIS2}
\begin{align}
\alpha'_L(a^3;q^3)&=a^L q^{L^2} \alpha_L(a;q), \\
\beta'_L(a^3;q^3)&=\frac{(aq;q)_{3L}}{(q^3;q^3)_L(a^3q^3;q^3)_{2L}}
\sum_{r=0}^L \frac{(q^{-3L};q^3)_r q^r}{(q^{-3L}/a;q)_r}\beta_r(a;q)
\end{align}
\end{subequations}
forms a Bailey pair relative to $a^3$ and $q^3$. 
\end{lemma}

To the above three lemmas we now add several new base-changing
Bailey lemmas. First is a Bailey-type lemma of a quadratic nature.
\begin{lemma}\label{lemBL2}
If $(\alpha(a;q),\beta(a;q))$ is a Bailey pair relative to $a$ and $q$,
then so is $(\alpha'(a;q),\beta'(a;q))$ given by
\begin{subequations}\label{BL2}
\begin{align}
\alpha'_{2L}(a;q)&=(-1)^L b^L q^{L^2}\frac{(aq/b;q^2)_L}
{(bq;q^2)_L}\alpha_L(a;q^2), \qquad  \alpha'_{2L+1}(a;q)=0, \\[2mm]
\beta'_L(a;q)&=\frac{(b;q^2)_L}{(q,b;q)_L(aq;q^2)_L}
\sum_{r=0}^{\lfloor L/2\rfloor}\frac{(aq/b;q^2)_r(q^{-L};q)_{2r}q^{2r}}
{(q^{2-2L}/b;q^2)_r}\beta_r(a;q^2).  \label{BL2b}
\end{align}
\end{subequations}
\end{lemma}
For $b\to 0$, $b\to\infty$, $b\to -a^{1/2}$ and $b\to -a^{1/2}q$
this corresponds to the even $j$ case of \eqref{t2}, \eqref{t2dual},
\eqref{t2b} and \eqref{t2c}.

\begin{proof}
Writing the nontrivial part of \eqref{BL2} as 
\begin{align*}
\alpha'_{2L}(a;q)&=h_L(a,b)\alpha_L(a;q^2), \\
\beta'_L(a;q)&=\sum_{r=0}^{\lfloor L/2\rfloor} f_{L,r}(a,b)\beta_r(a;q^2),
\end{align*}
the claim of the lemma boils down to showing that
\begin{equation*}
\sum_{r=s}^{\lfloor L/2\rfloor}\frac{f_{L,r}(a,b)}
{(q^2;q^2)_{r-s}(aq^2;q^2)_{r+s}}
=\frac{h_s(a,b)}{(q;q)_{L-2s}(aq;q)_{L+2s}}.
\end{equation*}
After shifting $r\to r+s$ this follows from \eqref{phi32sum} with
$(a,b,n)\to(aq^{4s+1},bq^{2s},L-2s)$.
\end{proof}

Next, \eqref{t4} and \eqref{t4b} for even $j$ correspond to the
following two quartic Bailey lemmas.
\begin{lemma}
If $(\alpha(a;q),\beta(a;q))$ is a Bailey pair relative to $a$ and $q$,
then so is $(\alpha'(a;q),\beta'(a;q))$ given by
\begin{subequations}\label{BL4}
\begin{align}
\alpha'_{2L}(a;q)&=\alpha_L(a^2;q^4),\qquad \alpha'_{2L+1}(a;q)=0, \\[2mm]
\beta'_L(a;q)&=\frac{q^L(-q^{-1};q^2)_L}{(q^2,aq;q^2)_L}
\sum_{r=0}^{\lfloor L/2\rfloor}\frac{(-aq^2,q^{-2L};q^2)_{2r}q^{4r}}
{(-q^{3-2L};q^2)_{2r}}\beta_r(a^2;q^4).\label{BL4b}
\end{align}
\end{subequations}
\end{lemma}
\begin{proof}
Copying the proof of Lemma~\ref{lemBL2} this follows from 
\eqref{phi43sum} with $(a,n)\to(-aq^{4s+1},L-2s)$.
\end{proof}

\begin{lemma}
If $(\alpha(a;q),\beta(a;q))$ is a Bailey pair relative to $a$ and $q$,
then so is $(\alpha'(a;q),\beta'(a;q))$ given by
\begin{subequations}\label{BL4p}
\begin{align}
\alpha'_{2L}(a;q)&=q^{2L}\frac{1+a}{1+aq^{4L}}\alpha_L(a^2;q^4),\qquad
\alpha'_{2L+1}(a;q)=0, \\[2mm]
\beta'_L(a;q)&=\frac{(-q;q^2)_L}{(q^2,aq;q^2)_L}
\sum_{r=0}^{\lfloor L/2\rfloor}\frac{(-a,q^{-2L};q^2)_{2r}q^{4r}}
{(-q^{1-2L};q^2)_{2r}}\beta_r(a^2;q^4).\label{BL4pb}
\end{align}
\end{subequations}
\end{lemma}
\begin{proof}
This follows from \eqref{phi43sumb} with $(a,n)\to(-aq^{4s+1},L-2s)$.
\end{proof}

Finally there is cubic Bailey lemma corresponding to the even $j$ case of
the transformation \eqref{t3}.
\begin{lemma}\label{lemBL3}
If $(\alpha(a;q),\beta(a;q))$ is a Bailey pair relative to $a$ and $q$,
then so is $(\alpha'(a;q),\beta'(a;q))$ given by
\begin{subequations}\label{BL3}
\begin{align}
\alpha'_{3L}(a;q)&=a^L q^{3L^2}\alpha_L(a;q^3),\qquad
\alpha'_{3L\pm 1}(a;q)=0, \\[2mm]
\beta'_L(a;q)&=\frac{(a;q^3)_L}{(q;q)_L(a;q)_{2L}}
\sum_{r=0}^{\lfloor L/3\rfloor}\frac{(q^{-L};q)_{3r}q^{3r}}
{(q^{3-3L}/a;q^3)_r}\beta_r(a;q^3).\label{BL3b}
\end{align}
\end{subequations}
\end{lemma}
\begin{proof}
Writing the nontrivial part of \eqref{BL3} as
\begin{align*}
\alpha'_{3L}(a;q)&=h_L(a,b)\alpha_L(a;q^3), \\
\beta'_L(a;q)&=\sum_{r=0}^{\lfloor L/3\rfloor} f_{L,r}(a,b)\beta_r(a;q^3),
\end{align*}
we need to show that 
\begin{equation*}
\sum_{r=s}^{\lfloor L/3\rfloor}\frac{f_{L,r}(a,b)}
{(q^3;q^3)_{r-s}(aq^3;q^3)_{r+s}}
=\frac{h_s(a,b)}{(q;q)_{L-3s}(aq;q)_{L+3s}}.
\end{equation*}
Replacing $r\to r+s$ and defining $n=L-3s$, this follows from
\eqref{qPS} with $(a,b,c,d,q)\to (q^{-n},q^{1-n},q^{2-n},aq^{6s+3},q^3)$.
\end{proof}

\begin{remark}
The Lemmas~\ref{lemBL2}--\ref{lemBL3} correspond to the even $j$ instances of
the $q$-binomial transformations of Section~\ref{sec2}.
Equally well can one find Bailey-type lemmas corresponding to
$j$ being odd. Since we will not use these in the remainder of the paper
we only state the result related to the odd case of
\eqref{t2}, \eqref{t2dual}, \eqref{t2b} and \eqref{t2c}.
\begin{lemma}
If $(\alpha(a;q),\beta(a;q))$ is a Bailey pair relative to $a$ and $q$,
then so is $(\alpha'(a;q),\beta'(a;q))$ given by
\begin{align*}
\alpha'_{2L+1}(a;q)&=(-1)^L b^L q^{L^2}\frac{(aq^3/b;q^2)_L}
{(bq;q^2)_L}\alpha_L(aq^2;q^2), \qquad  \alpha'_{2L}(a;q)=0, \\[2mm]
\beta'_{L+1}(a;q)&=\frac{(b;q^2)_L}{(q,b;q)_L(aq^3;q^2)_L}
\sum_{r=0}^{\lfloor L/2\rfloor}\frac{(aq^3/b;q^2)_r(q^{-L};q)_{2r}q^{2r}}
{(aq;q)_2(q^{2-2L}/b;q^2)_r}\beta_r(aq^2;q^2).
\end{align*}
\end{lemma}
\end{remark}

\vskip 5mm 

Again it is important to also find the inverses of the
transformations \eqref{BL2}--\eqref{BL3}. Since all are of the form
\begin{align*}
\alpha'_{nL}(a;q)&=g_L(a,q)\alpha_L(a^k;q^l),\qquad
\alpha'_{nL+m}(a;q)=0, \\
\beta'_L(a;q)&=\sum_{r=0}^{\lfloor L/n\rfloor} f_{L,r}(a,q)
\beta_r(a^k;q^l),
\end{align*}
with 
$n\in\{2,3\}$ and $m\in\{1,\dots,n-1\}$ we can only find
left-inverses, defined as 
\begin{align*}
\alpha'_L(a^k;q^l)&=\frac{\alpha_{nL}(a;q)}{g_L(a,q)}, \\
\beta'_L(a^k;q^l)&=\sum_{r=0}^{nL}\f_{L,r}(a,q)\beta_r(a;q),
\end{align*}
with $\f_{L,r}(a,q)$ given by
\begin{equation}\label{leftinvrel}
\sum_{r=ns}^{nL}\f_{L,r}(a,q)f_{r,s}(a,q)=\delta_{L,s}.
\end{equation}
Such an $\f_{L,r}(a,q)$ is obviously not unique, but guided
by our inverse relations of the previous section it is not hard to
find an $\f_{L,r}(a,q)$ that can be expressed in simple closed form.

First we state a left-inverse of \eqref{BL2}.
\begin{lemma}
If $(\alpha(a;q),\beta(a;q))$ is a Bailey pair relative to $a$ and $q$,
then the pair $(\alpha'(a;q^2),\beta'(a;q^2))$ given by
\begin{subequations}
\begin{align}
\alpha'_L(a;q^2)&=(-1)^L b^{-L} q^{-L^2}\frac{(bq;q^2)_L}{(aq/b;q^2)_L}
\alpha_{2L}(a;q), \\[2mm]
\beta'_L(a;q^2)&=\frac{(1/b;q^2)_L}{(q;q)_{2L}(aq/b;q^2)_L}
\sum_{r=0}^{2L}\frac{(aq;q^2)_r(bq,q^{-2L};q)_rq^r}{(bq^{2-2L};q^2)_r}
\beta_r(a;q) \label{BL2invb}
\end{align}
\end{subequations}
yields a Bailey pair relative to $a$ and $q^2$.
\end{lemma}
\begin{proof}
Reading off $f_{L,r}(a,q)$ and $\f_{L,r}(a,q)$ from \eqref{BL2b}
and \eqref{BL2invb}, respectively, the inverse relation 
\eqref{leftinvrel} (with $n=2$) can be verified by \eqref{A17} 
with $c=a$ and $(a,n)\to (b^{1/2}q^s,2L-2s)$.
\end{proof}

Next we have left-inverses of the two quartic transformations
\eqref{BL4} and \eqref{BL4p}.
\begin{lemma}
If $(\alpha(a;q),\beta(a;q))$ is a Bailey pair relative to $a$ and $q$,
then the pair $(\alpha'(a^2;q^4),\beta'(a^2;q^4))$ given by
\begin{subequations}
\begin{align}
\alpha'_L(a^2;q^4)&=\alpha_{2L}(a;q), \\[2mm]
\beta'_L(a^2;q^4)&=
\frac{(-q;q^2)_{2L}}{(q^2,-aq^2;q^2)_{2L}}
\sum_{r=0}^{2L} \frac{(aq,q^{-4L};q^2)_r q^r}{(-q^{1-4L};q^2)_r}
\beta_r(a;q) \label{BL4invb}
\end{align}
\end{subequations}
yields a Bailey pair relative to $a^2$ and $q^4$.
\end{lemma}
\begin{proof}
Reading off $f_{L,r}(a,q)$ and $\f_{L,r}(a,q)$ from \eqref{BL4b}
and \eqref{BL4invb}, equation \eqref{leftinvrel} (with $n=2$)
follows from \eqref{phi21del} with $(a,n)\to (-q,2L-2s)$.
\end{proof}

\begin{lemma}
If $(\alpha(a;q),\beta(a;q))$ is a Bailey pair relative to $a$ and $q$,
then the pair $(\alpha'(a^2;q^4),\beta'(a^2;q^4))$ given by
\begin{subequations}
\begin{align}
\alpha'_L(a^2;q^4)&=q^{-2L}\frac{1+aq^{4L}}{1+a}\alpha_{2L}(a;q), \\[2mm]
\beta'_L(a^2;q^4)&=
\frac{(-q^{-1};q^2)_{2L}}{(q^2,-a;q^2)_{2L}}
\sum_{r=0}^{2L}\frac{(aq,q^{-4L};q^2)_rq^{2r}}{(-q^{3-4L};q^2)_r}
\beta_r(a;q).\label{BL4pinvb}
\end{align}
\end{subequations}
yields a Bailey pair relative to $a^2$ and $q^4$.
\end{lemma}
\begin{proof}
Reading off $f_{L,r}(a,q)$ and $\f_{L,r}(a,q)$ from \eqref{BL4pb}
and \eqref{BL4pinvb}, equation \eqref{leftinvrel} (with $n=2$)
is \eqref{phi21del} with $(a,n)\to (-q^{-1},2L-2s)$.
\end{proof}

Finally we give a left-inverse of \eqref{BL3}.
\begin{lemma}
If $(\alpha(a;q),\beta(a;q))$ is a Bailey pair relative to $a$ and $q$,
then the pair $(\alpha'(a;q^3),\beta'(a;q^3))$ given by
\begin{subequations}
\begin{align}
\alpha'_L(a;q^3)&=a^{-L} q^{-3L^2}\alpha_{3L}(a;q), \\[2mm]
\beta'_L(a;q^3)&=\frac{(1/a;q^3)_L}{(q;q)_{3L}}
\sum_{r=0}^{3L} \frac{(aq;q)_{2r}(q^{-3L};q)_r q^r}{(aq^{3-3L};q^3)_r}
\beta_r(a;q)\label{BL3invb}
\end{align}
\end{subequations}
yields a Bailey pair relative to $a$ and $q^3$.
\end{lemma}
\begin{proof}
Reading off $f_{L,r}(a,q)$ and $\f_{L,r}(a,q)$ from \eqref{BL3b}
and \eqref{BL3invb}, equation \eqref{leftinvrel} (with $n=3$)
follows from \eqref{del} with $(b,n)\to (aq^{6s},3L-3s)$.
\end{proof}

\section{$q$-Hypergeometric transformations}\label{sec6}

\subsection{Applications of base-changing Bailey lemmas}\label{sec61}
In the following we have compiled a list of quadratic, cubic
and quartic transformation formulas obtained by
twice iterating the unit Bailey pair \cite{Andrews84}
\begin{equation}\label{ubp}
\alpha_L=(-1)^L q^{\binom{L}{2}}\frac{1-aq^{2L}}{1-a}
\frac{(a;q)_L}{(q;q)_L},\qquad \beta_L=\delta_{L,0}
\end{equation}
using the Lemmas~\ref{lemBL0}--\ref{lemBL3}.
Before stating the resulting transformations two remarks are in order.

First we note that \eqref{BL4p} applied to the unit Bailey pair yields
\begin{equation*}
\alpha'_{2L}(a;q)=(-1)^L q^{2L^2}\frac{1-aq^{4L}}{1-a}
\frac{(a^2;q^4)_L}{(q^4;q^4)_L},\quad
\beta'_L(a;q)=\frac{(-q;q^2)_L}{(q^2,aq;q^2)_L},
\end{equation*}
whereas \eqref{BL2} applied to the unit Bailey pair leads to
\begin{equation*}
\alpha'_{2L}(a;q)=b^L q^{2L^2-1}\frac{1-aq^{4L}}{1-a}
\frac{(a,aq/b;q^2)_L}{(q^2,bq;q^2)_L},\quad
\beta'_L(a;q)=\frac{(b;q^2)_L}{(q,b;q)_L(aq;q^2)_L},
\end{equation*}
where in both cases $\alpha'_{2L+1}(a;q)=0$. Since the second result includes
the first as the special case $b=-q$, we need not consider those identities
obtain by first applying \eqref{BL4p} to the unit Bailey pair.

A second remark is that taking the unit Bailey pair and applying the
transformation (5.i) followed by \eqref{BL0} and then using a standard 
polynomial argument yields a result that 
implies the identity obtained by applying (5.i) followed by (5.k).
Here $i\in\{1,\dots,7\}$ and $k\in\{2,3,4,7\}$. So, for example,
we will not consider the identity obtained by successive application
of (5.3) and (5.7) because, modulo a polynomial argument,
it is implied by the application of (5.3) followed by (5.1).

Taking both of the above comments into account we will derive
18 different results, obtained by first applying 
(5.i) with $i\in\{1,2,3,4,5,7\}$ and then (5.k) with $k\in\{1,5,6\}$.
Five of the six identities that arise by application of (5.i)
followed by \eqref{BL0} are not new. We nevertheless have chosen 
to state these known results as they will be needed later
to prove several of the claims made in Section~\ref{sec33}.

\subsubsection{Transformation \textup{(5.i)} followed by \eqref{BL0}}
Applying \eqref{BL0} twice to the unit Bailey pair and
replacing $q^{-n}$ by $f$ we find Watson's transformation
\cite[Eq. (III.17)]{GR90}
\begin{multline}\label{Watson}
{_8W_7}(a;b,c,d,e,f;q,a^2q^2/bcdef) \\
=\frac{(aq,aq/de,aq/df,aq/ef;q)_{\infty}}
{(aq/d,aq/e,aq/f,aq/def;q)_{\infty}}
{_4\phi_3}\biggl[\genfrac{}{}{0pt}{}{aq/bc,d,e,f}
{aq/b,aq/c,def/a};q,q\biggr],
\end{multline}
provided both series terminate. 
(Watson's transformation actually holds under slightly weaker conditions,
but these do not follow from the above derivation.)
The derivation of \eqref{Watson} using Bailey's lemma is
of course well-known, see e.g., \cite{Andrews84}.

Applying \eqref{BIS1} and then \eqref{BL0} to the unit Bailey pair and
replacing $q^{-2n}$ by $e$ we obtain a quadratic transformation due to 
Verma and Jain \cite[Eq. (1.3)]{VJ80};
\begin{multline}\label{VJ13}
{_{10}W_9}(a;b,c^{1/2},-c^{1/2},d^{1/2},-d^{1/2},e^{1/2},-e^{1/2};
q,-a^3q^3/bcde) \\
=\frac{(a^2q^2,a^2q^2/cd,a^2q^2/ce,a^2q^2/de;q^2)_{\infty}}
{(a^2q^2/c,a^2q^2/d,a^2q^2/e,a^2q^2/cde;q^2)_{\infty}} \\[1mm]
\times
{_5\phi_4}\biggl[\genfrac{}{}{0pt}{}{-aq/b,-a q^2/b,c,d,e}
{-aq,-aq^2,a^2q^2/b^2,cde/a^2};q^2,q^2\biggr],
\end{multline}
provided both series terminate.
We remark that Verma and Jain stated the above identity for
$e=q^{-2n}$ only. In the calculations of Section~\ref{sec62}, the above, 
slightly more general form, will however be crucial.
Similar remarks apply to all the subsequent identities of Verma and Jain.

Applying \eqref{BIS2} and then \eqref{BL0} to the unit Bailey pair and
replacing $q^{-3n}$ by $d$ we obtain the following
cubic transformation of Verma and Jain \cite[Eq. (1.5)]{VJ80}:
\begin{multline}\label{VJ15}
{_{12}W_{11}}(a;b^{1/3},b^{1/3}\omega,b^{1/3}\omega^2,
c^{1/3},c^{1/3}\omega,c^{1/3}\omega^2,d^{1/3},d^{1/3}\omega,d^{1/3}\omega^2;
q,a^4q^4/bcd) \\
=\frac{(a^3q^3,a^3q^3/bc,a^3q^3/bd,a^3q^3/cd;q^3)_{\infty}}
{(a^3q^3/b,a^3q^3/c,a^3q^3/d,a^3q^3/bcd;q^3)_{\infty}} \\[1mm]
\times
{_6\phi_5}\biggl[\genfrac{}{}{0pt}{}{aq,aq^2,aq^3,b,c,d}
{(aq)^{3/2},-(aq)^{3/2},a^{3/2}q^3,-a^{3/2}q^3,bcd/a^3};q^3,q^3\biggr],
\end{multline}
provided both series terminate.

Applying \eqref{BL2} and then \eqref{BL0} to the unit Bailey pair and
replacing $q^{-n}$ by $e$ we find a second quadratic transformation of
Verma and Jain \cite[Eq. (1.4)]{VJ80};
\begin{multline}\label{VJ14}
{_{10}W_9}(a;aq/b,c,cq,d,dq,e,eq;q^2,a^2bq^2/c^2d^2e^2) \\
=\frac{(aq,aq/cd,aq/ce,aq/de;q)_{\infty}}{(aq/c,aq/d,aq/e,aq/cde;q)_{\infty}}
{_5\phi_4}\biggl[\genfrac{}{}{0pt}{}{b^{1/2},-b^{1/2},c,d,e}
{(aq)^{1/2},-(aq)^{1/2},b,cde/a};q,q\biggr],
\end{multline}
provided that both series terminate. 

Next, \eqref{BL3} followed by \eqref{BL0} yields a second
cubic transformation of Verma and Jain \cite[Eq. (1.6)]{VJ80};
\begin{multline}\label{VJ16}
{_{12}W_{11}}(a;b,bq,bq^2,c,cq,cq^2,d,dq,dq^2;q^3,a^4q^3/b^3c^3d^3) \\
=\frac{(aq,aq/bc,aq/bd,aq/cd;q)_{\infty}}
{(aq/b,aq/c,aq/d,aq/bcd;q)_{\infty}} \\
\times {_6\phi_5}\biggl[\genfrac{}{}{0pt}{}{
a^{1/3},a^{1/3}\omega,a^{1/3}\omega^2,b,c,d}
{a^{1/2},-a^{1/2},(aq)^{1/2},-(aq)^{1/2},bcd/a};q,q\biggr],
\end{multline}
provided both series terminate.

Finally, the identity obtained from \eqref{BL4} followed by \eqref{BL0}
appears to be new;
\begin{multline}\label{BL4BL0}
\sum_{k=0}^{\infty} \frac{1-a^2q^{8k}}{1-a^2}
\frac{(a^2;q^4)_k}{(q^4;q^4)_k}
\frac{(b,c,d;q)_{2k}}{(aq/b,aq/c,aq/d;q)_{2k}}
\Bigl(-\frac{a^2q}{b^2c^2d^2}\Bigr)^k \\
=\frac{(aq,aq/bc,aq/bd,aq/cd;q)_{\infty}}{(aq/b,aq/c,aq/d,aq/bcd;q)_{\infty}}
{_5\phi_4}\biggl[\genfrac{}{}{0pt}{}{\text{i}q^{-1/2},-\text{i}q^{-1/2},b,c,d}
{-q,(aq)^{1/2},-(aq)^{1/2},bcd/a};q,q^2\biggr],
\end{multline}
provided both series terminate.

Interesting summations occur by making `singular' specializations
in the above six identities. To illustrate the idea, consider \eqref{Watson}
and put the prefactor of the ${_4\phi_3}$ series to the left-hand side.
If $k$ is the summation variable of the ${_8W_7}$ series, then the
summand on the left contains 
\begin{equation*}
\frac{(aq/e;q)_{\infty}(f;q)_k}{(aq/e;q)_k}=(aq^{k+1}/e;q)_{\infty}(f;q)_k
\end{equation*}
as a factor. By taking $e=aq^n$ and $f=q^{-n}$, this becomes
\begin{equation*}
(q^{k-n+1};q)_{\infty}(q^{-n};q)_k
\end{equation*}
which vanishes unless $k=n$. The resulting identity is the
$q$-Pfaff--Saalsch\"utz formula \eqref{qPS} with $(a,b,c,d)\to 
(aq/bc,aq^n,q^{-n},aq/b)$.
Similarly, by taking $d=a^2q^{2n}$ and $e=q^{-2n}$ in \eqref{VJ13} and
then negating $a$ we find
\begin{equation*}
{_4\phi_3}\biggl[\genfrac{}{}{0pt}{}{aq/b,aq^2/b,a^2q^{2n},q^{-2n}}
{aq,aq^2,a^2q^2/b^2};q^2,q^2\biggr]
=\Bigl(\frac{aq}{b}\Bigr)^n
\frac{1-a}{1-aq^{2n}}\frac{(-q,b;q)_n}{(a,-aq/b;q)_n},
\end{equation*}
and by taking $c=a^3q^{3n}$ and $d=q^{-3n}$ in \eqref{VJ15} we get
\begin{multline*}
{_5\phi_4}\biggl[\genfrac{}{}{0pt}{}{aq,aq^2,aq^3,a^3q^{3n},q^{-3n}}
{(aq)^{3/2},-(aq)^{3/2},a^{3/2}q^3,-a^{3/2}q^3};q^3,q^3\biggr] \\
=(aq)^n \frac{1-aq^{2n}}{1-a^3q^{6n}}\frac{(q^3;q^3)_n(aq;q)_{n-1}}
{(q;q)_n(a^3q^3;q^3)_{n-1}}.
\end{multline*}
Next consider \eqref{VJ14} and \eqref{BL4BL0} and again put the prefactor
of the ${_5\phi_4}$ series to the other side to obtain the factor
\begin{equation*}
(aq^{2k+1}/c;q)_{\infty}(d;q)_{2k}
\end{equation*}
in the summand on the left. By specializing
$c=aq^n$ and $d=q^{-n}$ this vanishes unless $2k=n$.
The two ensuing identities are \eqref{A17} and, after the change
$a\to a^2q^{1-n}$,
\begin{equation*}
{_4\phi_3}\biggl[\genfrac{}{}{0pt}{}
{\text{i}q^{-1/2},-\text{i}q^{-1/2},a^2q,q^{-n}}
{-q,aq^{1-n/2},-aq^{1-n/2}};q,q^2\biggr]
=\frac{1+a^2q^{n+1}}{a^2q+q^n}
\frac{(q,-q/a^2;q^2)_{n/2}}{(-q^2,1/a^2;q^2)_{n/2}}\chi(n\equiv 0\;(2)),
\end{equation*}
respectively.
Finally, removing the prefactor on the right and
specializing $c=q^{-n}$ and $d=aq^n$ in \eqref{VJ16},
the summand on the left will contain the factor 
\begin{equation*}
(q^{3k-n+1};q)_{\infty}(q^{-n};q)_{3k}
\end{equation*}
so that the only nonvanishing contribution comes from $3k=n$. The resulting 
identity is \eqref{GS432} with order of summation reversed.

Equation \eqref{BL4BL0} has another noteworthy specialization.
Taking $b=-(aq)^{1/2}$ and $c=-d=q^{-n}$ the left simplifies to
${_5W_4}(a^2;q^{-2n},q^{2-2n};q^4,-aq^{4n})$ which sums to
$(a^2q^2;q^2)_n(-a;q^4)_n/((a^2q^2;q^4)_n(-a;q^2)_n)$ 
by Rogers' $q$-Dougall sum \cite[Eq. (II.20)]{GR90}. The
resulting identity is \eqref{phi43new} with $a\to (aq)^{1/2}$.

\subsubsection{Transformation \textup{(5.i)} followed by 
\eqref{BL4} or \eqref{BL4p}}
Applying \eqref{BL0} and then \eqref{BL4} leads to
\begin{multline}\label{BL0BL4}
\sum_{k=0}^{\lfloor n/2\rfloor} \frac{1-a^2q^{8k}}{1-a^2}
\frac{(a^2,b,c;q^4)_k}{(q^4,a^2q^4/b,a^2q^4/c;q^4)_k}
\frac{(q^{-n};q)_{2k}}{(aq^{n+1};q)_{2k}}
\Bigl(-\frac{a^2q^{2n+3}}{bc}\Bigr)^k \\
=q^n\frac{(-q^{-1};q^2)_n(aq;q)_n}{(-q;q)_n(aq;q^2)_n}
{_5\phi_4}\biggl[\genfrac{}{}{0pt}{}{a^2q^4/bc,-aq^2,-aq^4,q^{-2n},q^{2-2n}}
{a^2q^4/b,a^2q^4/c,-q^{3-2n},-q^{5-2n}};q^4,q^4\biggr].
\end{multline}
For $b=1$ this simplifies to \eqref{phi43sum}.
Likewise, applying \eqref{BL0} and then \eqref{BL4p} generalizes 
\eqref{phi43sumb};
\begin{multline}\label{BL0BL4p}
{_{10}W_9}(a;-a,b^{1/2},-b^{1/2},c^{1/2},-c^{1/2},q^{-n},q^{1-n};
q^2,-a^2q^{2n+5}/bc) \\
=\frac{(-q;q^2)_n(aq;q)_n}{(-q;q)_n(aq;q^2)_n}
{_5\phi_4}\biggl[\genfrac{}{}{0pt}{}{a^2q^4/bc,-a,-aq^2,q^{-2n},q^{2-2n}}
{a^2q^4/b,a^2q^4/c,-q^{1-2n},-q^{3-2n}};q^4,q^4\biggr].
\end{multline}
This transformation seems to be a hybrid of \eqref{VJ13} and \eqref{VJ14}.

{}From \eqref{BIS1} followed by \eqref{BL4} we obtain
a ${_6\phi_5}$ to ${_4\phi_3}$ transformation. Summing the
${_6\phi_5}$ by \cite[Eq. (II.20)]{GR90} we recover \eqref{phi43sum}.
A similar kind of situation, but with a much happier outcome, arises
if we apply \eqref{BIS1} followed by \eqref{BL4p}.
In first instance this yields
\begin{multline*}
{_8W_7}(a;\text{i}a^{1/2},-\text{i}a^{1/2},b,q^{-n},q^{1-n};
q^2,aq^{2n+3}/b) \\
=\frac{(-q;q^2)_n(aq;q)_n}{(-q;q)_n(aq;q^2)_n}
{_5\phi_4}\biggl[\genfrac{}{}{0pt}{}{-a,-aq^2/b,-aq^4/b,q^{-2n},q^{2-2n}}
{a^2q^4/b^2,-aq^4,-q^{1-2n},-q^{3-2n}};q^4,q^4\biggr],
\end{multline*}
but thanks to \eqref{Watson} the left-hand side may be simplified to
a ${_4\phi_3}$ series. After also replacing $a$ by $-a^2$ this gives
\begin{multline*}
{_4\phi_3}\biggl[\genfrac{}{}{0pt}{}{q^2,b,q^{-n},q^{1-n}}
{aq^2,-aq^2,-bq^{1-2n}/a^2};q^2,q^2\biggr] \\
=\frac{(-q;q^2)_n(-a^2q/b;q)_n}{(-q;q)_n(-a^2q/b;q^2)_n}
{_5\phi_4}\biggl[\genfrac{}{}{0pt}{}{a^2,a^2q^2/b,a^2 q^4/b,q^{-2n},q^{2-2n}}
{a^2 q^4,a^4q^4/b^2,-q^{1-2n},-q^{3-2n}};q^4,q^4\biggr].
\end{multline*}
This result, which for $b=1$ this reduces to \eqref{phi43sumb},
should be compared with \eqref{comp}.

Next, \eqref{BIS2} followed by \eqref{BL4} yields
\begin{multline*}
\sum_{k=0}^{\lfloor n/2\rfloor} \frac{1-a^2q^{8k}}{1-a^2}
\frac{(a^2;q^4)_k}{(q^4;q^4)_k}
\frac{(q^{-3n};q^3)_{2k}}{(a^3q^{3n+3};q^3)_{2k}}(-a^2 q^{6n+1})^k \\
=q^{3n}\frac{(-q^{-3};q^6)_n(a^3q^3;q^3)_n}{(-q^3;q^3)_n(a^3q^3;q^6)_n}
{_5\phi_4}\biggl[\genfrac{}{}{0pt}{}
{a^2q^4,a^2q^8,a^2 q^{12},q^{-6n},q^{6-6n}}
{a^3 q^6,a^3q^{12},-q^{9-6n},-q^{15-6n}};q^{12},q^{12}\biggr]
\end{multline*}
and \eqref{BIS2} followed by \eqref{BL4p} yields
\begin{multline*}
\sum_{k=0}^{\lfloor n/2\rfloor} \frac{1-a^2q^{8k}}{1-a^2}
\frac{(a^2;q^4)_k}{(q^4;q^4)_k}
\frac{(-a^3;q^{12})_k}{(-a^3q^{12};q^{12})_k}
\frac{(q^{-3n};q^3)_{2k}}{(a^3q^{3n+3};q^3)_{2k}}(-a^2 q^{6n+7})^k \\
=\frac{(-q^3;q^6)_n(a^3q^3;q^3)_n}{(-q^3;q^3)_n(a^3q^3;q^6)_n}
{_6\phi_5}\biggl[\genfrac{}{}{0pt}{}
{-a^3,a^2q^4,a^2q^8,a^2 q^{12},q^{-6n},q^{6-6n}}
{a^3 q^6,a^3q^{12},-a^3q^{12},-q^{3-6n},-q^{9-6n}};q^{12},q^{12}\biggr].
\end{multline*}

{}From \eqref{BL2} followed by \eqref{BL4} we obtain yet another
generalization of \eqref{phi43sum} (obtained by setting $b=a^2 q^4$);
\begin{multline*}
\sum_{k=0}^{\lfloor n/4\rfloor} \frac{1-a^2q^{16k}}{1-a^2}
\frac{(a^2,a^2q^4/b;q^8)_k}{(q^8,bq^4;q^8)_k}
\frac{(q^{-n};q)_{4k}}{(aq^{n+1};q)_{4k}}(b q^{4n-2})^k \\
=q^n\frac{(-q^{-1};q^2)_n(aq;q)_n}{(-q;q)_n(aq;q^2)_n}
{_5\phi_4}\biggl[\genfrac{}{}{0pt}{}{
b^{1/2},-b^{1/2},-a q^4,q^{-2n},q^{2-2n}}
{b,aq^2,-q^{3-2n},-q^{5-2n}};q^4,q^4\biggr],
\end{multline*}
and from \eqref{BL2} followed by \eqref{BL4p} we obtain
our last generalization of \eqref{phi43sumb};
\begin{multline*}
{_{10}W_9}(a;-a,aq^2/b^{1/2},-aq^2/b^{1/2},q^{-n},q^{1-n},q^{2-n},q^{3-n};
q^4,bq^{4n+2})\\
=\frac{(-q;q^2)_n(aq;q)_n}{(-q;q)_n(aq;q^2)_n}
{_5\phi_4}\biggl[\genfrac{}{}{0pt}{}{b^{1/2},-b^{1/2},-a,q^{-2n},q^{2-2n}}
{b,aq^2,-q^{1-2n},-q^{3-2n}};q^4,q^4\biggr].
\end{multline*}

Applying \eqref{BL3} and then \eqref{BL4} or \eqref{BL4p} gives
\begin{multline*}
\sum_{k=0}^{\lfloor n/6\rfloor} \frac{1-a^2q^{24k}}{1-a^2}
\frac{(a^2;q^{12})_k}{(q^{12};q^{12})_k}
\frac{(q^{-n};q)_{6k}}{(aq^{n+1};q)_{6k}}(-a^2 q^{6n-3})^k \\
=q^n\frac{(-q^{-1};q^2)_n(aq;q)_n}{(-q;q)_n(aq;q^2)_n}
{_6\phi_5}\biggl[\genfrac{}{}{0pt}{}{
a^{2/3},a^{2/3}\omega,a^{2/3}\omega^2,-a q^4,q^{-2n},q^{2-2n}}
{a,-a,aq^2,-q^{3-2n},-q^{5-2n}};q^4,q^4\biggr].
\end{multline*}
and
\begin{multline*}
{_{10}W_9}(a;-a,q^{-n},q^{1-n},q^{2-n},q^{3-n},q^{4-n},q^{5-n};q^6,
-a^2 q^{6n+3}) \\
=\frac{(-q;q^2)_n(aq;q)_n}{(-q;q)_n(aq;q^2)_n}
{_5\phi_4}\biggl[\genfrac{}{}{0pt}{}{
a^{2/3},a^{2/3}\omega,a^{2/3}\omega^2,q^{-2n},q^{2-2n}}
{a,aq^2,-q^{1-2n},-q^{3-2n}};q^4,q^4\biggr],
\end{multline*}
respectively.

Finally, applying \eqref{BL4} twice leads to
\begin{multline*}
\sum_{k=0}^{\lfloor n/4\rfloor} \frac{1-a^4q^{32k}}{1-a^4}
\frac{(a^4;q^{16})_k}{(q^{16};q^{16})_k}
\frac{(q^{-n};q)_{4k}}{(aq^{n+1};q)_{4k}}(-q^{4n-6})^k \\
=q^n\frac{(-q^{-1};q^2)_n(aq;q)_n}{(-q;q)_n(aq;q^2)_n}
{_5\phi_4}\biggl[\genfrac{}{}{0pt}{}
{\text{i}q^{-2},-\text{i}q^{-2},-aq^4,q^{-2n},q^{2-2n}}
{-q^4,aq^2,q^{3-2n},q^{5-2n}};q^4,q^8\biggr]
\end{multline*}
whereas \eqref{BL4} followed by \eqref{BL4p} yields
\begin{multline*}
\sum_{k=0}^{\lfloor n/4\rfloor} \frac{1-a^4q^{32k}}{1-a^4}
\frac{(a^4;q^{16})_k}{(q^{16};q^{16})_k}
\frac{(-a;q^8)_k}{(-aq^8;q^8)_k}
\frac{(q^{-n};q)_{4k}}{(aq^{n+1};q)_{4k}}(-q^{4n-2})^k \\
=\frac{(-q;q^2)_n(aq;q)_n}{(-q;q)_n(aq;q^2)_n}
{_5\phi_4}\biggl[\genfrac{}{}{0pt}{}
{\text{i}q^{-2},-\text{i}q^{-2},-a,q^{-2n},q^{2-2n}}
{-q^4,aq^2,q^{1-2n},q^{3-2n}};q^4,q^8\biggr].
\end{multline*}

\subsection{Generalized $\sum FF=\sum FF$ identities}\label{sec62}
Many of the results of Section~\ref{sec61} may be further exploited to
yield base-changing transformations between balanced or `almost'
balanced series. The idea is to take two of the transformations
from the previous section and to specialize the respective left-hand
sides such that they coincide. As a result the corresponding right-hand
sides may be equated leading to a new transformation.
This way one can for example rederive all of the transformations 
implied by the $\sum FF=\sum FF$ relations of Section~\ref{sec33}.
Instead, however, we will only prove those identities of
Section~\ref{sec33} that generalize $\sum FF=\sum FF$ transformations.

As a first example we consider Watson's 
transformation \eqref{Watson} and transformation
\eqref{VJ14} of Verma and Jain.
The respective left-hand sides are given by
\begin{equation*}
{_8W_7}(a;b,c,d,e,f;q,a^2q^2/bcdef)
\end{equation*}
and
\begin{equation*}
{_{10}W_9}(a;aq/b,c,cq,d,dq,e,eq;q^2,a^2bq^2/c^2d^2e^2).
\end{equation*}
By the substitution $(b,d,e,f,q)\to (cq,aq/b,q^{-n},q^{1-n},q^2)$
in the first and $(d,e)\to (-(aq)^{1/2},q^{-n})$
in the second expression,
both become
\begin{equation*}
{_8W_7}(a;aq/b,c,cq,q^{-n},q^{1-n};q^2,abq^{2n+1}/c^2).
\end{equation*}
Hence under these substitutions the respective right-hand sides
may be equated, resulting in the quadratic transformation
\begin{multline*}
{_4\phi_3}\biggl[\genfrac{}{}{0pt}{}{b^{1/2},-b^{1/2},c,q^{-n}}
{(aq)^{1/2},b,-cq^{1-n}/(aq)^{1/2}};q,q\biggr] \\
=\frac{(aq/c;q)_n(b;q^2)_n}{((aq)^{1/2},b,-(aq)^{1/2}/c;q)_n}
{_4\phi_3}\biggl[\genfrac{}{}{0pt}{}{aq/c^2,aq/b,q^{-n},q^{1-n}}
{aq/c,aq^2/c,q^{2-2n}/b};q^2,q^2\biggr].
\end{multline*}
By the variable change $(a,b,c)\to (a^2/q,c^2,b)$ this becomes \eqref{trafo1}.

Next consider the pair of identities \eqref{VJ13} and \eqref{VJ14}.
If in \eqref{VJ13} we let $(b,c,d,e,q)\to (aq/b,-aq^2,q^{-2n},q^{2-2n},q^2)$
and in \eqref{VJ14} we let $(c,d,e)\to (-(aq)^{1/2},q^{-n},-q^{-n})$
then both left-hand sides become
\begin{equation}\label{W87}
{_8W_7}(a;aq/b,q^{-n},-q^{-n},q^{1-n},-q^{1-n};q^2,abq^{4n+1}).
\end{equation}
Again we conclude that under the above substitutions
the right-hand sides of \eqref{VJ13} and \eqref{VJ14} may be equated.
The resulting transformation is \eqref{quartic} with $a$ replaced by
$(aq)^{1/2}$. Equation \eqref{eqfail} is found by noting that 
\eqref{W87} also arises from \eqref{VJ13} 
by letting $(b,c,d,e,q)\to (-bq,a^2q^2/b^2,q^{-2n},q^{2-2n},q^2)$.

In our following example we again equate appropriately specialized
right-hand sides of \eqref{VJ13} and \eqref{VJ14}, but
this time $(b,c,d,e,q)\to (aq/b,\infty,q^{-2n},q^{2-2n},q^2)$
in \eqref{VJ13} and $(c,d,e)\to (\infty,q^{-n},-q^{-n})$ in \eqref{VJ14}.
Since both the left sides become
\begin{equation*}
{_8\phi_9}\biggl[\genfrac{}{}{0pt}{}
{a,a^{1/2}q^2,-a^{1/2}q^2,aq/b,q^{-n},-q^{-n},q^{1-n},-q^{1-n}}
{a^{1/2},-a^{1/2},bq,aq^{n+1},-aq^{n+1},aq^{n+2},-aq^{n+2},0,0}
;q^2,a^2bq^{4n+3}\biggr]
\end{equation*}
we may again equate the right sides leading to \eqref{quartic2} with 
$a\to (aq)^{1/2}$.

Next we equate \eqref{VJ14} with \eqref{BL0BL4p}. In order to do so we
need to choose
$(b,c,d,e)\to (-q,b^{1/2},-b^{1/2},q^{-n})$ in the former
and $c\to bq^2$ in the latter.
Both left-hand sides then simplify to
\begin{equation*}
{_{10}W_9}(a;-a,b^{1/2},-b^{1/2},b^{1/2}q,-b^{1/2}q,q^{-n},q^{1-n};
q^2,-a^2q^{2n+3}/b^2).
\end{equation*}
The corresponding identity obtained by equating the respective
right-hand sides is \eqref{tr1p} with $a\to aq$.

Finally we treat the pair of identities \eqref{BL4BL0} and \eqref{BL0BL4}.
In the first we let $(b,c,d)\to (b^{1/2},-b^{1/2},q^{-n})$ and in the 
second we let $c\to bq^2$.
Then both the left sides become
\begin{equation*}
\sum_{k=0}^{\infty} \frac{1-a^2q^{8k}}{1-a^2}
\frac{(a^2;q^4)_k}{(q^4;q^4)_k}
\frac{(b;q^2)_{2k}}{(a^2q^2/b;q^2)_{2k}}
\frac{(q^{-n};q)_{2k}}{(aq^{n+1};q)_{2k}}
\Bigl(-\frac{a^2q^{2n+1}}{b^2}\Bigr)^k.
\end{equation*}
Accordingly, we may equate right-hand sides to find \eqref{tr1} with
$a\to aq$.

\section{Rogers--Ramanujan type identities}\label{sec7}

In this section the transformations of Section~\ref{sec2} are
applied to yield identities of the Rogers--Ramanujan type.

A first remark is that most of the single-sum Rogers--Ramanujan
identities that result when applying our new $q$-binomial
transformations are well-known and can nearly all be found in Slater's
compendium of 130 such identities \cite{Slater52}.
This should come as no surprise since in the large $L$ limit
most of our transformations reduce to sums implied by the ordinary
Bailey lemma \eqref{BL0}.
In view of recent work on a polynomial analogue of the Slater list
by Sills \cite{Sills02}, we note that our transformations 
give rise to rather natural polynomial versions of
many of the single-sum identities, different from those in \cite{Sills02}.
For example, all of the Rogers--Ramanujan identities in Slater's list that
have a product side of the form 
\begin{equation*}
\frac{(\pm q^a,\pm q^{b-a},q^b;q^b)_{\infty}}{(q^2;q^2)_{\infty}}
\quad \text{or} \quad
(-q;q^2)_{\infty}
\frac{(\pm q^a,\pm q^{b-a},q^b;q^b)_{\infty}}{(q^2;q^2)_{\infty}}
\end{equation*}
can be given a polynomial analogue by using the transformations \eqref{t2}
and \eqref{t2b}, respectively.
To give just one example of this we take \eqref{seed} as starting point
and first apply the quadratic transformation \eqref{qtrafo3} to get
a polynomial identity equivalent to G(1) in Slater's list of Bailey 
pairs \cite{Slater51}. Then applying \eqref{t2} or \eqref{t2b} we obtain
\begin{equation*}
\sum_{j=-L}^L (-1)^j q^{j(7j+1)/2}\qbin{2L}{L-2j}=
\sum_{n=0}^{\infty}
\frac{q^{2n^2}(q^2;q^2)_L(-q;q)_{L-2n}}{(q^2;q^2)_n(-q;q)_{2n}(q^2;q^2)_{L-2n}}
\end{equation*}
and
\begin{equation*}
\sum_{j=-L}^L (-1)^j q^{j(5j+1)/2}\qbin{2L}{L-2j}=
(1+q^L)\sum_{n=0}^{\infty} \frac{q^{n^2}(q;q)_L(-q^2;q^2)_{L-n-1}}
{(q^4;q^4)_n(q;q)_{L-2n}}.
\end{equation*}
Assuming $\abs{q}<1$, taking the large $L$ limit and using the Jacobi triple 
product identity \cite[Eq. (II.28)]{GR90}
\begin{equation}\label{Jacobi}
\sum_{k=-\infty}^{\infty}(-1)^k a^k q^{\binom{k}{2}}=(a,q/a,q;q)_{\infty},
\end{equation}
we find the Rogers--Selberg identity
\begin{equation}\label{S33}
\sum_{n=0}^{\infty}\frac{q^{2n^2}}{(q^2;q^2)_n(-q;q)_{2n}}=
\frac{(q^3,q^4,q^7;q^7)_{\infty}}{(q^2;q^2)_{\infty}}
\end{equation}
and Rogers' 
\begin{equation*}
\sum_{n=0}^{\infty} \frac{q^{n^2}}{(q^4;q^4)_n}=
(-q;q^2)_{\infty}\frac{(q^2,q^3,q^5;q^5)_{\infty}}{(q^2;q^2)_{\infty}},
\end{equation*}
labelled (33) and (20) in Slater's list, respectively.

To actually find results that are new one has to work a little
harder. One particularly nice example of a result that appears to
be new is a the following `perfect' Rogers--Ramanujan identity 
involving the bases $q$, $q^2$, $q^3$ and $q^6$:
\begin{equation}\label{perfect}
\sum_{n=0}^{\infty}\frac{q^{n(n+1)/2}(q^3;q^3)_n}
{(q;q)^2_n(q^3;q^2)_n}
=\frac{(q^6;q^6)_{\infty}}{(q;q)_{\infty}(q^3;q^2)_{\infty}}.
\end{equation}
To prove this we take two polynomial identities equivalent to
the Bailey pairs pairs J(1) and J(1)$-$J(2) \cite{Slater52}
\begin{equation}\label{J1}
\sum_{j=-\infty}^{\infty}(-1)^j q^{3j(3j+1)/2}\qbin{2L}{L-3j}=
\begin{cases}
1 & \text{$L=0$,} \\ \displaystyle 
(1+q^L)\frac{(q^3;q^3)_{L-1}}{(q;q)_{L-1}} & \text{$L>0$}
\end{cases}
\end{equation}
and
\begin{equation}\label{J12}
\sum_{j=-\infty}^{\infty}(-1)^j q^{9j(j+1)/2}\qbin{2L}{L-3j-1}=
q^{L-1}\frac{(q^3;q^3)_{L-1}}{(q;q)_{L-1}}\chi(L>0),
\end{equation}
and calculate the sum $\eqref{J1}+q^{L+1}\eqref{J12}$ to get
\begin{equation*}
\sum_{j=-\infty}^{\infty}(-1)^j q^{3j(3j+1)/2}\qbin{2L+1}{L-3j}=
\frac{(q^3;q^3)_L}{(q;q)_L}.
\end{equation*}
Replacing $q\to q^2$ and then applying the quadratic transformation
\eqref{t2b} gives
\begin{multline*}
\sum_{j=-\infty}^{\infty}(-1)^j q^{6j(3j+1)}\qbin{2L}{L-6j-1} \\
=(1+q^L)(q;q)_L\sum_{n=0}^{\infty}
\frac{q^{n(n+1)}(q^6;q^6)_n(-q;q^2)_{L-n-1}}
{(q^2;q^2)_n^2(q^2;q^4)_{n+1}(q;q)_{L-2n-1}}.
\end{multline*}
By \eqref{Jacobi} this yields \eqref{perfect} with $q\to q^2$ in the 
large $L$ limit.

An identity similar to \eqref{perfect} that is in
Slater's list as item (78) results if we apply \eqref{t2c} to \eqref{J1}. 
We include its derivation here to demonstrate that also
transformations of the type \eqref{t2c} and \eqref{t4b} may
be successfully exploited to derive Rogers--Ramanujan identities,
despite the factor $1+q^{aj}$ in the denominator on the right.
First, by \eqref{t2c} and \eqref{J1}
\begin{multline*}
2\sum_{j=-\infty}^{\infty}\frac{(-1)^j q^{6j(3j+1)}}{1+q^{6j}}
\qbin{2L}{L-6j} \\
=(-q;q^2)_L+
2\sum_{n=1}^{\infty} \frac{q^{n(n+1)}(q^6;q^6)_{n-1}(-q;q^2)_{L-n}}
{(-q;q^2)_n(q^2;q^2)_{n-1}}\qbin{L}{2n}.
\end{multline*}
By negating $j$ it follows that the term $1+q^{6j}$ on the left
cancels the prefactor $2$. Then taking the large $L$ limit 
and replacing $q^2$ by $q$ yields
\begin{equation}\label{S78}
1+2\sum_{n=1}^{\infty}
\frac{q^{n(n+1)/2}(q^3;q^3)_{n-1}}{(q;q)_n(q;q^2)_n(q;q)_{n-1}}
=\frac{(q^9;q^9)_{\infty}(q^9;q^{18})_{\infty}}
{(q;q)_{\infty}(q;q^2)_{\infty}}.
\end{equation}

For our final single-sum Rogers--Ramanujan result we first establish
the truth of a family of polynomial identities obtained previously by
Andrews \cite[Eq. (4.5)]{Andrews95} for $k=3$.
\begin{proposition}\label{prop1}
For integers $L\geq 0$, $k\geq 2$ and $i\in\{1,\dots,k-1\}$ there holds
\begin{multline*}
\sum_{j=-\infty}^{\infty}(-1)^j q^{j(3kj-k+2i)/2}\qbin{2L}{L-3j}_{q^{k/3}} \\
=\sum_{n=0}^{\lfloor L/3\rfloor} \frac{q^{kn^2}(q^i,q^{k-i};q^k)_n
(1-q^{2kL/3})(q^k;q^k)_{L-n-1}}{(q^k;q^k)_{2n}(q^{k/3};q^{k/3})_{L-3n}}.
\end{multline*}
\end{proposition}

\begin{proof}
According to the finite form of Jacobi's triple product identity 
\cite[p. 49]{Andrews76}
\begin{equation}\label{FJTP}
\sum_{j=-L}^L (-1)^j a^j q^{\binom{j}{2}}\qbin{2L}{L-j}=(a,q/a;q)_L.
\end{equation}
Replacing $q\to q^k$, applying the cubic transformation
\eqref{t3} with $q\to q^{k/3}$,
and specializing $a=q^i$ yields the claimed proposition.
\end{proof}
In the above, \eqref{t3} was used for even values of $L$ and $j$ even. 
Needed for the odd case is the polynomial identity
\begin{equation*}
\sum_{j=-L}^L (-1)^j a^j q^{\binom{j+1}{2}}\qbin{2L+1}{L-j}=
(1-q^L/a)(1/a,q/a;q)_L
\end{equation*}
which easily follows from \eqref{FJTP}.
Mimicing the earlier proof and then taking $a=q^{-i}$ results in
the odd version of Proposition~\ref{prop1}.
\begin{proposition}
For integers $L\geq 0$, $k\geq 2$ and $i\in\{1,\dots,k-1\}$ there holds
\begin{multline*}
\sum_{j=-\infty}^{\infty}(-1)^j q^{j(3kj+3k-2i)/2}
\qbin{2L+1}{L-3j-1}_{q^{k/3}} \\
=\sum_{n=0}^{\lfloor(L-1)/3\rfloor} \frac{q^{kn(n+1)}
(1-q^{kn+i})(q^i,q^{k-i};q^k)_n(1-q^{2kL/3+k/3})(q^k;q^k)_{L-n-1}}
{(q^k;q^k)_{2n+1}(q^{k/3};q^{k/3})_{L-3n-1}}.
\end{multline*}
\end{proposition}
Letting $L$ tend to infinity and using the Jacobi triple product identity
\eqref{Jacobi} gives Rogers--Ramanujan-type identities for modulus $3k$.
\begin{corollary}
For $k\geq 2$ and $i\in\{1,\dots,k-1\}$,
\begin{equation}\label{coreq1}
\sum_{n=0}^{\infty} \frac{q^{kn^2}(q^i,q^{k-i};q^k)_n}{(q^k;q^k)_{2n}}
=\frac{(q^{k+i},q^{2k-i},q^{3k};q^{3k})_{\infty}}
{(q^k;q^k)_{\infty}}
\end{equation}
and
\begin{equation*}
\sum_{n=0}^{\infty} \frac{q^{kn(n+1)}(1-q^{kn+i})
(q^i,q^{k-i};q^k)_n}{(q^k;q^k)_{2n+1}}
=\frac{(q^i,q^{3k-i},q^{3k};q^{3k})_{\infty}}{(q^k;q^k)_{\infty}}.
\end{equation*}
\end{corollary}
For $k=3$ this yields three modulus $9$ identities due to Bailey
\cite[Eqs. (1.6)--(1.8)]{Bailey47}.

Another natural choice for $a$ in all of the above would have been $a=-q^i$.
This would for example give the following companion to \eqref{coreq1}:
\begin{equation*}
\sum_{n=0}^{\infty} \frac{q^{kn^2}(-q^i,-q^{k-i};q^k)_n}{(q^k;q^k)_{2n}}
=\frac{(-q^{k+i},-q^{2k-i},q^{3k};q^{3k})_{\infty}}{(q^k;q^k)_{\infty}}.
\end{equation*}

The power of the transformations of Section~\ref{sec2}
in deriving new Rogers--Ramanu\-jan-type identities becomes
fully clear when considering multisum identities.
Because of the iterative nature of the Lemmas~\ref{lem1}--\ref{lem6}
a sheer endless number of elegant new multisum identities may
be obtained. In particular, by combining the results \eqref{qtrafo},
\eqref{qtrafo2}, \eqref{qtrafo3}--\eqref{qtrafo5} with \eqref{t2},
\eqref{t2b} and \eqref{t3}, each seed (initial $q$-binomial identity)
becomes the root of a octonary tree (modulo redundancies implied
by the relations of Section~\ref{sec3}) of polynomial Rogers--Ramanujan 
identities. 
As such roots one can take the identities implied by the
A--M families of Bailey pairs as tabulated by Slater \cite{Slater51,Slater52}
with the provision that only `independent' or `irreducible' Bailey
pairs should be employed. For example, the root \eqref{seed} has among its
sons or successors (polynomial identities equivalent to) the pairs B(1) 
(by \eqref{qtrafo}),
H(2) (by \eqref{qtrafo2}), G(1) (by \eqref{qtrafo3}),
L(2) (by \eqref{qtrafo4}), C(1) (by \eqref{t2}), I(7)+I(8) (by \eqref{t2c})
and J(1) of \eqref{J1} (by \eqref{t3}).
Moreover, depending on the fine-details of the polynomial identity 
associated with a particular node of the tree, one may also invoke the
Lemmas~\ref{lem3} and \ref{lem5} (and of course Lemma \ref{lem4}). 
Indeed, as we have already seen in the derivation of \eqref{S78},
the undesirable (from a Rogers--Ramanujan identities 
point of view) denominator terms $1+q^{aj}$ may actually cancel, permitting
the use of the Jacobi triple product identity in the large $L$ limit.
Also, a polynomial identity may arise that is well-suited for
further iteration.
For example, if we once again take \eqref{seed} as starting point an apply
\eqref{t2c} we get after simplification
\begin{equation*}
\sum_{j=-\infty}^{\infty}(-1)^j q^{2j^2}\qbin{2L}{L-2j}=(-q;q^2)_L.
\end{equation*}
which is a companion to the Bailey pair I(12). This same identity also 
follows by applying \eqref{t4b} to \eqref{seed}.
Lastly we note that the tree can be 
further enhanced by noting that the polynomial identities
obtained by application of \eqref{qtrafo}, \eqref{qtrafo2},
\eqref{qtrafo3} and \eqref{t2} are not `self-dual'.
That is, by replacing $q\to 1/q$ a different
identity results that yet again may be further iterated.
In the case of \eqref{t2} this is equivalent to also using its
companion \eqref{t2dual}. This way it can for example be seen that
also the Bailey pairs H(4) (dual to B(1)), G(4) (dual to G(1))
and C(4) (dual to C(1)) become part of the tree rooted in \eqref{seed}. 

Below we shall only give a representative sample of the multisum identities 
contained in the tree with \eqref{seed} as root, and
we encourage the reader to exploit their own favourite combination of 
transformations. 
In all our Rogers--Ramanujan series we assume $\abs{q}<1$. Furthermore,
unless stated otherwise we adopt the convention that $n_0:=L$ and $n_k:=0$.
\begin{theorem}\label{thmRR1}
For $k\geq 2$ there holds
\begin{multline*}
\sum_{n_1,\dots,n_{k-1}\geq 0}
\frac{q^{n_1^2+2n_2^2+\cdots+2^{k-2}n_{k-1}^2}}{(q;q)_{2n_1}}
\prod_{j=2}^k (-q^{2^{j-2}};q^{2^{j-2}})_{n_{j-1}-2n_j}
\qbin{n_{j-1}}{2n_j}_{q^{2^{j-1}}} \\
=\frac{(q^{\frac{1}{2}(4^k-2^k)},q^{\frac{1}{2}\cdot 4^k},q^{4^k-2^{k-1}};
q^{4^k-2^{k-1}})_{\infty}}{(q;q)_{\infty}}.
\end{multline*}
\end{theorem} 
For $k=2$ this is item (61) of Slaters list.
\begin{proof}
Let 
\begin{equation}\label{G}
G(L;\alpha,\beta,K;q)=G(L;\alpha,\beta,K)=
\sum_{j=-\infty}^{\infty}(-1)^j q^{Kj((\alpha+\beta)j+\alpha-\beta)/2}
\qbin{2L}{L-Kj}.
\end{equation}
A $k$-fold application of \eqref{t2} to \eqref{seed} yields the polynomial
identity
\begin{equation*}
G(L;2^k-1,2^k,2^k)
=\sum_{n_1,\dots,n_{k-1}\geq 0}\prod_{j=1}^k
q^{2^j n_j^2}(-q^{2^{j-1}};q^{2^{j-1}})_{n_{j-1}-2n_j}
\qbin{n_{j-1}}{2n_j}_{q^{2^j}}.
\end{equation*}
For $k=1$ this is the Bailey pair identity C(1).
Taking the large $L$ limit and replacing $q^2$ by $q$ yields the theorem
thanks to Jacobi's triple product identity \eqref{Jacobi}.
\end{proof}

\begin{theorem}
For $k\geq 2$ there holds
\begin{multline*}
\sum_{n_1,\dots,n_{k-1}\geq 0}
\frac{q^{n_1^2+2n_2^2+\cdots+2^{k-2}n_{k-1}^2}}
{2(q;q)_{2n_1}(-q^2;q^2)_{n_1-1}}
\prod_{j=2}^k
(-q^{2^j n_j};q^{2^j})_{n_{j-1}-2n_j}
\qbin{n_{j-1}}{2n_j}_{q^{2^{j-1}}} \\
=\frac{(q^{\frac{1}{2}(4^k-2^k)},q^{\frac{1}{2}(4^k+2^k)},
q^{4^k};q^{4^k})_{\infty}}{(q;q^2)_{\infty}(q^4;q^4)_{\infty}}
\end{multline*}
\end{theorem}
For $k=2$ this is item (71) of Slater's list.

\begin{proof}
A $k$-fold application of \eqref{t2b} to \eqref{seed} yields the polynomial
identity
\begin{multline}\label{Gplus}
G(L;\tfrac{1}{2}(2^k-1),\tfrac{1}{2}(2^k+1),2^k) \\
=\tfrac{1}{2}(1+q^L)\sum_{n_1,\dots,n_{k-1}\geq 0}\prod_{j=1}^k
q^{2^{j-1}n_j^2}(-q^{2^j n_j};q^{2^j})_{n_{j-1}-2n_j}
\qbin{n_{j-1}}{2n_j}_{q^{2^{j-1}}}.
\end{multline}
The large $L$ limit yields the desired theorem.
\end{proof}

\begin{theorem}
For $k\geq 2$ there holds
\begin{multline*}
\sum_{n_1,\dots,n_{k-1}\geq 0}
\frac{q^{n_1^2+3n_2^2+\cdots+3^{k-2}n_{k-1}^2}}{(q;q)_{2n_1}}\prod_{j=2}^k
\frac{(q^{3^{j-1}};q^{3^{j-1}})_{n_{j-1}-n_j-1}(1-q^{2.3^{j-2}n_{j-1}})}
{(q^{3^{j-2}};q^{3^{j-2}})_{n_{j-1}-3n_j}(q^{3^{j-1}};q^{3^{j-1}})_{2n_j}} \\
=\frac{(q^{\frac{1}{6}(9^k-3^k)},q^{\frac{1}{6}(9^k+3^k)},
q^{\frac{1}{3}\cdot 9^k};q^{\frac{1}{3}\cdot 9^k})_{\infty}}{(q;q)_{\infty}}.
\end{multline*}
\end{theorem}
For $k=2$ this is item (93) of Slater's list.

\begin{proof}
A $k$-fold application of \eqref{t3} to \eqref{seed} yields 
\begin{multline*}
G(L;\tfrac{1}{2}(3^k-1),\tfrac{1}{2}(3^k+1),3^k) \\
=\sum_{n_1,\dots,n_{k-1}\geq 0}\prod_{j=1}^k
\frac{q^{3^j n_j^2}(q^{3^j};q^{3^j})_{n_{j-1}-n_j-1}
(1-q^{2.3^{j-1}n_{j-1}})}
{(q^{3^{j-1}};q^{3^{j-1}})_{n_{j-1}-3n_j}(q^{3^j};q^{3^j})_{2n_j}}.
\end{multline*}
For $k=1$ this is \eqref{J1}. Taking the limit and replacing
$q^3\to q$ completes the proof.
\end{proof}

So far we have only presented Rogers--Ramanujan identities 
obtained by iterating one and the same transformation. 
Finally we state eight more theorems
that arise when alternating \eqref{qtrafo3} and \eqref{t2},
or \eqref{qtrafo} and \eqref{t2}.
Eight and not two families of identities result because it not only
matters with which transformation one starts, but also whether an
even or odd number of iterations is carried out.

The first four theorems occur by alternating \eqref{qtrafo3} and \eqref{t2}.
\begin{theorem}
For $k$ an odd integer such that $k\geq 3$ there holds
\begin{multline}\label{rth}
\sum_{n_1,\dots,n_{k-1}\geq 0}
\frac{q^{n_1^2}}{(q;q^2)_{2n_1}}
\prod_{\substack{j=2 \\j\equiv 0\;(2)}}^{k-1}
\frac{q^{2n_j^2+n_{j+1}^2}}
{(q;q)_{n_{j-1}-2n_j}(-q;q)_{2n_j}
(q^2;q^2)_{n_j-n_{j+1}}(q;q^2)_{n_{j+1}}} \\
=\frac{(q^{2^k-1},q^{2^k},q^{2^{k+1}-1};q^{2^{k+1}-1})_{\infty}}
{(q;q)_{\infty}}.
\end{multline}
\end{theorem}
\begin{proof}
Take \eqref{seed}, and in alternating fashion apply 
\eqref{qtrafo3} $(k+1)/2$ times and \eqref{t2} $(k-1)/2$ times,
starting with \eqref{qtrafo3}. This yields an identity
for $G(L;2^{(k-1)/2}(1-2^{-k}),2^{(k-1)/2},2^{(k-1)/2};q^2)$ 
which implies the theorem in the large $L$ limit.
\end{proof}

To concisely state the next theorem we depart from our earlier convention
that $n_0=L$ and below the term $(q;q)_{n_{j-1}-2n_j}$ for $j=1$ should
be interpreted as $1$.
\begin{theorem}\label{thmRR2}
For $k$ an even integer such that $k\geq 2$ there holds
\begin{multline*}
\sum_{n_1,\dots,n_{k-1}\geq 0}
\prod_{\substack{j=1 \\j\equiv 1\;(2)}}^{k-1}
\frac{q^{2n_j^2+n_{j+1}^2}}{(q;q)_{n_{j-1}-2n_j}(-q;q)_{2n_j}
(q^2;q^2)_{n_j-n_{j+1}}(q;q^2)_{n_{j+1}}} \\ 
=\frac{(q^{2^k-1},q^{2^k},q^{2^{k+1}-1};q^{2^{k+1}-1})_{\infty}}
{(q^2;q^2)_{\infty}}.
\end{multline*}
\end{theorem}
For $k=2$ this is \eqref{S33}.

\begin{proof}
Take \eqref{seed}, and in alternating fashion apply 
\eqref{qtrafo3} and \eqref{t2} each $k/2$ times
starting with \eqref{qtrafo3}. This yields an identity
for $G(L;2^{k/2}(1-2^{-k}),2^{k/2},2^{k/2})$ 
which implies the theorem in the large $L$ limit.
\end{proof}

\begin{theorem}
For $k$ an odd integer such that $k\geq 3$ there holds
\begin{multline*}
\sum_{n_1,\dots,n_{k-1}\geq 0}
\frac{q^{2n_1^2}}{(-q;q)_{2n_1}}
\prod_{\substack{j=2 \\j\equiv 0\;(2)}}^{k-1}
\frac{q^{n_j^2+2n_{j+1}^2}}{(q^2;q^2)_{n_{j-1}-n_j}(q;q^2)_{n_j}
(q;q)_{n_j-2n_{j+1}}(-q;q)_{2n_{j+1}}} \\
=\frac{(q^{2^{k+1}-2},q^{2^{k+1}},q^{2^{k+2}-2};q^{2^{k+2}-2})_{\infty}}
{(q^2;q^2)_{\infty}}.
\end{multline*}
\end{theorem}
We note that if we replace $q^2\to q$ then the right-hand side equals
the right-hand side of \eqref{rth}.

\begin{proof}
Take \eqref{seed}, and in alternating fashion apply 
\eqref{t2} $(k+1)/2$ times and \eqref{qtrafo3} $(k-1)/2$ times,
starting with \eqref{t2}. This yields an identity
for $G(L;2^{(k+1)/2}(1-2^{-k}),2^{(k+1)/2},2^{(k+1)/2})$ 
which implies the theorem in the large $L$ limit.
\end{proof}

In the next theorem $(q^2;q^2)_{n_0-n_1}=1$.
\begin{theorem}
For $k$ an even integer such that $k\geq 2$ there holds
\begin{multline*}
\sum_{n_1,\dots,n_{k-1}\geq 0}
\prod_{\substack{j=1 \\j\equiv 1\;(2)}}^{k-1}
\frac{q^{n_j^2+2n_{j+1}^2}}{(q^2;q^2)_{n_{j-1}-n_j}(q;q^2)_{n_j}
(q;q)_{n_j-2n_{j+1}}(-q;q)_{2n_{j+1}}} \\
=\frac{(q^{2^{k+1}-2},q^{2^{k+1}},q^{2^{k+2}-2};q^{2^{k+2}-2})_{\infty}}
{(q;q)_{\infty}}.
\end{multline*}
\end{theorem}
For $k=2$ this coincides with and Theorem~\ref{thmRR1}.

\begin{proof}
Take \eqref{seed}, and in alternating fashion apply 
\eqref{t2} and \eqref{qtrafo3} each $k/2$ times
starting with \eqref{t2}. This yields an identity
for $G(L;2^{k/2}(1-2^{-k}),2^{k/2},2^{k/2};q^2)$ 
which implies the theorem in the large $L$ limit.
\end{proof}

The next four results are obtained by alternating \eqref{qtrafo} 
and \eqref{t2}.
\begin{theorem}
For $k$ an odd integer such that $k\geq 3$ there holds
\begin{multline*}
\sum_{n_1,\dots,n_{k-1}\geq 0}
\frac{q^{n_1^2}}{(q;q^2)_{n_1}} \\ \times
\prod_{\substack{j=2 \\j\equiv 0\;(2)}}^{k-1}
\frac{q^{2^{j/2}(n_j^2+n_{j+1}^2)}}
{(q^{2^{j/2-1}};q^{2^{j/2-1}})_{n_{j-1}-2n_j}
(q^{2^{j/2}};q^{2^{j/2}})_{n_j-n_{j+1}}
(q^{2^{j/2}};q^{2^{j/2+1}})_{n_{j+1}}} \\
=\frac{(q^{3\cdot 2^{k-1}-2^{(k+1)/2}},q^{3\cdot 2^{k-1}-2^{(k-1)/2}},
q^{3\cdot 2^k-3.2^{(k-1)/2}};q^{3\cdot 2^k-3.2^{(k-1)/2}})_{\infty}}
{(q;q)_{\infty}}.
\end{multline*}
\end{theorem}
\begin{proof}
Take \eqref{seed}, and in alternating fashion apply 
\eqref{qtrafo} $(k+1)/2$ times and \eqref{t2} $(k-1)/2$ times,
starting with \eqref{qtrafo}. This yields an identity
for $G(L;3\cdot 2^{(k-1)/2}-2,3\cdot 2^{(k-1)/2}-1,2^{(k-1)/2};q)$ 
which implies the theorem in the large $L$ limit.
\end{proof}

Below, $(q^{1/2};q^{1/2})_{n_0-2n_1}=1$.
\begin{theorem}
For $k$ an even integer such that $k\geq 2$ there holds
\begin{multline*}
\sum_{n_1,\dots,n_{k-1}\geq 0} \\
\prod_{\substack{j=0 \\j\equiv 0\;(2)}}^{k-2}
\frac{q^{2^{j/2}(n_{j+1}^2+n_{j+2}^2)}}
{(q^{2^{j/2-1}};q^{2^{j/2-1}})_{n_j-2n_{j+1}}
(q^{2^{j/2}};q^{2^{j/2}})_{n_{j+1}-n_{j+2}}
(q^{2^{j/2}};q^{2^{j/2+1}})_{n_{j+2}}} \\
=\frac{(q^{2^k-2^{k/2}},q^{2^k-2^{k/2-1}},
q^{2^{k+1}-3\cdot 2^{k/2-1}};q^{2^{k+1}-3\cdot 2^{k/2-1}})_{\infty}}
{(q;q)_{\infty}}.
\end{multline*}
\end{theorem}
For $k=2$ this is the first Rogers--Ramanujan identity.

\begin{proof}
Take \eqref{seed}, and in alternating fashion apply 
\eqref{qtrafo} and \eqref{t2} each $k/2$ times
starting with \eqref{qtrafo}. This yields an identity
for $G(L;2^{k/2+1}-2,2^{k/2+1}-1,2^{k/2})$ 
which implies the theorem with $q\to q^2$ in the large $L$ limit.
\end{proof}

\begin{theorem}
For $k$ an odd integer such that $k\geq 3$ there holds
\begin{multline*}
\sum_{n_1,\dots,n_{k-1}\geq 0}
q^{n_1^2} \\ \times
\prod_{\substack{j=0 \\j\equiv 0\;(2)}}^{k-3}
\frac{q^{2^{j/2}(n_{j+2}^2+2n_{j+3}^2)}}
{(q^{2^{j/2}};q^{2^{j/2}})_{n_{j+1}-n_{j+2}}
(q^{2^{j/2}};q^{2^{j/2+1}})_{n_{j+2}}
(q^{2^{j/2}};q^{2^{j/2}})_{n_{j+2}-2n_{j+3}}} \\
=\frac{(q^{2^{k+1}-3\cdot 2^{(k-1)/2}},q^{2^{k+1}-2^{(k+1)/2}},
q^{2^{k+2}-5\cdot 2^{(k-1)/2}};q^{2^{k+2}-5\cdot 2^{(k-1)/2}})_{\infty}}
{(q;q)_{\infty}}.
\end{multline*}
\end{theorem}

\begin{proof}
Take \eqref{seed}, and in alternating fashion apply 
\eqref{t2} $(k+1)/2$ times and \eqref{qtrafo} $(k-1)/2$ times,
starting with \eqref{t2}. This yields an identity
for $G(L;2^{(k+3)/2}-3,2^{(k+3)/2}-2,2^{(k+1)/2})$ 
which implies the theorem with $q\to q^2$ in the large $L$ limit.
\end{proof}

Below, $(q;q)_{n_0-n_1}=1$.
\begin{theorem}
For $k$ an even integer such that $k\geq 2$ there holds
\begin{multline*}
\sum_{n_1,\dots,n_{k-1}\geq 0}
\prod_{\substack{j=0 \\j\equiv 0\;(2)}}^{k-2}
\frac{q^{2^{j/2}(n_{j+1}^2+2n_{j+2}^2)}}
{(q^{2^{j/2}};q^{2^{j/2}})_{n_j-n_{j+1}}
(q^{2^{j/2}};q^{2^{j/2+1}})_{n_{j+1}}
(q^{2^{j/2}};q^{2^{j/2}})_{n_{j+1}-2n_{j+2}}} \\
=\frac{(q^{3(2^{k}-2^{k/2})},q^{3\cdot 2^k-2^{k/2+1}},
q^{3\cdot 2^{k+1}-5\cdot 2^{k/2}};
q^{3\cdot 2^{k+1}-5\cdot 2^{k/2}})_{\infty}}{(q;q)_{\infty}}.
\end{multline*}
\end{theorem}
For $k=2$ this is identity (61) of Slater.

\begin{proof}
Take \eqref{seed}, and in alternating fashion apply 
\eqref{t2} and \eqref{qtrafo} each $k/2$ times
starting with \eqref{t2}. This yields an identity
for $G(L;3(2^{k/2}-1),3\cdot 2^{k/2}-2,2^{k/2})$ 
which implies the theorem in the large $L$ limit.
\end{proof}

\section{Rogers--Szeg\"o polynomials}\label{sec8}
The Bailey pairs F(3) and F(4) in Slater's list \cite{Slater51}
may be concisely written as the polynomial identity
\begin{equation}\label{F34}
\sum_{\substack{j=-L \\ j\equiv L\;(2)}}^L 
q^j\qbin{L}{\frac{1}{2}(L-j)}_{q^4}=q^{-L} (-q^2;q^2)_L.
\end{equation}
Hence, by the quartic transformation \eqref{t4},
\begin{equation*}
\sum_{j=-L}^L q^j\qbin{2L}{L-j}=
\sum_{r=0}^L q^{L-2r}(-q^2;q^2)_r(-q^{-1};q^2)_{L-r}\qbin{L}{r}_{q^2},
\end{equation*}
and by its cousin \eqref{t4b},
\begin{equation*}
\sum_{j=-L}^L \frac{q^{2j}}{1+q^{2j}}\qbin{2L}{L-j}=
\frac{1}{2}\sum_{r=0}^L (-1;q^2)_r(-q;q^2)_{L-r}\qbin{L}{r}_{q^2}.
\end{equation*}
By negating $j$ it follows that the left-hand side of this last identity
simplifies to $$\frac{1}{2}\sum_{j=-L}^L \qbin{2L}{L-j}.$$
If we then replace $j$ by $L-j$ in both formulas and $r\to L-r$
in the second formula, we get
\begin{equation*}
\sum_{j=0}^{2L}q^{-j}\qbin{2L}{j}=
\sum_{r=0}^L q^{-2r}(-q^2;q^2)_r(-q^{-1};q^2)_{L-r}\qbin{L}{r}_{q^2}
\end{equation*}
and
\begin{equation*}
\sum_{j=0}^{2L}\qbin{2L}{j}=
\sum_{r=0}^L (-q;q^2)_r(-1;q^2)_{L-r}\qbin{L}{r}_{q^2}.
\end{equation*}
Recalling the definition of the Rogers--Szeg\"o (RS) polynomials
$H_n(t)$ given by \cite{Andrews76}
\begin{equation}\label{RSdef}
H_n(t)=\sum_{j=0}^n t^j \qbin{n}{j},
\end{equation}
one may recognize the above as identities for $H_{2L}(q^{-1})$
and $H_{2L}(1)$. This suggests the following more general result.
\begin{theorem}
The Rogers--Szeg\"o polynomials can be expressed as
\begin{equation}\label{RS}
H_n(t)=\sum_{r=0}^{\lfloor n/2\rfloor} 
t^{2r}(-q/t;q^2)_r(-t;q^2)_{\lfloor (n+1)/2\rfloor -r}
\qbin{\lfloor n/2\rfloor}{r}_{q^2}.
\end{equation}
\end{theorem}
This new representation for the RS polynomials manifests the well-known
facts that
$H_{2n}(-1)=(q;q^2)_n$, $H_{2n+1}(-1)=0$ and
$H_n(-q)=(q;q^2)_{\lfloor (n+1)/2\rfloor}$, which are not immediately 
clear from the standard definition \eqref{RSdef}.

In the following we will give two proofs of \eqref{RS}. In the first we 
show that \eqref{RS} satisfies the recurrence \cite{Andrews76}
\begin{equation}\label{Hrec}
H_{n+1}(t)=(1+t)H_n(t)-(1-q^n)tH_{n-1}(t)
\end{equation}
which determines the RS polynomials together with the
initial conditions $H_0(t)=1$ and $H_1(t)=1+t$.
In the second more complicated and interesting proof, we establish
equality between \eqref{RSdef} and \eqref{RS} using
basic hypergeometric series manipulations.

\begin{proof}[First proof of \eqref{RS}]
Take \eqref{Hrec}, replace $n$ by $2n$ and substitute \eqref{RS}.
Then use $(-t;q^2)_{n-r+1}=(-t;q^2)_{n-r}(1+t q^{2n-2r})$ on the
left and
\begin{equation*}
(1-q^{2n})\qbin{n-1}{r}_{q^2}=(1-q^{2n-2r})\qbin{n}{r}_{q^2}
\end{equation*}
on the right. All terms now pairwise cancel.

Next take \eqref{Hrec}, replace $n$ by $2n-1$ and substitute \eqref{RS}.
Then use 
\begin{equation}\label{qbinrec}
\qbin{n}{r}_{q^2}=\qbin{n-1}{r}_{q^2}+q^{2n-2r}\qbin{n-1}{r-1}_{q^2}
\end{equation}
and cancel one of the two terms on the left with one of the
terms on the right. To proceed use
$(-t;q^2)_{n-r}=(-t;q^2)_{n-r-1}(1+t q^{2n-2r-2})$ on the right
and replace $r$ by $r+1$ on the left.
All resulting terms again pairwise cancel.

The final checking of the initial conditions $H_0(t)=1$ and $H_1(t)=1+t$
is left to the reader.
\end{proof}

\begin{proof}[Second proof of \eqref{RS}]
As a first step we twice use the $q$-binomial theorem \eqref{qbthm}
to expand the $q$-shifted factorials on the right. After this one can extract
coefficients of $t^j$ on both sides leading to the double sum
\begin{multline*}
\sum_{r=0}^{\lfloor n/2\rfloor} \sum_{k=0}^{\lfloor (n+1)/2\rfloor -r}
q^{(2r+k-j)^2+k(k-1)} \\ \times
\qbin{r}{2r+k-j}_{q^2} \qbin{\lfloor (n+1)/2\rfloor -r}{k}_{q^2}
\qbin{\lfloor n/2\rfloor}{r}_{q^2}=\qbin{n}{j}.
\end{multline*}
Now introduce two new summation variables $r'$ and $k'$ by
$r'=2r+k-j$ and $k'=j-k-r$. Eliminating $k$ and $r$ in favour of their
primed counterparts and then dropping the primes yields
\begin{multline*}
\sum_{k=0}^{\lfloor j/2\rfloor}
q^{(j-2k)(j-2k-1)}\qbin{\lfloor (n+1)/2\rfloor-k}{j-2k}_{q^2}
\qbin{\lfloor n/2\rfloor}{k}_{q^2} \\ \times
{_2\phi_2}\biggl[\genfrac{}{}{0pt}{}
{q^{-2(\lfloor n/2\rfloor-k)},q^{-2(j-2k)}}
{q^{-2(\lfloor (n+1)/2\rfloor-k)},0};q^2,q^{3-2\sigma}\biggr]=
\qbin{n}{j},
\end{multline*}
where $\sigma\in\{0,1\}$ is fixed by $n+\sigma\equiv 0\pmod{2}$.
We note that the lower bound on $k$ may be optimized
to $\max(0,j-\lfloor (n+1)/2\rfloor)$.

When $n$ is even the ${_2\phi_2}$ series becomes a ${_1\phi_1}$ which 
sums to $q^{-\binom{j-2k}{2}}(-q;q)_{j-2k}$ by the $a\to\infty$ limit of
\eqref{phi21sumb}
or by \eqref{F34} with $q^2\to 1/q$ and $j\to L-2j$.
When $n$ is odd the ${_2\phi_2}$ sums to $$q^{-\binom{j-2k}{2}}
\frac{1-q^{n-j+1}}{1-q^{n-2k+1}}(-q;q)_{j-2k}$$
by the $a\to\infty$ limit of
\begin{equation}\label{phi32odd}
{_3\phi_2}(a,bq^2,q^{-2n};b,q^{2-2n}/a;q^2,q/a)=
q^{-n}\frac{1-bq^n}{1-b}\frac{(-q,a;q)_n}{(a;q^2)_n}.
\end{equation}
For even $n$ we are thus left with (after replacing $n$ by $2n$)
\begin{equation}\label{even}
\sum_{k=0}^{\lfloor j/2\rfloor} q^{\binom{j-2k}{2}}(-q;q)_{j-2k}
\qbin{n-k}{j-2k}_{q^2} \qbin{n}{k}_{q^2}=\qbin{2n}{j}
\end{equation}
and for odd $n$ (after replacing $n$ by $2n+1$) with
\begin{equation*}
\sum_{k=0}^{\lfloor j/2\rfloor} q^{\binom{j-2k}{2}}(-q;q)_{j-2k}
\frac{1-q^{2n-j+2}}{1-q^{2n-2k+2}}\qbin{n-k+1}{j-2k}_{q^2}
\qbin{n}{k}_{q^2}=\qbin{2n+1}{j}.
\end{equation*}
Multiplying both sides by $(1-q^{2n+2})/(1-q^{2n-j+2})$ this can easily 
be seen to correspond to \eqref{even} with $n\to n+1$.
Hence we only need to consider \eqref{even}. But this is
nothing but \eqref{t2dual} with $(L,j,r)\to (n,n-j,n-j+2k)$.

It remains to prove \eqref{phi32odd}, which for $b=a$ reduces to
\eqref{phi21sumb} and for $b=0$ yields a companion thereof.
Now by the contiguous relation \cite[Eq. (3.2)]{Krattenthaler93}
\begin{multline*}
{_{r+1}\phi_r}\biggl[\genfrac{}{}{0pt}{}{a,bq,(A)}{(B)};q,z\biggr]
={_{r+1}\phi_r}\biggl[\genfrac{}{}{0pt}{}{aq,b,(A)}{(B)};q,z\biggr] \\
+z(b-a)\frac{\prod_{i=1}^{r-1}(1-A_i)}{\prod_{i=1}^r(1-B_i)}
{_{r+1}\phi_r}\biggl[\genfrac{}{}{0pt}{}{aq,bq,(Aq)}{(Bq)};q,z\biggr]
\end{multline*}
with $(a,b,(A),(B),z,q)\to (a,b,(q^{-2n}),(b,q^{2-2n}/a),q/a,q^2)$
the ${_3\phi_2}$ series on the left-hand side of \eqref{phi32odd}
splits into two ${_2\phi_1}$ series, both of which
are summable by \eqref{phi21sumb}. After some simplifications this
yields the claimed right-hand side.
\end{proof}

\section{The generalized Borwein conjecture}\label{sec9}

\subsection{The Borwein conjecture and Bressoud's generalization}
Some years ago Peter Borwein conjectured \cite{Andrews95} that the polynomials
$A_n(q)$, $B_n(q)$ and $C_n(q)$ given by
\begin{equation*}
A_n(q^3)-q B_n(q^3)-q^2 C_n(q^3)=(q,q^2;q^3)_n
\end{equation*}
have nonnegative coefficients.
Defining 
\begin{multline*}
G(N,M;\alpha,\beta,K;q)=
G(N,M;\alpha,\beta,K) \\
=\sum_{j=-\infty}^{\infty}(-1)^j q^{Kj((\alpha+\beta)j+\alpha-\beta)/2}
\qbin{M+N}{N-Kj}
\end{multline*}
it follows from \eqref{FJTP} that \cite{Andrews95}
\begin{align*}
A_n(q)&=G(n,n;4/3,5/3,3) \\
B_n(q)&=G(n+1,n-1;2/3,7/3,3) \\
C_n(q)&=G(n+1,n-1;1/3,8/3,3).
\end{align*}
This led Bressoud \cite{Bressoud96} to a more general conjecture concerning 
the nonnegativity of the coefficients of $G(N,M;\alpha,\beta,K;q)$.
Since we will only be concerned with the case $N=M$ we write
$G(N;\alpha,\beta,K;q)$ instead of $G(N,N;\alpha,\beta,K;q)$,
in accordance with \eqref{G}. We also write $P(q)\geq 0$ for $P(q)$ a 
polynomial in $q$ with nonnegative coefficients.
Then the $N=M$ case of Bressoud's generalized Borwein conjecture can
be stated as follows.
\begin{conjecture}\label{BC}
Let $K,L,\alpha K,\beta K$ be integers such that 
$0\leq \alpha\leq K$, $0\leq \beta\leq K$ and
$1\leq \alpha+\beta\leq 2K-1$. 
Then $G(L;\alpha,\beta,K;q)\geq 0$.
\end{conjecture}
When $\alpha$ and $\beta$ are integers the conjecture becomes
\cite[Thm. 1]{ABBBFV87} of Andrews \textit{et al.}
For fractional values of $\alpha$ and or $\beta$ several cases
of the conjecture have been proven in \cite{Bressoud96,IKS99,W01,W03}.

Without loss of generality one may in fact put stronger restrictions 
on the parameters $\alpha$ and $\beta$ in Conjecture~\ref{BC}. By
\begin{equation*}
G(L;\alpha,\beta,K)=G(L;\beta,\alpha,K)
\end{equation*}
we may assume that 
$0\leq \alpha\leq \beta\leq K$. Furthermore, by
\begin{equation*}
G(L;\alpha,\beta,K;1/q)=q^{-L^2} G(L;K-\beta,K-\alpha,K;q)
\end{equation*}
we may also assume that $1\leq \alpha+\beta\leq K$.
Hence we obtain $\max(0,1-\beta)\leq\alpha\leq\min(\beta,K-\beta)$.

Because of the positivity preserving nature of the $q$-binomial 
transformations of Section~\ref{sec2} we can easily prove many
cases of Conjecture~\ref{BC}. For example, iterating \eqref{seed}
using \eqref{t2b} yields the polynomial identity \eqref{Gplus},
which implies that 
\begin{equation}\label{Gplus2}
G(L;\tfrac{1}{2}(2^k-1),\tfrac{1}{2}(2^k+1),2^k)\geq 0.
\end{equation}
For $k=1$ the polynomial identity \eqref{Gplus} and the according
nonnegativity of the coefficients of $G(L;1/2,3/2,2)$ are due to 
Ismail \textit{et al.} \cite[Prop. 2 (3)]{IKS99}.

A reformulation of the positivity preserving transformations
of Section~\ref{sec2} in the language of the generalized Borwein conjecture
reads as follows.
\begin{lemma}\label{lemGtoG}
If $G(L;\alpha,\beta,K)\geq 0$ then $G(L;\alpha',\beta',K')\geq 0$
with
\begin{subequations}
\begin{align}\label{GtoG1}
\alpha'&=\alpha+K,& \beta'&=\beta+K, & K'&=2K, \\
\alpha'&=\alpha,& \beta'&=\beta, & K'&=2K,
\label{GtoG2} \\
\alpha'&=\alpha+K/2,& \beta'&=\beta+K/2, & K'&=2K,
\label{GtoG3} \\
\alpha'&=2\alpha,& \beta'&=2\beta, & K'&=2K,
\label{GtoG4} \\
\alpha'&=\alpha+K,& \beta'&=\beta+K, & K'&=3K.
\label{GtoG5}
\end{align}
\end{subequations}
\end{lemma}

\begin{proof}
We will only prove \eqref{GtoG1}. All other cases proceed along similar
lines, be it that instead of \eqref{t2} one needs to employ
\eqref{t2dual}, \eqref{t2b}, \eqref{t4} and \eqref{t3}.
Adopting the notation of Section~\ref{sec3} and
writing $F^{\eqref{t2}}_{L,r}(q)$ for $q^{r^2/2}f_{L,r}(q)$, with $f_{L,r}(q)$
given by \eqref{fLr}, we have by Lemma~\ref{lem1} and the assumption that
$G(L;\alpha,\beta,K;q)\geq 0$,
\begin{align*}
0&\leq \sum_{r=0}^{\infty}F^{\eqref{t2}}_{L,2r}(q)
G(r;\alpha,\beta,K;q^2) \\
&=\sum_{j=-\infty}^{\infty} (-1)^j
q^{Kj((\alpha+\beta)j+\alpha-\beta)}
\sum_{r=0}^{\infty}F^{\eqref{t2}}_{L,2r}(q)\qbin{2r}{r-Kj}_{q^2} \\
&=\sum_{j=-\infty}^{\infty} (-1)^j
q^{Kj((2K+\alpha+\beta)j+\alpha-\beta)}\qbin{2L}{L-2Kj} \\
&=G(L;K+\alpha,K+\beta,2K;q).   \qedhere
\end{align*}
\end{proof}

Lemma~\ref{lemGtoG} may be complemented by two more results.
First we remark that the equation following Theorem~2.5 of \cite{W03} 
can be recast as the following lemma.
\begin{lemma}\label{lemW}
For $L$ and $j$ integers there holds
\begin{equation*}
\sum_{\substack{r=0\\r\equiv j\;(2)}}^L
q^{\frac{1}{2}r^2} f_{L,r}(q)\qbin{r}{\frac{1}{2}(r-j)}
=q^{\frac{1}{2}j^2}\qbin{2L}{L-j}
\end{equation*}
with
\begin{equation*}
f_{L,r}(q)=\qbin{L}{r}\sum_{n=0}^{L-r}q^{n(n+r)}\qbin{L-r}{n}.
\end{equation*}
\end{lemma}
By the substitution $q\to 1/q$ this implies a related result.
\begin{corollary}
For $L$ and $j$ integers there holds
\begin{equation*}
\sum_{\substack{r=0\\r\equiv j\;(2)}}^L
q^{\frac{1}{4}r^2} f_{L,r}(q)\qbin{r}{\frac{1}{2}(r-j)}
=q^{\frac{1}{4}j^2}\qbin{2L}{L-j}
\end{equation*}
with
\begin{equation*}
f_{L,r}(q)=\qbin{L}{r}\sum_{n=0}^{L-r} q^{Ln}\qbin{L-r}{n}.
\end{equation*}
\end{corollary}
Note that these results correspond to \eqref{sym2} with $k=1$ and $\gamma=2$
or $\gamma=1$, but that, unlike the solutions to \eqref{sym2} presented
in Section~\ref{sec2}, $f_{L,r}(q)$ is non-factorizable.
{}From Lemma~\ref{lemW} and its corollary we get \cite[Lemma 6.7]{W03}.
\begin{lemma}\label{lemGtoG2}
If $G(L;\alpha,\beta,K)\geq 0$ then $G(L;\alpha',\beta',K')\geq 0$
with
\begin{subequations}
\begin{align}\label{GtoG6}
\alpha'&=\alpha/2+K,& \beta'&=\beta/2+K, & K'&=2K, \\
\alpha'&=(\alpha+K)/2,& \beta'&=(\beta+K)/2, & K'&=2K.
\label{GtoG7}
\end{align}
\end{subequations}
\end{lemma}

The results of Lemmas~\ref{lemGtoG} and \ref{lemGtoG2} can be iterated to
yield a tree of conditional nonnegativity results, subject to various
internal relations. For example, applying \eqref{GtoG1} and then \eqref{GtoG2}
is equivalent to applying \eqref{GtoG2} and then \eqref{GtoG3} (this fact
corresponds to \eqref{rel11a} with $r\to 2r$), but is also equivalent to
application of \eqref{GtoG4} followed by \eqref{GtoG7}. Indeed in each case,
starting with $G(L;\alpha,\beta,K)$ one obtains $G(L;\alpha+K,\beta+K,4K)$.
Even when ignoring these degeneracies it is very complicated to give
a complete description of an arbitrary node of the tree, and for the example 
of the subtree generated by Lemma~\ref{lemGtoG2}, which requires the
theory of continued fractions, we refer to \cite[Prop. 6.8]{W03}.
Instead of trying to achieve maximum generality we take the easy way out
and restrict ourselves to several easy to state and prove examples. 
\begin{proposition}
For $K,\alpha,\beta$ and $k$ integers such that
$\max(0,1-\beta)\leq\alpha\leq\min(\beta,K-\beta)$ and $k\geq 0$,
\begin{equation*}
G(L;\alpha+\tfrac{1}{2}(2^k-1)K,\beta+\tfrac{1}{2}(2^k-1)K,2^kK)\geq 0.
\end{equation*}
\end{proposition}
For $\alpha=0$, $\beta=1$ and $K=1$ this yields \eqref{Gplus2}.

\begin{proof}
The proposition is true for $k=0$ by the remark following Conjecture~\ref{BC}.
The rest follows from \eqref{GtoG3} and induction.
\end{proof}

\begin{proposition}\label{propGG}
For $K,\alpha,\beta$ and $k$ integers such that
$\max(0,1-\beta)\leq\alpha\leq\min(\beta,K-\beta)$ there
holds $G(L;\alpha',\beta',K')\geq 0$ with
\begin{equation*}
\alpha'=2^{-k}(\alpha+\tfrac{5}{14}(8^k-1)K),\quad
\beta'=2^{-k}(\beta+\tfrac{5}{14}(8^k-1)K),\quad
K'=4^k K
\end{equation*}
for $k\geq 0$, and
\begin{multline*}
\alpha'=2^{1-k}(\alpha+\tfrac{1}{28}(3\cdot 8^k-10)K),\quad
\beta'=2^{1-k}(\beta+\tfrac{1}{28}(3\cdot 8^k-10)K), \\
K'=\tfrac{1}{2}\cdot 4^k K
\end{multline*}
for $k\geq 1$.
\end{proposition}

\begin{proof}
The first equation for $k=0$ is true by the remark following 
Conjecture~\ref{BC}. 
Thanks to $\tfrac{5}{14}(8^k-1)+2^{3k-1}=\tfrac{1}{28}(3\cdot 8^k-10)$,
applying \eqref{GtoG3} to $G(L;\alpha',\beta',K')$
with $(\alpha',\beta',K')$ given by the first equation
gives $G(L;\alpha',\beta',K')$ with $(\alpha',\beta',K')$ given by the 
second equation where $k\to k+1$.
Thanks to $\tfrac{1}{28}(3\cdot 8^k-10)+2^{3k-2}=\tfrac{5}{14}(8^k-1)$,
applying \eqref{GtoG7} to $G(L;\alpha',\beta',K')$
with $(\alpha',\beta',K')$ given by the second equation
gives $G(L;\alpha',\beta',K')$ with $(\alpha',\beta',K')$ given by the 
first equation. These observations suffice to conclude the proposition
by induction.
\end{proof}

Reversing the order of \eqref{GtoG3} and \eqref{GtoG7} in the above leads
to the following modification of Proposition~\ref{propGG}
\begin{proposition}
For $K,\alpha,\beta$ and $k$ integers such that
$\max(0,1-\beta)\leq\alpha\leq\min(\beta,K-\beta)$ there
holds $G(L;\alpha',\beta',K')\geq 0$ with 
\begin{equation*}
\alpha'=2^{-k}(\alpha+\tfrac{3}{7}(8^k-1)K),\quad
\beta'=2^{-k}(\beta+\tfrac{3}{7}(8^k-1)K),\quad
K'=4^k K
\end{equation*}
for $k\geq 0$, and
\begin{multline*}
\alpha'=2^{-k}(\alpha+\tfrac{1}{28}(5\cdot 8^k-12)K),\quad
\beta'=2^{-k}(\beta+\tfrac{1}{28}(5\cdot 8^k-12)K), \\
K'=\tfrac{1}{2}\cdot 4^k K
\end{multline*}
for $k\geq 1$.
\end{proposition}
Since we went through quite a bit of trouble to show that the 
coefficients of the polynomial in \eqref{fcubic} are positive we should
at least include one example that makes use of \eqref{GtoG5}.
The next result is obtained by replacing \eqref{GtoG3} in
the proof of Proposition~\ref{propGG} by \eqref{GtoG5}.
\begin{proposition}
For $K,\alpha,\beta$ and $k$ integers such that
$\max(0,1-\beta)\leq\alpha\leq\min(\beta,K-\beta)$ there
holds $G(L;\alpha',\beta',K')\geq 0$ with
\begin{equation*}
\alpha'=2^{-k}(\alpha+\tfrac{4}{11}(12^k-1)K),\quad
\beta'=2^{-k}(\beta+\tfrac{4}{11}(12^k-1)K),\quad
K'=6^k K
\end{equation*}
for $k\geq 0$, and
\begin{multline*}
\alpha'=2^{1-k}(\alpha+\tfrac{1}{44}(5\cdot 12^k-16)K),\quad
\beta'=2^{1-k}(\beta+\tfrac{1}{44}(5\cdot 12^k-16)K), \\
K'=\tfrac{1}{2}\cdot 6^k K
\end{multline*}
for $k\geq 1$.
\end{proposition}
We leave it to the reader to derive more examples of the above kind.

Our next examples use the result \cite[Cor. 3.2]{W01}.
\begin{theorem}\label{thmW}
$G(L;b-1/a,b,a)\geq 0$ for $a,b$ coprime integers such that $0<b<a$.
\end{theorem}

First we note that \eqref{GtoG1} and \eqref{GtoG2}
nicely combine to the following statement.
\begin{lemma}\label{lemGtoG3}
Let $k,i$ be integers such that $0\leq i<2^k$.
Then $G(L;\alpha,\beta,K)\geq 0$ implies that
$G(L;\alpha+iK,\beta+iK,2^k K)\geq 0$.
\end{lemma}
\begin{proof}
For $k=0$ the lemma is trivially true and for $k=1$ it corresponds to
\eqref{GtoG1} when $i=1$ and to \eqref{GtoG2} when $i=0$.
Induction now does the rest since application of \eqref{GtoG1} to $(i,k)$
yields $(k',i')$ with $k'=k+1$ and $2^k\leq i'<2^{k+1}=2^{k'}$ and
application of \eqref{GtoG2} to $(i,k)$ yields 
$(k',i')$ with $k'=k+1$ and $0\leq i'<2^k$. When combined this results in
$(i',k')$ with $k'=k+1$ and $0\leq i'<2^{k'}$.
\end{proof}
In the same vein one can also show that \eqref{GtoG2} and \eqref{GtoG4}
combine in a simple manner.
\begin{lemma}\label{lemGtoG4}
Let $k,i$ be integers such that $0\leq i\leq k$.
Then $G(L;\alpha,\beta,K)\geq 0$ implies that
$G(L;2^i \alpha,2^i \beta,2^k K)\geq 0$.
\end{lemma}
If we now invoke Theorem~\ref{thmW} we obtain the following two
Theorems.
\begin{theorem}\label{thmB1}
$G(L;b-1/a,b,2^k a)\geq 0$ for $k$ a nonnegative integer and
$a,b$ coprime integers such that $0<b<2^ka$.
\end{theorem}
\begin{proof}
Obviously, by Lemma~\ref{lemGtoG3} and Theorem~\ref{thmW}, 
$G(L;b+ia-1/a,b+ia,2^k a)\geq 0$ for $a,b$ coprime integers such that 
$0<b<a$ and $k,i$ integers such that $0\leq i<2^k$. Replacing $b+ia$ by $b$
this becomes
$G(L;b-1/a,b,2^k a)\geq 0$ for $a,b$ coprime integers such that 
$ia<b<(i+1)a$ and $k,i$ integers such that $0\leq i<2^k$.
But since $a$ and $b$ are coprime the condition
$ia<b<(i+1)a$ with $0\leq i<2^k$ may be replaced by
$0<b<2^ka$.
\end{proof}
\begin{theorem}
$G(L;2^i(b-1/a),2^i b,2^{k+k'}a)\geq 0$ for $k$ a nonnegative integer,
$a,b$ coprime integers such that $0<b<2^ka$ and
$k',i$ integers such that $0\leq i<2^{k'}$.
\end{theorem}
Note that this contains the previous theorem as special case.
\begin{proof}
Simply apply Lemma~\ref{lemGtoG4} with $k\to k'$ to Theorem~\ref{thmB1}.
\end{proof}
Our final example arises by repeating the proof of Proposition~\ref{propGG}
but with Theorem~\ref{thmB1} as seed.
\begin{theorem}
Let $k'$ a nonnegative integer and
$a,b$ a pair of coprime integers such that $0<b<2^{k'}a$. Then
$G(L,\alpha,\beta,K)\geq 0$ with
\begin{equation*}
\alpha=2^{-k}(b-1/a+\tfrac{5}{14}(8^k-1)a),\quad
\beta=2^{-k}(b+\tfrac{5}{14}(8^k-1)a),\quad
K=2^{2k+k'}a
\end{equation*}
for $k\geq 0$, and
\begin{multline*}
\alpha=2^{1-k}(b-1/a+\tfrac{1}{28}(3\cdot 8^k-10)a),\quad
\beta=2^{1-k}(b+\tfrac{1}{28}(3\cdot 8^k-10)a), \\
K=2^{2k+k'-1}K
\end{multline*}
for $k\geq 1$.
\end{theorem}

\subsection{New representations of the Borwein polynomials}
Unfortunately the positivity preserving transformations of this paper are 
inadequate for proving the original Borwein conjecture.
Indeed, the only way to for example obtain $A_L(q)$ is by applying the
cubic transformation \eqref{t3} to \eqref{FJTP}, as was done
in the proof of Proposition~\ref{prop1}. But the resulting
\begin{equation*}
A_L(q)=(1-q^{2L})\sum_{n=0}^{\lfloor L/3\rfloor} \frac{q^{3n^2}(q;q)_{3n}
(q^3;q^3)_{L-n-1}}{(q^3;q^3)_{2n}(q^3;q^3)_n(q;q)_{L-3n}}
\end{equation*}
which was first found by Andrews \cite[Eq. (4.5)]{Andrews95}
is insufficient for proving that $A_L(q)\geq 0$.
Another representation follows from application of \eqref{qtrafo5} to
\eqref{J1}, \cite[Thm 2]{Stanton01},
but since \eqref{qtrafo5} is not positivity preserving this
again fails to prove that $A_L(q)\geq 0$.
To conclude we prove alternative representations
for $A_L(q)$ and $C_L(q)$ as triple sums, and use this to formulate
a refinement of the Borwein conjecture.
\begin{theorem}\label{thmAn}
Let $L$ be a nonnegative integer and 
$N_1=n_1+n_2+n_3$, $N_2=n_2+n_3$, $N_3=n_3$.
Then
\begin{equation*}
A_L(q)=\sum_{\substack{n_1,n_2,n_3\geq 0 \\[0.5mm] N_1+N_2+N_3\leq L}}
\frac{q^{N_1^2+N_2^2+N_3^2}(q;q)_{L-N_1}(q;q)_{2L-N_1-N_2}}
{(q;q)_{2L-2N_1}(q;q)_{L-N_1-N_2-N_3}
(q;q)_{n_1}(q;q)_{n_2}(q;q)_{n_3}}.
\end{equation*}
\end{theorem}

\begin{proof}
Define
\begin{equation}\label{Bdef}
\B(L,M,a,b)=\qbin{L+M+a-b}{L+a}\qbin{L+M-a+b}{L-a}.
\end{equation}
Then, according to \cite[Eq. (5.33)]{ASW99}, the following doubly-bounded
analogues of the Rogers--Ramanujan identities hold:
\begin{equation}\label{ASW}
\sum_{j=-\infty}^{\infty}
(-1)^j q^{j(5j+2\sigma+1)/2}\B(L,M,j,j)
=\sum_{n\geq 0} \frac{q^{n(n+\sigma)}(q;q)_{L+M}}
{(q;q)_{L-n}(q;q)_{M-n}(q;q)_n},
\end{equation}
where $\sigma\in\{0,1\}$.
Now for $L,M,a,b$ integers such that not $-L+a\leq -b\leq L+a<b\leq M$
or $-L-a\leq b\leq L-a<-b\leq M$, there holds
\begin{equation}\label{Btrafo}
\sum_{i=0}^M q^{i^2}\qbin{2L+M-i}{2L}\B(L-i,i,a,b)=q^{b^2}\B(L,M,a+b,b).
\end{equation}
This result is known as the Burge transform \cite{Burge93,FLW97,SW00}
and can be applied to \eqref{ASW}.
First let us show that the conditions on the parameters 
(for their origin see \cite{SW00}) are harmless.
{}From $-L+a\leq -b\leq L+a<b\leq M$ one can extract the three conditions
\begin{equation*}
\text{(i)}\;\, L\geq 0, \qquad \text{(ii)}\;\, b>0, \qquad
\text{(iii)}\;\, L+a-b<0
\end{equation*}
and from $-L+a\leq -b\leq L+a<b\leq M$ it follows that
\begin{equation*}
\text{(iv)}\;\, L\geq 0, \qquad \text{(v)}\;\, b<0, \qquad
\text{(vi)}\;\, L-a+b<0.
\end{equation*}
Now in order to transform \eqref{ASW} by the Burge transform
we need to take $a=b=j$.
Since the inequalities (i) and (iii), and also (iv) and (vi) 
become mutually exclusive we can indeed utilize \eqref{Btrafo} to get
\begin{multline*}
\sum_{j=-\infty}^{\infty}
(-1)^j q^{j(7j+2\sigma+1)/2}\B(L,M,2j,j) \\
=\sum_{N_1,N_2\geq 0} \frac{q^{N_1^2+N_2^2+\sigma N_2}(q;q)_L}
{(q;q)_{L-N_1-N_2}(q;q)_{N_1-N_2}(q;q)_{N_2}}\qbin{2L+M-N_1}{2L}.
\end{multline*}
Here $i$ and $n$ have been replaced by $N_1$ and $N_2$, respectively.
We need to apply \eqref{Btrafo} one more time. Since now $a=2j$ and $b=j$
the inequalities (i)--(iii) become $L\geq 0$, $j\geq 0$, and
$L+j<0$ which cannot occur. Similarly, the inequalities (iv)--(v)
are now $L\geq 0$, $j<0$ and $L-j<0$ which again is impossible.
As a result we get
\begin{multline}\label{ALM}
\sum_{j=-\infty}^{\infty}
(-1)^j q^{j(9j+2\sigma+1)/2}\B(L,M,3j,j)= \\
\sum_{\substack{n_1,n_2,n_3\geq 0 \\[0.5mm] N_1+N_2+N_3\leq L}}
\frac{q^{N_1^2+N_2^2+N_3^2+\sigma N_3}(q;q)_{L-N_1}}
{(q;q)_{L-N_1-N_2-N_3}(q;q)_{n_2}(q;q)_{n_3}}\qbin{2L+M-N_1}{2L}
\qbin{2L-N_1-N_2}{2L-2N_1},
\end{multline}
where we have replaced $(i,N_1,N_2)\to (N_1,N_2,N_3)$ and have used
$n_1=N_1-N_2$, $n_2=N_2-N_3$ and $n_3=N_3$.
Letting $M$ tend to infinity, the above identity simplifies to
the theorem when $\sigma=0$ and to an identity for 
$G(L;1,2,3)=(1+q^L)(q^3;q^3)_{L-1}/(q;q)_{L-1}$ when $\sigma=1$.
\end{proof}

By modifying the above proof using an asymmetric version of the Burge
transform the following companion to Theorem~\ref{thmAn} can be shown
to hold.
\begin{theorem}\label{thmCn}
Let $L$ be a nonnegative integer and
$N_1=n_1+n_2+n_3$, $N_2=n_2+n_3$, $N_3=n_3$.
Then
\begin{equation*}
C_L(q)=\sum_{\substack{n_1,n_2,n_3\geq 0 \\[0.5mm] N_1+N_2+N_3\leq L-1}}
\frac{q^{N_1^2+N_2^2+N_3^2+N_1+N_2+N_3}(q;q)_{L-N_1-1}(q;q)_{2L-N_1-N_2-1}}
{(q;q)_{2L-2N_1-1}(q;q)_{L-N_1-N_2-N_3-1}
(q;q)_{n_1}(q;q)_{n_2}(q;q)_{n_3}}.
\end{equation*}
\end{theorem}

\begin{proof}
We extend definition \eqref{Bdef} to
\begin{equation*}
\B_{r,s}(L,M,a,b)=\qbin{L+M+a-b}{L+a}\qbin{L+M-a+b+r+s}{L-a+r}
\end{equation*}
and wish to first show that
\begin{equation}\label{ASWvar}
\sum_{j=-\infty}^{\infty}
(-1)^j q^{j(5j+3)/2}\B_{0,1}(L,M,j,j)
=\sum_{n\geq 0} \frac{q^{n(n+1)}(q;q)_{L+M}}
{(q;q)_{L-n}(q;q)_{M-n}(q;q)_n}.
\end{equation}
To do so we subtract \eqref{ASW} with $\sigma=1$ to get
\begin{equation*}
\sum_{j=-\infty}^{\infty}
(-1)^j q^{j(5j+3)/2}\qbin{L+M}{L+j}
\biggl(\qbin{L+M+1}{L-j}-\qbin{L+M}{L-j}\biggr)=0.
\end{equation*}
By \eqref{qbinrec} this becomes
\begin{equation*}
\sum_{j=-\infty}^{\infty}
(-1)^j q^{5j(j+1)/2}\qbin{L+M}{L+j}\qbin{L+M}{L-j-1}=0.
\end{equation*}
Since the expression on the left is negated after replacing the summation
variable $j$ by $-j-1$ this obviously is true.

Now that \eqref{ASWvar} has been established we apply the following
asymmetric version of the Burge transform \cite{Burge93,SW00}.
For $L,M,a,b,r,s$ integers such that not 
$-L+a-r\leq -b \leq L+a<b+s\leq M+s$
or $-L-a\leq b\leq L-a+r<-b-s\leq M$, there holds
\begin{multline}\label{Btrafo2}
\sum_{i=\max(b,-b-s)}^M q^{i(i+b)}
\qbin{2L+M+r-i}{2L+r}\B_{r+s,s}(L-i-s,i,a,b) \\
=q^{b(b+s)}\B_{r,s}(L,M,a+b,b).
\end{multline}
Before we use this transform we replace $L$ by $L+1$ in \eqref{ASWvar}
and use that $\B_{0,1}(L+1,M,j,j)=\B_{2,1}(L,M,j+1,j)$.
Then we apply \eqref{Btrafo2} with $r=s=1$, $a=j+1$, $b=j$.
This then yields
\begin{multline*}
\sum_{j=-\infty}^{\infty}
(-1)^j q^{j(7j+5)/2}\B_{1,1}(L,M,2j+1,j) \\
=\sum_{N_1,N_2\geq 0} \frac{q^{N_1^2+N_2^2+N_1+N_2}(q;q)_L}
{(q;q)_{L-N_1-N_2}(q;q)_{N_1-N_2}(q;q)_{N_2}}\qbin{2L+M-N_1+1}{2L+1}.
\end{multline*}
We still need to check that the conditions imposed on the asymmetric Burge 
transform do hold. The first condition is certainly satisfied
if not simultaneously
\begin{equation*}
\text{(i)}\;\, L\geq -r, \qquad \text{(ii)}\;\, 2b>s, \qquad
\text{(iii)}\;\, L+a-b-s<0
\end{equation*}
and the second condition is satisfied if not simultaneously
\begin{equation*}
\text{(iv)}\;\, L\geq -r, \qquad \text{(v)}\;\, 2b<-s, \qquad
\text{(vi)}\;\, L-a+b+r+s<0.
\end{equation*}
If $r=s=1$, $a=j+1$, $b=j$ this is easily seen to be the case.

Next we choose $r=0$, $s=1$, $a=2j+1$ and $b=j$,
check that neither (i)--(iii) nor (iv)--(vi) are all satisfied,
and apply \eqref{Btrafo2} to find that
\begin{multline}\label{CLM}
\sum_{j=-\infty}^{\infty}
(-1)^j q^{j(9j+7)/2} \B_{0,1}(L,M,3j+1,j) \\
=\sum_{\substack{n_1,n_2,n_3\geq 0 \\[0.5mm] N_1+N_2+N_3\leq L-1}}
\frac{q^{N_1^2+N_2^2+N_3^2+N_1+N_2+N_3}(q;q)_{L-N_1-1}}
{(q;q)_{L-N_1-N_2-N_3-1}(q;q)_{n_2}(q;q)_{n_3}} \\
\times
\qbin{2L+M-N_1}{2L}
\qbin{2L-N_1-N_2-1}{2L-2N_1-1}.
\end{multline}
In the large $M$ limit this implies Theorem~\ref{thmCn}.
\end{proof}
Theorems~\ref{thmAn} and \ref{thmCn} are insufficient
to conclude that $A_n(q)\geq $ and $C_n(q)\geq 0$. 
It does in fact appear that the polynomials $\A_{L,M}(q)$
and $\C_{L,M}(q)$ given by (the right or left-hand side of)
\eqref{ALM} with $a=1$ and \eqref{CLM} have nonnegative
coefficients. 
\begin{conjecture}
For $L,M$ nonnegative integers 
\begin{equation*}
\A_{L,M}(q)=\sum_{j=-\infty}^{\infty}
(-1)^j q^{j(9j+1)/2}\B(L,M,3j,j)\geq 0
\end{equation*}
and
\begin{equation*}
\C_{L,M}(q)=\sum_{j=-\infty}^{\infty}
(-1)^j q^{j(9j+7)/2}\B_{0,1}(L,M,3j+1,j)\geq 0.
\end{equation*}
\end{conjecture}
However, since
\begin{equation*}
A_L(q)=(q;q)_{2L}\lim_{M\to\infty}\A_{L,M}(q) \quad \text{and}
\quad C_L(q)=(q;q)_{2L}\lim_{M\to\infty}\C_{L,M}(q),
\end{equation*}
this does not imply that $A_L(q)\geq 0$ and $C_L(q)\geq 0$.
Furthermore, it is certainly not true that 
$(q;q)_{2L}\A_{L,M}(q)\geq 0$ or $(q;q)_{2L}\C_{L,M}(q)\geq 0$
for general $L$ and $M$. (We actually believe this never to be true for
positive $L$ and finite $M$).
Nevertheless, the Theorems~\ref{thmAn} and \ref{thmCn} do give rise 
to a two-parameter refinement of the Borwein conjecture.
To describe this we need the polynomials
\begin{equation*}
A_{L,m}(q)=\sum_{\substack{n_1,n_2,n_3\geq 0 \\[0.5mm] N_1+N_2+N_3=m}}
\frac{q^{N_1^2+N_2^2+N_3^2}(q;q)_{L-N_1}(q;q)_{2L-N_1-N_2}}
{(q;q)_{L-m}(q;q)_{2L-2N_1}(q;q)_{n_1}(q;q)_{n_2}(q;q)_{n_3}}
\end{equation*}
and
\begin{equation*}
C_{L,m}(q)=\sum_{\substack{n_1,n_2,n_3\geq 0 \\[0.5mm] N_1+N_2+N_3=m}}
\frac{q^{N_1^2+N_2^2+N_3^2+m}(q;q)_{L-N_1-1}(q;q)_{2L-N_1-N_2-1}}
{(q;q)_{L-m-1}(q;q)_{2L-2N_1-1}(q;q)_{n_1}(q;q)_{n_2}(q;q)_{n_3}}.
\end{equation*}
Comparing these definitions with Theorems \ref{thmAn} and 
\ref{thmCn} shows that
\begin{equation*}
A_L(q)=\sum_{m=0}^L A_{L,m}(q) \quad \text{and}\quad
C_L(q)=\sum_{m=0}^{L-1} C_{L,m}(q).
\end{equation*}
\begin{conjecture}
The polynomials $A_{L,m}(q)$ for $0\leq m\leq L$ 
and $C_{L,m}(q)$ for $0\leq m\leq L-1$
have nonnegative coefficients.
\end{conjecture}
Since this conjecture obviously implies that $A_L(q)\geq 0$ and 
$C_L(q)\geq 0$,
we hope it provides a new way of tackling the
Borwein conjecture. Here it should be noted that proving $A_L(q)\geq 0$ and 
$C_L(q)\geq 0$ is sufficient since $B_L(q)=q^{L^2-1}C_L(1/q)$.

\subsection*{Acknowledgements}
We gratefully acknowledge helpful discussions with George Andrews, 
David Bressoud, George Gasper, Christian Krattenthaler,
Richard Stanley and Dennis Stanton.
We also greatly benefited from Christian Krattenthaler's
Mathematica package HYPQ, a tool for carrying out 
computations involving basic hypergeometric series.

\bibliographystyle{amsplain}

\end{document}